# A COBORDISM CATEGORY ATTACHED TO KHOVANOV–ROZANSKY LINK HOMOLOGIES BASED ON OPERADS

GISA SCHÄFER AND YASUYOSHI YONEZAWA

Abstract. We consider colored operads and their actions on categories. As a special example we construct a cobordism category with a colored operad action arising from oriented planar arc diagrams. This is used to construct an invariant of oriented tangle diagrams with values in the homotopy category attached to the cobordism category. Motivated by Bar-Natan's categorification of the Jones polynomial, it categorifies the quantum $sl_n$ quantum invariants and is adapted to the categorification of the $sl_n$ quantum invariants by Khovanov and Rozansky using matrix factorizations. We conjecture to exist the consistency of the cobordism category and to have an explicit functor from the cobordism category to a category of matrix factorizations.

Contents

## 1. Introduction

Khovanov constructed a categorification of the Jones polynomial [Kho00], a (Laurent) polynomial invariant of links arising from the quantum group attached to the Lie algebra $sl_2$. More precisely he constructed a link homology which assigns to any link a complex of graded $\mathbb{Z}$-modules such that its homology class is an invariant of the link and its graded Euler characteristics is the Jones polynomial. The notion *link homology* was coined by Khovanov [Kho06] and means that a classical polynomial link invariant $p(L) \in \mathbb{Z}[q, q^{-1}]$ is categorified in the sense that it is realized as a graded Euler characteristics of some complex.

After his seminal work the same homology theory was constructed in many, very different, contexts using for instance topological cobordisms [BN05], Lie theory [Str05] or symplectic geometry [SS05].

By now, for the $sl_n$ link invariant, there are also several constructions of $sl_n$ link homology which is a categorification of the $sl_n$ link invariant. For instance, we have a construction using matrix factorizations [KR08, Wu14, Yon11], a geometric construction [CK08], a Lie theoretic construction [Sus07, MS09], a diagrammatic construction [Web17], a construction using foams [MSV09] and a Howe duality construction [CKL10, CK18, MY19].

We are interested in defining a cobordism category for the categorification of the $sl_n$ link invariant in the sense of Bar-Natan. In the study of the functoriality of $sl_n$ link homology, such a cobordism category already appears. The functoriality of the $sl_2$ link homology was proved [Cap08, DCW09, Bla10, EST17]. Subsequently, the functoriality was proved in the case $sl_3$ [Cla09] and in the case $sl_n$[MEW18].

Mackaay, Stosic and Vaz define a $\mathbb{Z}$-graded $\mathbb{Q}$-linear cobordism category [MSV09] and Morrison-Nieh define a cobordism category for $n = 3$ using Khovanov's form construction [MN08, Kho04] but our cobordism category is a $\mathbb{Z}\oplus\mathbb{Z}/2\mathbb{Z}$-graded $\mathbb{Z}[\frac{1}{n}]$-linear category which is derived from the structure of matrix factorizations. Therefore, our relations between cobordisms differ from those of Mackaay-Stosic-Vaz by coefficients to better adapt to the matrix factorizations.

Moreover, each set of morphisms in Morrison-Nieh and Mackaay-Stosic-Vaz categories is finite dimensional [MSV09, Proposition 6.2] but sets of morphisms in our cobordism category are not finite dimensional in general (See Lemma 5.17 (Ra) (Rb) (Rc)). Because the set of cobordism relations in our cobordism category could be minimal for giving chain complexes invariance under Reidemeister moves.

Subsequently, Queffelec-Rose define a cobordism category consisting of ladder-type colored diagrams [QR16] and Robert-Wagner define a cobordism category consisting of closed colored diagrams [RW17]. Both these categories induce the colored Khovanov-Rozansky homology and especially there exists a functor from the Khovanov-Lauda categoried quantum group [KL09, KL10] to Queffelec-Rose's category using the categorical skew Howe duality.

In the categorification using the matrix factorizations, we can concern not only ladder-type and closed colored diagrams but also general colored diagrams. Therefore, we need generators of colored saddle cobordisms in the category. In further study, we will construct a cobordism category consisting of general colored diagrams whose subcategory could be related to Queffelec-Rose and Robert-Wagner cobordism categories.

Moreover, our paper is motivated by defining a new concept of a category with a colored operad action and more concretely constructing the cobordism category for $sl_n$ link homology using it. We develop this framework in detail and then apply it to the special case of oriented planar arc diagrams and trivalent diagrams.

Operads (also called multi-categories, see e.g. [EM06]) were introduced in the 1970s for the purpose of studying homotopy invariant algebraic structures on topological spaces (see e.g. [BV73, May72]) and were subsequently extended to colored operads [BM07]. We generalize this notion to twisted colored operads and connect it to category theory. The main tool here is the notion of colored operads acting on categories.

We first introduce a new category of cobordisms $\mathcal{COB}^{gr}_{G_n/L}$. The category $\mathcal{COB}^{gr}_{G_n/L}$ is the quotient category of $\mathcal{COB}^{gr}_{G_n}$ whose morphism sets consist of the morphism sets of $\mathcal{COB}^{gr}_{G_n}$ modulo the local relations $L$ associated to the category of matrix factorizations. Then we construct an invariant of links with values in the homotopy category of complexes in our category of cobordisms.

We conjecture to control the 2-morphisms of the Khovanov-Rozansky's category of matrix factorizations $\mathrm{HMF}^{KR_n}$ using an operad structure completely and to obtain a well-defined functor from $\mathcal{COB}^{gr}_{G_n/L}$ to $\mathrm{HMF}^{KR_n}$. This existence of a functor naturally gives us a consistency of the cobordism category $\mathcal{COB}^{gr}_{G_n/L}$. Moreover, we expect that our cobordism category is universal in the sense that there should be a functor into any other $sl_n$ link homology theory.

Two problem arises here for defining the functor: The first is the sign ambiguity of morphisms addressed in [KR08, Proposition 36]. A similar issue has already appeared in the Khovanov homology [Jac04, Kho05]. We expect that there exists an operad action on $\mathrm{HMF}^{KR_n}$ solving the sign problem. See Remark 8.19 and Problem 8.20.

The second problem is the functoriality of the assignment $\mathcal{F}$ from $\mathcal{COB}^{gr}_{G_n/L}$ to $\mathrm{HMF}^{KR_n}$. In other words, we require that the assignment $\mathcal{F}$ has the compatibility with the cobordism local relation of $\mathcal{COB}^{gr}_{G_n/L}$ and the structure of morphisms of matrix factorizations. The cobordism local relations except the cobordism isotopy arises from morphisms of matrix factorizations. We need to show that the assignment $\mathcal{F}$ preserves the cobordism isotopy relations in $\mathcal{COB}^{gr}_{G_n/L}$.

In a subsequent study, we should decide what the generators of the isotopy relations are. For instance, the following cobordisms are composed of different cobordisms but these are isotopic.

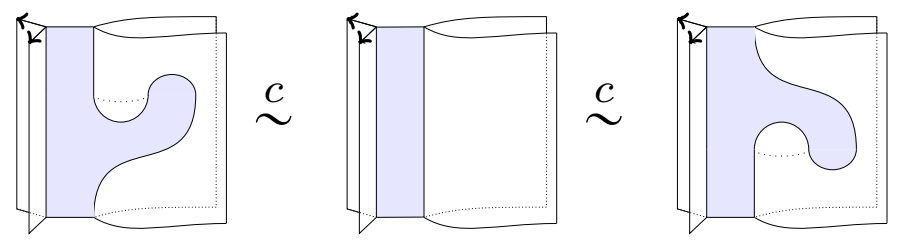

These problems give us a difficulty to construct a functor. Note that the existence of the functor implies the consistency of the relations in our cobordism category, i.e. the category is not trivial. We believe that incorporating all isotopies the theorem still holds, but at the moment we are not able to give a generating set for all cobordism isotopies.

We finally discuss a relation between the cobordism category $\mathcal{COB}^{gr}_{G_n/L}$ and the category of matrix factorizations $\mathrm{HMF}^{KR_n}$.

The paper is organized as follows. In Section 2 we introduce the notion of a twisted colored operad C, C-sets (which are sets with an action of C) and C-categories (which means categories with an action of C). In Section 3 we introduce the notion of a C-category generated by elements and morphisms using an operad action. Then Section 4 defines the main player, a twisted colored operad $P$ consisting of oriented planar arc diagrams. We show in Theorem 4.6 the following result:

**Theorem.** P *is a braided* $(S_\pm, \iota)$*-colored operad.*

Our category of cobordisms $\mathcal{COB}^{gr}_{G_n/L}$ is then constructed in Section 5 as a P-category. The cobordism local relations $L$ in Definition 5.14 are derived from the morphism structure of matrix factorizations. The morphism of matrix factorizations corresponding to the saddle morphism has the sign problem described in Remark 8.19. Since we define $\mathcal{COB}^{gr}_{G_n/L}$ generated by small building blocks (cobordisms) using an operad action, we expect that this cobordism category has consistency.

**Conjecture** (Consistency conjecture)**.** *The* $\mathcal{COB}^{gr}_{G_n/L}$ *is a non trivial* P*-category.*

In Section 6 we define a complex of planar diagrams of $\mathcal{COB}^{gr}_{G_n/L}$ from an oriented tangle. In Appendix A, we show that this complex is invariant under Reidemeister moves in the homotopy category $\mathrm{K}^{\mathrm{b}}(\mathrm{Mat}(\mathcal{COB}^{gr}_{G_n/L}))$, where $\mathrm{Mat}(\mathcal{COB}^{gr}_{G_n/L})$ is the additive closure of $\mathcal{COB}^{gr}_{G_n/L}$. The homotopy equivalences required for the Reidemeister moves are given explicitly. We believe that this might help to do explicit computations for Khovanov-Rozansky-link homology which is hard in general [CM].

In Section 7, we define matrix factorizations corresponding to oriented arc planar diagrams in P. We find that the set $\mathrm{P}_{\mathrm{HMF}}$ composed of these matrix factorizations naturally has a colored operad structure and there exists an operad morphism from P to the set of matrix factorizations.

**Theorem.** $\mathrm{P}_{\mathrm{HMF}}$ *is a braided* $(S_\pm, \iota)$*-colored operad.*

In Section 8, we discuss a map $\mathcal{F}$ from the cobordism category $\mathcal{COB}^{gr}_{G_n/L}$ to the category of matrix factorizations $\mathrm{HMF}^{KR_n}$ and problems for getting a functor from $\mathcal{COB}^{gr}_{G_n/L}$ to $\mathrm{HMF}^{KR_n}$. We find a $\mathrm{P_{HMF}}$-set structure on the object sets of the category of matrix factorizations $\mathrm{HMF}^{KR_n}$.

**Theorem.** *The collection of object sets* $\{\mathrm{Ob}(\mathrm{HMF}^{KR_n}(s_i))\}$ *is a* $\mathrm{P_{HMF}}$*-set.*

We expect that there exists a $\mathrm{P_{HMF}}$-set structure on the morphism sets of $\mathrm{HMF}^{KR_n}$ which solves the sign problem and there exists a functor from $\mathcal{COB}^{gr}_{G_n/L}$ to $\mathrm{HMF}^{KR_n}$ of P-categories.

**Acknowledgements**: The authors would like to thank Catharina Stroppel for suggesting the construction of a cobordism category using operads and some helpful discussions. The first author would like to thank the GRK 1150 and the Max Planck Institute for Mathematics for their support during the creation of the main part of this paper. The second author would like to thank the Hausdorff Center of Mathematics in Bonn and the Max Planck Institute for Mathematics for their support during his visit, when the main part of this paper was studied.

## 2. Colored operad and category with colored operad action

We introduce the notion of twisted colored operads, a generalization of colored operads in the sense of e.g. [BM07], by adding an automorphism $\varphi : S \to S$ which twists the colors $S$ in the multiplication of the colored operad. We then consider sets, modules, algebras and categories over these operads.

### 2.1. $(S, \varphi)$-colored operad.

Let $S$ be a set with an automorphism $\varphi : S \to S$.

**Definition 2.1.** An $(S, \varphi)$*-colored operad* $\mathrm{C} = (\mathrm{C}, \alpha, I)$ consists of the data

- For each $m$-tuple $\underline{s} = (s_1, ..., s_m) \in S^m$ $(m \geqslant 0)$, called *input*, and $s_0 \in S$, called *output*, a set
$$\mathrm{C}_m(s_0; \underline{s}) = \mathrm{C}(s_0; s_1, ..., s_m), \tag{2.1}$$
- a distinguished element $I_s \in \mathrm{C}(s; \varphi(s))$ for any $s \in S$ where $\mathrm{C}(s; \varphi(s)) \neq \varnothing$,
- an operation $\alpha = \{\alpha_{\underline{s}^0}\}_{\underline{s}^0 \in \bigsqcup_{k \geqslant 0} S^{k+1}}$ (also called a multiplication) composed of
$$\alpha_{\underline{s}^0} : \mathrm{C}(\underline{s}^0) \times \mathrm{C}(\underline{s}^1) \times \cdots \times \mathrm{C}(\underline{s}^k) \longrightarrow \mathrm{C}(\underline{s}) \tag{2.2}$$
for non-empty sets $\mathrm{C}(\underline{s}^i)$, where
$$\begin{aligned} \underline{s}^0 &= (s_0; \varphi(s_0^{(1)}), ..., \varphi(s_0^{(k)})) \in S^{k+1}, \\ \underline{s}^i &= (s_0^{(i)}; s_1^{(i)}, ..., s_{m_i}^{(i)}) \in S^{m_i+1} \quad (i = 1, ..., k), \\ \underline{s} &= (s_0; s_1^{(1)}, ..., s_{m_1}^{(1)}, ..., s_1^{(k)}, ..., s_{m_k}^{(k)}) \in S^{m_1+\cdots+m_k+1}, \end{aligned}$$

such that $\alpha$ is associative and $I = \{I_s | s \in S\}$ is a unit with respect to $\alpha$ in the sense defined below.

**Remark 2.2.** It might be helpful to visualize an element from $\mathrm{C}(s_0; s_1, ..., s_m)$ as a transporter between the input data $(s_1, ..., s_m)$ and the output data $s_0$.

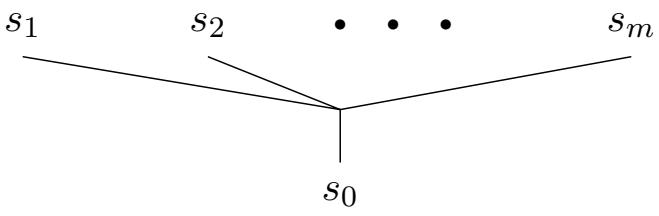

$\in \quad \mathrm{C}(s_0; s_1, ..., s_m)$ and 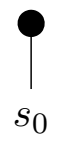

$\in \quad \mathrm{C}(s_0)$

The multiplication, say $\alpha : (f_0; f_1, ..., f_k) \mapsto f$, assigns a new transporter $f$ which first creates from the input data
$$(s_1^{(1)}, ..., s_{m_1}^{(1)}, s_1^{(2)}, ..., s_1^{(k)}, ..., s_{m_k}^{(k)})$$
the output $(s_0^{(1)}, s_0^{(2)}, ..., s_0^{(k)})$ using $f_1$,...,$f_k$, and then takes this output, twisted by $\varphi$, as a new input to create via $f_0$ the final output $s_0$ (the twist map is indicated by the green

symbols):

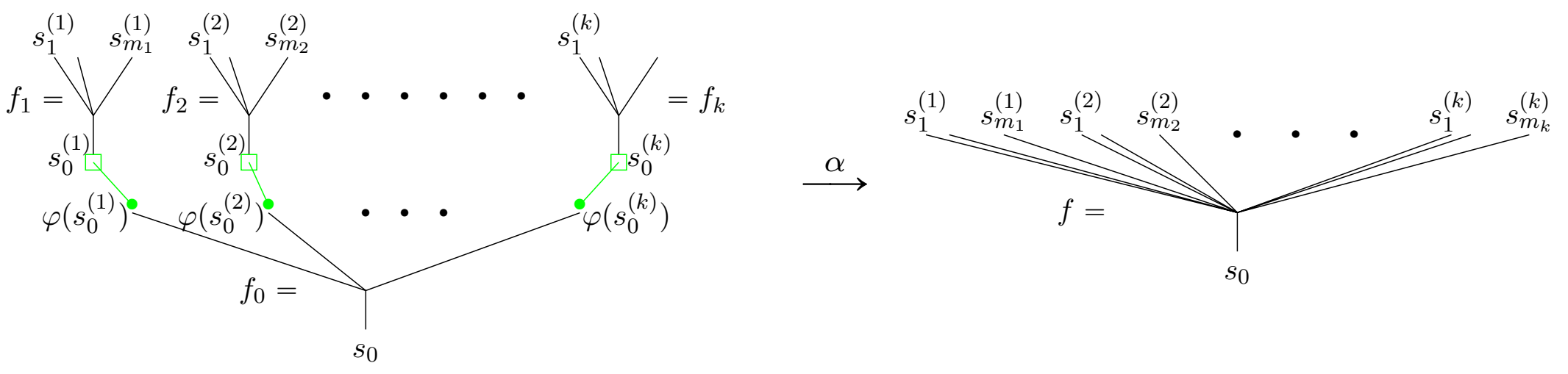

**Definition 2.3.** The multiplications $\alpha$ is *associative* if for non empty sets the following diagram commutes:

$$\begin{array}{ccccc}
\mathrm{C}(\underline{s}) \times \prod_{i=1}^{k} \left( \mathrm{C}(\underline{s}^i) \times \prod_{j=1}^{m_i} \mathrm{C}(\underline{t}^{ij}) \right) & \xrightarrow{\mathrm{id}_{\mathrm{C}(\underline{s})} \times \alpha_{\underline{s}^1} \times \cdots \times \alpha_{\underline{s}^k}} & \mathrm{C}(\underline{s}) \times \prod_{i=1}^{k} \mathrm{C}(\underline{u}^i) & \xrightarrow{\alpha_{\underline{s}}} & \mathrm{C}(\underline{w}) \\
\downarrow \xi \, \cong & & & & \uparrow \alpha_{\underline{v}} \\
\mathrm{C}(\underline{s}) \times \prod_{i=1}^{k} \mathrm{C}(\underline{s}^i) \times \prod_{i=1}^{k} \prod_{j=1}^{m_i} \mathrm{C}(\underline{t}^{ij}) & & \xrightarrow{\alpha_{\underline{s}} \times \mathrm{id}} & & \mathrm{C}(\underline{v}) \times \prod_{i=1}^{k} \prod_{j=1}^{m_i} \mathrm{C}(\underline{t}^{ij})
\end{array}$$

where $\xi$ is the obvious bijection of sets and

$$\begin{aligned}
&\underline{s} := (s_0; \varphi(s_1), ..., \varphi(s_k)) \in S^{k+1}, \quad \underline{s}^i := (s_i; \varphi(t_1^{(i)}), ..., , \varphi(t_{m_i}^{(i)})) \in S^{m_i+1} \quad (i = 1, ..., k), \\
&\underline{t}^{ij} := (t_j^{(i)}; t_1^{(i,j)}, ..., t_{l_{ij}}^{(i,j)}) \in S^{l_{ij}+1} \quad (j = 1, ..., m_i), \\
&\underline{u}^i := (s_i; t_1^{(i,1)}, ..., t_{l_{i1}}^{(i,1)}, ..., t_1^{(i,m_i)}, ..., t_{l_{im_i}}^{(i,m_i)}) \in S^{l_{i1}+\cdots+l_{im_i}+1}, \\
&\underline{v} := (s_0; \varphi(t_1^{(1)}), ..., \varphi(t_{m_1}^{(1)}), ..., \varphi(t_1^{(k)}), ..., \varphi(t_{m_k}^{(k)})) \in S^{m_1+\cdots+m_k+1}, \\
&\underline{w} := (s_0; t_1^{(1,1)}, ..., t_{l_{11}}^{(1,1)}, ..., t_1^{(1,m_1)}, ..., t_{l_{1m_1}}^{(1,m_1)}, ..., t_1^{(l,m_k)}, ..., t_{l_{km_k}}^{(k,m_k)}).
\end{aligned}$$

**Definition 2.4.** The elements $I_s$, $s \in S$, form a *unit* with respect to $\alpha$ if each $I_s$ satisfies that for any $c \in \mathrm{C}(s_0; \varphi(s_1), ..., \varphi(s_k))$

$$\alpha_{(s_0;\varphi(s_0))}(I_{s_0}, c) = c, \quad \alpha_{(s_0;\varphi(s_1),...,\varphi(s_k))}(c, I_{s_1}, ..., I_{s_k}) = c. \tag{2.3}$$

**Remark 2.5.** If $\varphi$ is the identity we get the usual notion of a colored operad [BM07]. In particular, if the coloring set $S$ has only one color, we get the usual notion of an operad with $P(m) = C(1; \underbrace{1, \ldots, 1}_{m})$.

**Example 2.6.** Let $S = \{1\}$ and $V$ a $\Bbbk$-vector space, then $\mathrm{C}(1; \overbrace{1, 1, ..., 1}^{m}) = \mathrm{Hom}_{\Bbbk}(V^{\otimes m}, V)$ with the obvious composition $\alpha$ and $I = \{\mathrm{id}_V \in \mathrm{C}(1; 1)\}$ is an $(S, \mathrm{id})$-colored operad, called the *endomorphism operad*. More generally, we can take $S = \{1, 2, ..., r\}$ and a $\Bbbk$-vector space $V_s$ for each $s \in S$. Then

$$\mathrm{C}(s_0; s_1, ..., s_m) = \mathrm{Hom}_{\Bbbk}(V_{s_1} \otimes \cdots \otimes V_{s_m}, V_{s_0})$$

with unit $I = \{\mathrm{id}_{V_s} | s \in S\}$ is an $(S, \mathrm{id})$-colored operad. If we take $S = \{\pm 1, \pm 2, ..., \pm r\}$ and fix a $\Bbbk$-vector space $V_s$ with an isomorphism $\delta_s : V_{+s} \to V_{-s}$ for each $+s \in S$, then the obvious composition with the identification $\delta_s$ defines an $(S, \varphi)$-colored operad with $\varphi(s) = -s$.

**Example 2.7.** Another standard example of an ( $(\{1\}, \mathrm{id})$-colored) operad is the little interval operad with

$$\mathrm{C}(1;\underbrace{1,...,1}_{m}) = P(m) = \left\{(\lambda_1, \lambda_2, ..., \lambda_m) \,\middle|\, \begin{array}{c} \lambda_i : [0,1] \to [0,1] \text{ linear embeddings such that} \\ \lambda_i([0,1]) \cap \lambda_j([0,1]) = \varnothing \text{ for } i \neq j \end{array} \right\}$$

and the unit $I = \{I_1 = \mathrm{id}_{[0,1]} \in P(1)\}$ with composition

$$\begin{array}{rccc} \alpha : & P(m) \times P(k_1) \times \cdots \times P(k_m) & \longrightarrow & P(k_1 + \cdots + k_m) \\ & \cup\!\!\!\!\!\!\!\! \quad & & \\ & (\lambda, \lambda^{(1)}, ..., \lambda^{(m)}) & \longmapsto & \lambda', \end{array}$$

where

$\lambda = (\lambda_1, ..., \lambda_m), \lambda^{(i)} = (\lambda^{(i)}_1, ..., \lambda^{(i)}_{k_i})$ for $1 \leqslant i \leqslant m, \lambda' = (\lambda_1\lambda^{(1)}_1, ..., \lambda_1\lambda^{(1)}_{m_1}, \lambda_2\lambda^{(2)}_1, ..., \lambda_k\lambda^{(k)}_{m_k})$.

**Example 2.8.** Any $\Bbbk$-algebra $A$ with unit $1_A$ is a $(\{1\}, \mathrm{id})$-colored operad via

$$\mathrm{C}(1;\underbrace{1,...,1}_{m}) = \begin{cases} A & \text{if } m = 1, \\ \varnothing & \text{otherwise,} \end{cases}$$

with $I = \{I_1 = 1_A \in \mathrm{C}(1;1)\}$ and $\alpha : \mathrm{C}(1;1) \times \mathrm{C}(1;1) \to \mathrm{C}(1;1)$ is the multiplication in $A$.

**Example 2.9.** Let $S = \{1, 2, ..., r\}$ and fix $\sigma \in \mathbb{S}_r$.

$$\begin{aligned} \mathrm{C}(s_0; s_1, \ldots, s_m) = \{&\text{planar rooted trees where the root is colored by } s_0 \in S \\ &\text{and the leaves are colored by } s_1, \ldots, s_k \in S\}. \end{aligned}$$

The multiplication is defined by gluing $\sigma(s)$-colored leaves and $s$-colored roots and then forgetting the color there. $I_s \in \mathrm{C}(s; \sigma(s))$ is the line with one endpoint labeled $s$ and the other labeled $\sigma(s)$.

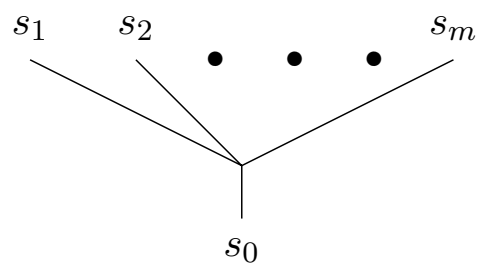

FIGURE 1. An element in $\mathrm{C}(s_0; s_1, \ldots, s_m)$

**Example 2.10.** Any $(S, \varphi)$-colored operad C gives rise to another $(S, \varphi)$-colored operad $\overline{\mathrm{C}}$, called the *identities operad*, with

$$\overline{\mathrm{C}}(s_0; s_1, ..., s_m) := \{\mathrm{id}_c \,|\, c \in \mathrm{C}(s_0; s_1, ..., s_m)\}, \mathrm{id}_c(c') = \begin{cases} c' & \text{if } c' = c \\ 0 & \text{otherwise,} \end{cases}$$

and with the obvious multiplication.

**Example 2.11.** The planar diagram operad defined by Jones [Jon] is an $(S, \iota)$-colored operad, where the colors are sequences of sign

$$S = \{\epsilon(2k) = (\epsilon_1, \epsilon_2, ..., \epsilon_{2k}) | k \geqslant 1, \epsilon_j \in \{+, -\}\} \sqcup \{\varnothing\}$$

and the automorphism $\iota$ defined by $\iota(\varnothing) = \varnothing$, $\iota(+) = -$, $\iota(-) = +$ and $\iota(\epsilon(2k)) = (\iota(\epsilon_1), ..., \iota(\epsilon_{2k}))$.

## 2.2. Morphisms of colored operads.

**Definition 2.12.** A *morphism* between two $(S,\varphi)$-colored operads $\mathrm{C}=(\mathrm{C},\alpha^{\mathrm{C}},I^{\mathrm{C}})$ and $\mathrm{D}=(\mathrm{D},\alpha^{\mathrm{D}},I^{\mathrm{D}})$ is a collection $\Phi$ of maps $\Phi_{(s_0;s_1,\dots s_n)}:\mathrm{C}(s_0;s_1,\dots,s_n)\to\mathrm{D}(s_0;s_1,\dots,s_n)$ which is compatible with

- the multiplication:

$$\Phi_{\underline{s}}\circ\alpha^{\mathrm{C}}_{\underline{s}^0}=\alpha^{\mathrm{D}}_{\underline{s}}\circ(\Phi_{\underline{s}^0}\times\cdots\times\Phi_{\underline{s}^k}):\mathrm{C}(\underline{s}^0)\times\mathrm{C}(\underline{s}^1)\times\cdots\times\mathrm{C}(\underline{s}^k)\to\mathrm{D}(\underline{s})$$

  where

$$\begin{aligned}\underline{s}^0 &= (s_0;\varphi(s_0^{(1)}),...,\varphi(s_0^{(k)}))\in S^{k+1},\\ \underline{s}^i &= (s_0^{(i)};s_1^{(i)},...,s_{m_i}^{(i)})\in S^{m_i+1}\quad(i=1,...,k),\\ \underline{s} &= (s_0;s_1^{(1)},...,s_{m_1}^{(1)},...,s_1^{(k)},...,s_{m_k}^{(k)})\in S^{m_1+\cdots+m_k+1},\end{aligned}$$

- the unit: $\Phi_{(s;\varphi(s))}(I^{\mathrm{C}}_s)=I^{\mathrm{D}}_s$ for all $s\in S$.

Moreover, a collection $\Phi$ is an *isomorphism* if $\Phi_{(s_0;s_1,\dots,s_m)}$ is a bijection for all $m+1$-tuple $(s_0;s_1,\dots,s_m)\in S^{m+1}$. Two $(S,\varphi)$-colored operads C and D are *isomorphic* if there exists an isomorphism between C and D.

**Lemma 2.13.** *Any $(S,\varphi)$-colored operad is isomorphic to the associated identities operad.*

*Proof.* See Example 2.10. For $\underline{s}\in\bigsqcup_{k\geqslant 0}S^{k+1}$ and $c\in\mathrm{C}(\underline{s})$, the assignment $c\overset{\Phi_{\underline{s}}}{\longmapsto}\mathrm{id}_c$ defines the isomorphism. □

## 2.3. Braided colored operad.

**Definition 2.14.** Let $B_m$ be the braid group

$$\langle b_1,...,b_{m-1}|b_ib_{i+1}b_i=b_{i+1}b_ib_{i+1},b_ib_j=b_jb_i(|j-i|>1)\rangle.$$

Let $\mathbb{S}_m$ be the symmetric group

$$\langle s_1,...,s_{m-1}|s_i^2=1,s_is_{i+1}s_i=s_{i+1}s_is_{i+1},s_is_j=s_js_i(|j-i|>1)\rangle.$$

The group homomorphism $B_m\to\mathbb{S}_m$ is defined by the generator $b_i\in B_m$ mapping to $s_i\in\mathbb{S}_m$. We denote the image of $b\in B_m$ by $\hat{b}\in\mathbb{S}_m$.

An $(S,\varphi)$-colored operad $\mathrm{C}=(\mathrm{C},\alpha,I)$ is *braided* if there exists an action of the braid groups $B_m$, $m\geqslant 2$, by invertible automorphisms

$$b_*:\mathrm{C}(s_0;s_1,\dots,s_m)\longrightarrow\mathrm{C}(s_0;s_{\hat{b}(1)},\dots,s_{\hat{b}(m)})\tag{2.4}$$

for $b\in B_m$ satisfying the following compatibility with the multiplication

$$\begin{array}{ccc}\mathrm{C}(s;\varphi(s_1),...,\varphi(s_m))\times\mathrm{C}(\underline{s}^1)\times\cdots\times\mathrm{C}(\underline{s}^m) & \xrightarrow{\ \alpha\ } & \mathrm{C}(\underline{s})\\ \big\downarrow{\scriptstyle b_*\times b^*} & & \big\downarrow{\scriptstyle b(k_1,...,k_m)_*}\\ \mathrm{C}(s;\varphi(s_{\hat{b}(1)}),...,\varphi(s_{\hat{b}(m)}))\times\mathrm{C}(\underline{s}^{\hat{b}(1)})\times\cdots\times\mathrm{C}(\underline{s}^{\hat{b}(m)}) & \xrightarrow{\ \alpha\ } & b(k_1,...,k_m)_*(\mathrm{C}(\underline{s})),\end{array}\tag{2.5}$$

where $\underline{s}^i=(s_i;s_1^{(i)},...,s_{k_i}^{(i)})$ $(i=1,...,m)$, $\underline{s}=(s;s_1^{(1)},...,s_{k_1}^{(1)},...,s_1^{(m)},...,s_{k_m}^{(m)})$, and

$$b^*:\mathrm{C}(\underline{s}^1)\times\cdots\times\mathrm{C}(\underline{s}^m)\to\mathrm{C}(\underline{s}^{\hat{b}(1)})\times\cdots\times\mathrm{C}(\underline{s}^{\hat{b}(m)})$$

is the bijective transformation of sets and $b(k_1,...,k_m)\in B_{k_1+\cdots+k_m}$ is the block braiding associated to $b$.

An $(S,\varphi)$-colored operad $\mathrm{C}=(\mathrm{C},\alpha,I)$ is *symmetric* if C is a braided $(S,\varphi)$-colored operad and invertible automorphisms $b_*$ satisfy $(b_i)_*^2=1_{B_m*}$ $(i=1,...,m-1)$.

Examples 2.6, 2.7 are symmetric operads. The colored operad $P$ in Section 4 is a braided operad.

A collection $\Phi$ is a *morphism of braided colored operads*, if additionally the morphisms $\Phi_{(s_0;s_1,...s_m)}$ between C and D as a colored operad are compatible with the action of braid group:

$$b_* \circ \Phi_{(s_0;s_1,...,s_k)} = \Phi_{(s_0;s_{\hat{b}(1)},...,s_{\hat{b}(k)})} \circ b_* : \mathrm{C}(s_0; s_1, ..., s_k) \to \mathrm{D}(s_0; s_{\hat{b}(1)}, ..., s_{\hat{b}(k)}),$$

where $b \in B_k$.

2.4. **Operad action on a collection of sets with algebraic structure.** We define a collection of $S$-colored sets $\{X(s)\}_{s\in S}$ with an action of an $(S,\varphi)$-colored operad C, simply called C-set. We also define C-sets with several algebraic structures, a C-semigroup, a C-monoid, a C-module and a C-algebra.

**Definition 2.15.** Let $S$ be a set and let $\mathrm{C} = (\mathrm{C}, \alpha, I)$ be an $(S,\varphi)$-colored operad.

A C-*set* $(X, \beta)$ is an $S$-colored set $X = \{X(s)\}_{s\in S}$ together with multiplication maps $\beta_{\underline{s}}$ for $\underline{s} \in S^{m+1}$

$$\beta_{(s_0;\varphi(s_1),...,\varphi(s_m))} : \mathrm{C}(s_0; \varphi(s_1), ..., \varphi(s_m)) \times X(s_1) \times \cdots \times X(s_m) \to X(s_0),$$

satisfying

(C1) the compatibility with the unit: $\beta_{(s;\varphi(s))}(I_s, x_s) = x_s$ for $s \in S$ and $x_s \in X(s)$

(C2) the following compatibility with $\alpha$:

$$\mathrm{C}(\underline{s}) \times \prod_{i=1}^{k}\left(\mathrm{C}(\underline{s}^i) \times \prod_{j=1}^{m_i} X(s_j^{(i)})\right) \xrightarrow{\mathrm{id}_{\mathrm{C}(\underline{s})} \times \beta_{\underline{s}^1} \times \cdots \times \beta_{\underline{s}^k}} \mathrm{C}(\underline{s}) \times \prod_{i=1}^{k} X(s_i) \xrightarrow{\beta_{\underline{s}}} X(s)$$

$$\mathrm{C}(\underline{s}) \times \prod_{i=1}^{k}\left(\mathrm{C}(\underline{s}^i) \times \prod_{j=1}^{m_i} X(s_j^{(i)})\right) \xrightarrow[\cong]{\zeta} \mathrm{C}(\underline{s}) \times \prod_{i=1}^{k} \mathrm{C}(\underline{s}^i) \times \prod_{i=1}^{k}\prod_{j=1}^{m_i} \mathrm{C}(s_j^{(i)})$$

$$\mathrm{C}(\underline{s}) \times \prod_{i=1}^{k} \mathrm{C}(\underline{s}^i) \times \prod_{i=1}^{k}\prod_{j=1}^{m_i} \mathrm{C}(s_j^{(i)}) \xrightarrow{\alpha_{\underline{s}} \times \mathrm{id}} \mathrm{C}(\underline{s'}) \times \prod_{i=1}^{k}\prod_{j=1}^{m_i} \mathrm{C}(s_j^{(i)}) \xrightarrow{\beta_{\underline{s'}}} X(s)$$

where $\zeta$ is the obvious bijection and

$$\underline{s} = (s; \varphi(s_1), ..., \varphi(s_k)) \in S^{k+1}, \underline{s}^i = (s_i; \varphi(s_1^{(i)}), ..., \varphi(s_{m_i}^{(i)})) \in S^{m_i+1} \quad (i = 1, ..., k),$$
$$\underline{s'} = (s; \varphi(s_1^{(1)}), ..., \varphi(s_{m_1}^{(1)}), ..., \varphi(s_1^{(k)}), ..., \varphi(s_{m_k}^{(k)})) \in S^{m_1+\cdots+m_k+1}.$$

For any colored operad C and C-set $X$, the set $\overline{X} := \{\mathrm{id}_x | x \in X\}$ is a C-set whose C-action factors through the operad $\overline{\mathrm{C}}$.

Note that a C-set is usually called C-algebra in the literature. For non-twisted operads it is the same as an operad morphism from C to the endomorphism operad of $X$. For twisted colored operads, the definition of the endomorphism operad is more involved and thus we chose to generalize the other definition.

**Definition 2.16.** A *semigroup over* $\mathrm{C} = (\mathrm{C}, \alpha, I)$, simply called C-*semigroup*, $(X, \beta)$ is an $S$-colored semigroup $\{X(s)\}_{s\in S}$ such that $(X, \beta)$ is a C-set and the C-action is compatible with the semigroup structure: For all $c \in \mathrm{C}(s_0; \varphi(s_1), ..., \varphi(s_m))$ and $x_j, x_i' \in X(s_j)$, $1 \leqslant j \leqslant m$,

$$\beta(c, x_1x_1', x_2x_2', ..., x_mx_m') = \beta(c, x_1, x_2, ..., x_m)\beta(c, x_1', x_2', ..., x_m'), \tag{2.6}$$

where on the right hand side the product is given by the semigroup structure. It is a C-*monoid* or a *monoid over* C if additionally for all $c \in \mathrm{C}(s_0; \varphi(s_1), ..., \varphi(s_m))$ and units $u_i \in X(s_i)$, $0 \leqslant i \leqslant m$, we have $\beta(c, u_1, u_2, ...., u_m) = u_0$.

**Definition 2.17.** Let $R$ be a commutative ring. An *$R$-module over* $\mathrm{C} = (\mathrm{C}, \alpha, I)$, simply called a *C-$R$-module*, is an $S$-colored $R$-module $\{X(s)\}_{s\in S}$ such that $\{X(s)\}_{s\in S}$ is a C-set and the C-action is compatible with the $R$-module structure: for all $c \in \mathrm{C}(s_0; \varphi(s_1), ..., \varphi(s_m))$, $r, r' \in R$, $x_j \in X(s_j)$, $1 \leqslant j \leqslant m$, and $x_i' \in X(s_i)$, we have

$$\begin{aligned} & \beta(c, x_1, ..., x_{i-1}, rx_i + r'x_i', x_{i+1}, ...., x_m) && (2.7)\\ = \quad & r\beta(c, x_1, ..., x_{i-1}, x_i, ...., x_m) + r'\beta(c, x_1, ..., x_{i-1}, x_i', x_{i+1}, ...., x_m). && (2.8)\end{aligned}$$

It is an *$R$-algebra over* C, simply called a *C-$R$-algebra*, if additionally the C-action is compatible with the multiplication structure in the sense of (2.6).

**Definition 2.18.** A *morphism* of C-sets (resp. C-semigroups, C-monoids, C-modules, C-algebras) from $(X, \beta^X)$ to $(Y, \beta^Y)$ is a collection of morphisms of sets (resp. semigroups, monoids, modules, algebras) $\Psi = \{\Psi_s : X(s) \to Y(s)\}_{s\in S}$ such that for all $\underline{s} = (s_0; \varphi(s_1), ..., \varphi(s_m)) \in \bigsqcup_{k\geqslant 0} S^{k+1}$

$$\beta^Y \circ (\mathrm{id} \times \Psi_{s_1} \times \cdots \times \Psi_{s_m}) = \Psi_{s_0} \circ \beta^X : \mathrm{C}(\underline{s}) \times X(s_1) \times \cdots \times X(s_m) \to Y(s_0).$$

**Remark 2.19.** As usual, the "pullback" of a morphism $\Phi : \mathrm{C} = (\mathrm{C}, \alpha^{\mathrm{C}}, I^{\mathrm{C}}) \to \mathrm{D} = (\mathrm{D}, \alpha^{\mathrm{D}}, I^{\mathrm{D}})$ of $(S, \varphi)$-colored operads defines for each D-set $(X, \beta)$ the structure of a C-set $(X, \beta \circ \Phi)$. This implies in particular, for the choice $(\mathrm{C}, \mathrm{D}) = (\mathrm{C}, \overline{\mathrm{C}})$ with the isomorphism from Lemma 2.13 that C-sets are the same as $\overline{\mathrm{C}}$-sets.

2.5. **C-category (Category over an $(S, \varphi)$-colored operad).** Let $\mathrm{C} = (\mathrm{C}, \alpha, I)$ be an $(S, \varphi)$-colored operad. In this section we will use the previous definitions (C-sets, C-semigroups and C-monoids) to define the notion of colored categories with an operad C action. The C-set structure is used for the action on object sets and the $\overline{\mathrm{C}}$-semigroup and $\overline{\mathrm{C}}$-monoid structure is used for the action on morphism sets.

**Definition 2.20.** A *category over* C, simply called a C-*category*, $\mathcal{X} = (\mathcal{X}, \beta)$ is an $S$-colored category $\mathcal{X} = \{\mathcal{X}(s)\}_{s\in S}$ with

- structure maps $\beta = \{\beta_{\underline{s}}\}_{\underline{s}\in\bigsqcup_{k\geqslant 0} S^{k+1}}$

$$\beta : \quad \mathrm{C}(s_0; \varphi(s_1), ..., \varphi(s_m)) \times \mathrm{Ob}(\mathcal{X}(s_1)) \times \cdots \times \mathrm{Ob}(\mathcal{X}(s_m)) \to \mathrm{Ob}(\mathcal{X}(s_0)) \tag{2.9}$$

  such that $\{\mathrm{Ob}(\mathcal{X}(s))\}_{s\in S}$ turns into a C-set,
- structure maps $\overline{\beta} = \{\overline{\beta}_{\underline{s}}\}_{\underline{s}\in\bigsqcup_{k\geqslant 0} S^{k+1}}$

$$\begin{aligned}\overline{\beta} : \quad & \overline{\mathrm{C}}(s_0; \varphi(s_1), ..., \varphi(s_m)) \times \mathcal{X}(s_1)(A_1, B_1) \times \cdots \times \mathcal{X}(s_m)(A_m, B_m)\\ & \to \mathcal{X}(s_0)(\beta(c, A_1, ..., A_m), \beta(c, B_1, ..., B_m))\end{aligned} \tag{2.10}$$

  such that the $S$-colored semigroup

$$\{\mathcal{M}\mathrm{or}(s)\}_{s\in S} := \left\{ \bigcup_{(A,B)\in\mathrm{Ob}(\mathcal{X}(s))^2} \mathcal{X}(s)(A, B) \right\}_{s\in S} \tag{2.11}$$

  given by composition of maps is a $\overline{\mathrm{C}}$-semigroup and its restriction to $A = B$ is a $\overline{\mathrm{C}}$-monoid with units given by the identity morphisms.

**Example 2.21.** Consider the $(\mathbb{Z}_{\geqslant 1}, \mathrm{id})$-colored operad $D$ defined by

$$\mathrm{D}\left(\sum_{i=1}^{m} s_i; s_1, ..., s_m\right) = \left\{ \begin{array}{l} \text{The } m\text{-leaved planar rooted tree} \\ \text{with coloring } \sum_{i=1}^{m} s_i \text{ on the root} \\ \text{and the coloring } s_1, ..., s_m \text{ on the } m \text{ leaves} \end{array} \right\}$$

The multiplication $\alpha$ is analogous to Example 2.9.

Now fix a monoidal category $(\mathcal{C}, \otimes)$ and let $\mathcal{X}(s)$ for $s \in \mathbb{Z}_{\geqslant 1}$ be the full subcategory with objects $X_1 \otimes \cdots \otimes X_s$, $X_i \in \mathrm{Ob}(\mathcal{C})$. This gives a $\mathbb{Z}_{\geqslant 1}$-colored category. It is even a $D$-category with the structure map

$$\begin{array}{rccc} \beta: & \mathrm{D}(\underline{s}) \times \mathrm{Ob}(\mathcal{X}(s_1)) \times \cdots \times \mathrm{Ob}(\mathcal{X}(s_m)) & \twoheadrightarrow & \mathrm{Ob}(\mathcal{X}(\sum_{i=1}^{m} s_i)) \\ & \uplus & & \uplus \\ & (\mathrm{D}(\underline{s}), X^{(1)}, ..., X^{(m)}) & \longmapsto & X^{(1)} \otimes \cdots \otimes X^{(m)}, \end{array}$$

where $\underline{s} = (\sum_{i=1}^{m} s_i; s_1, ..., s_m)$. The structure map $\overline{\beta}$ is canonically induced by the structure map $\beta$. The compatibility condition between $\alpha$ and $\beta$ follows from the coherence of the monoidal category.

**Definition 2.22.** A functor $\Psi$ between two C-categories $(\mathcal{X}, \beta^{\mathcal{X}})$ and $(\mathcal{Y}, \beta^{\mathcal{Y}})$ is a collection of functors $\Psi = \{\Psi_s : \mathcal{X}(s) \to \mathcal{Y}(s)\}_{s \in S}$ such that:

- the collection of morphisms $\Psi_s : \mathrm{Ob}(\mathcal{X}(s)) \to \mathrm{Ob}(\mathcal{Y}(s))$ is a morphism of C-sets,
- for any $s \in S$ and $A \in \mathrm{Ob}(\mathcal{X}(s))$, the maps $\Psi_s : \mathcal{X}(s)(A, A) \to \mathcal{Y}(s)(\Psi(A), \Psi(A))$ define a morphism of $\overline{\mathrm{C}}$-monoids,
- for any $s \in S$ and $A, B \in \mathrm{Ob}(\mathcal{X}(s))$ the maps $\Psi_s : \mathcal{X}(s)(A, B) \to \mathcal{Y}(s)(\Psi(A), \Psi(B))$ define a morphism of $\overline{\mathrm{C}}$-semigroups.

2.6. **From preadditive C-categories to homotopy C-categories.** In the rest of this section, we define a preadditive C-category $\mathcal{X}^{p.add}$, an additive C-category $\mathcal{X}^{add}$, the category of bounded complexes over an additive C-category $\mathrm{Com}^{\mathrm{b}}(\mathcal{X}^{add})$ and finally its homotopy C-category $\mathrm{K}^{\mathrm{b}}(\mathcal{X}^{add})$. These categories are analogous to Bar-Natan's definition of categories of complexes over cobordism categories [BN05].

Recall that a preadditive category is a category where the morphisms sets are abelian groups and the composition maps are bilinear.

**Definition 2.23.** A C-category $(\mathcal{X}, \beta)$ is preadditive if $\mathcal{X}(s)$ is preadditive for all $s \in S$ and the structure maps are compatible with the preadditive structure:
For $\mathrm{id}_c \in \overline{\mathrm{C}}(s; \varphi(s_1), ..., \varphi(s_m))$ and $f_j = \sum_{i=1}^{n_j} f_j^i \in \mathcal{X}(s_j)(X_j, X'_j)$ $(1 \leqslant j \leqslant m)$, the structure map $\overline{\beta}$ satisfies

$$\begin{aligned} \overline{\beta}(\mathrm{id}_c, f_1, ..., f_m) &= \sum_{i_1=1}^{n_1} \cdots \sum_{i_m=1}^{n_m} \overline{\beta}(\mathrm{id}_c, f_1^{i_1}, ..., f_m^{i_m}) \\ &\in \mathcal{X}(s)(\beta(c, X_1, ..., X_m), \beta(c, X'_1, ..., X'_m)). \end{aligned} \tag{2.12}$$

**Remark 2.24.** Given a C-category $\mathcal{X} = (\mathcal{X}, \beta)$, its *associated preadditive* C*-category* $\mathcal{X}^{p.add} = (\mathcal{X}^{p.add}, \beta)$ is obtained by allowing formal $\mathbb{Z}$-linear sums of morphisms and extending the structure maps linearly. Similarly, for any $\mathbb{Z}$-algebra $A$, we obtain the *scalar extended preadditive* C*-category* $\mathcal{X}_A^{p.add} = (\mathcal{X}_A^{p.add}, \beta)$ associated to $(\mathcal{X}, \beta)$ by extending the scalars to $A$. To summarize:

$$\begin{aligned} \mathcal{X}, \text{a C} - \text{category} &\rightsquigarrow \mathcal{X}^{p.add}, \text{a preadditive C} - \text{category} \\ &\rightsquigarrow \mathcal{X}_A^{p.add}, \text{a preadditive C} - \text{category with scalars extended to } A. \end{aligned}$$

For objects $A_j, B_j \in \mathcal{X}(s_j)$ with $A_j = \bigoplus_{k=1}^{n_j} A_j^k$, $B_j = \bigoplus_{l=1}^{n'_j} B_j^l$ a morphism $f_j \in \mathcal{X}(s_j)(A_j, B_j)$ is given as a $(n_j, n'_j)$-matrix

$$f = (f_j^{k,l})_{k,l}. \tag{2.13}$$

**Definition 2.25.** A preadditive C-category is *additive* if the underlying categories $\mathcal{X}(s)$, $s \in S$, are additive and the structure maps are compatible with the additive structure:

• for $c \in \mathrm{C}(s_0;\varphi(s_1),...,\varphi(s_m))$ and $X_j^{i_j} \in \mathrm{Ob}(\mathcal{X}(s_j))$ $(1 \leqslant j \leqslant m,\ 1 \leqslant i_j \leqslant n_j)$, the structure maps $\beta$ satisfy

$$\beta(c, \bigoplus_{i_1=1}^{n_1} X_1^{i_1}, ..., \bigoplus_{i_m=1}^{n_m} X_m^{i_m}) = \bigoplus_{i_1=1}^{n_1} \cdots \bigoplus_{i_m=1}^{n_m} \beta(c, X_1^{i_1}, ..., X_m^{i_m}), \tag{2.14}$$

• the structure map applied to matrices equals the matrix of the structure maps applied to the entries, in formulas: for $\mathrm{id}_c \in \overline{\mathrm{C}}(s_0;\varphi(s_1),...,\varphi(s_m))$ and $f_j \in \mathcal{X}^{add}(s_j)(A_j, B_j)$ as in (2.13) for $1 \leqslant j \leqslant m$ we have

$$\overline{\beta}(\mathrm{id}_c, f_1, ..., f_m) = \left(\overline{\beta}(\mathrm{id}_c, f_1^{k_1,l_1}, ..., f_m^{k_m,l_m})\right) \tag{2.15}$$

with $1 \leqslant k_r \leqslant n_j$, $1 \leqslant l_r \leqslant n'_j$ for $1 \leqslant j \leqslant m$ holds.

A preadditive C-category $\mathcal{X}$ can formally be extended to an additive C-category:

**Definition 2.26.** The *additive closure of a preadditive* C*-category* $\mathcal{X}$, denoted by

$$\mathrm{Mat}(\mathcal{X}) := \{\mathrm{Mat}(\mathcal{X}(s))\}_{s \in S},$$

is defined as follows: The object set of each $S$-colored category $\mathrm{Mat}(\mathcal{X}(s))$ consists of the formal direct sums $\oplus_{i=1}^n X^i$ of objects $X^i$ in $\mathrm{Ob}(\mathcal{X}(s))$ with the zero object $0_s$ corresponding to the sum over the empty set. Morphisms in $\mathrm{Mat}(\mathcal{X}(s))(\oplus_{k=1}^n X^k, \oplus_{l=1}^{n'} X'^l)$ are $n' \times n$ matrices $(f^{l,k})$ of morphisms $f^{l,k} \in \mathcal{X}^{p.add}(s)(X^k, X'^l)$, as in (2.13) with composition defined in terms of matrix multiplication. The structure maps $\beta$, $\overline{\beta}$ are assumed to satisfy (2.14) and (2.15).

For an additive C-category $\{\mathcal{X}(s)\}_{s \in S}$ we have the usual categories of bounded complexes $\mathrm{Com^b}(\mathcal{X}(s))$, $s \in S$. We like to make them compatible by adding structure maps generalizing the construction of tensor product on complexes:

**Definition 2.27.** Let $\mathcal{X} = \{\mathcal{X}(s)\}_{s \in S}$ be an additive C-category. The associated C-*category of bounded complexes* is given by the $S$-colored category

$$\mathrm{Com^b}(\mathcal{X}) := \{\mathrm{Com^b}(\mathcal{X}(s))\}_{s \in S},$$

of bounded complexes with the following structure maps (with $1 \leqslant j \leqslant m$):
• for $c \in \mathrm{C}(s_0;\varphi(s_1),...,\varphi(s_m))$ and bounded complexes

$$X_j^\bullet = \cdots \xrightarrow{d_{X_j^{i-1}}} X_j^i \xrightarrow{d_{X_j^i}} X_j^{i+1} \xrightarrow{d_{X_j^{i+1}}} \cdots$$

in $\mathrm{Com^b}(\mathcal{X}(s_j))$ the structure map $\beta^{\mathrm{Com^b}(\mathcal{X})}$ defined by

$$\beta(c, X_1^\bullet, ..., X_m^\bullet) \quad = \quad \cdots \xrightarrow{d_{i-1}} X^i \xrightarrow{d_i} X^{i+1} \xrightarrow{d_{i+1}} \cdots,$$

where $X^i = \bigoplus_{i_1+\cdots+i_m=i} \beta(c, X_1^{i_1}, ..., X_m^{i_m})$ and

$$d_i \quad = \sum_{i_1+\cdots+i_m=i} \sum_{k=1}^m (-1)^{i_1+\cdots+i_{k-1}} \overline{\beta}(\mathrm{id}_c, \mathrm{id}_{X_1^{i_1}}, ..., \mathrm{id}_{X_{k-1}^{i_{k-1}}}, d_{X_k^{i_k}}, \mathrm{id}_{X_{k+1}^{i_{k+1}}}, ..., \mathrm{id}_{X_m^{i_m}});$$

• for $\mathrm{id}_c \in \overline{\mathrm{C}}(s_0;\varphi(s_1),...,\varphi(s_m))$ and $f_j^\bullet \in \mathrm{Com^b}(\mathcal{X}(s_j))(X_j^\bullet, X'^\bullet_j)$ $(1 \leqslant j \leqslant m)$, the structure map $\overline{\beta}$ defined by

$$\begin{aligned} \overline{\beta}(\mathrm{id}_c, f_1^\bullet, ..., f_m^\bullet) &= (..., f^{i-1}, f^i, f^{i+1}, ...) \\ &\in \mathrm{Com^b}(\mathcal{X}^{p.add}(s))(\beta(c, X_1^\bullet, ..., X_m^\bullet), \beta(c, X'^\bullet_1, ..., X'^\bullet_m)), \end{aligned} \tag{2.16}$$

$$\begin{aligned} \text{where } f^i &= \sum_{i_1+\cdots+i_m=i} \overline{\beta}(\mathrm{id}_c, f_1^{i_1}, ..., f_m^{i_m}) \in \mathcal{X}^{p.add}(s)(X^i, X'^i), \\ X^i &= \oplus_{i_1+\cdots+i_m=i}\beta(c, X_1^{i_1}, ..., X_m^{i_m}), \\ X'^i &= \oplus_{i_1+\cdots+i_m=i}\beta(c, X'^{i_1}_1, ..., X'^{i_m}_m). \end{aligned}$$

**Remark 2.28.** The above operad action on a category of complexes is analogous the tensor product of complexes in the sense that it works structurally the same way. As in the case of tensor products of complexes, the choice of signs ensures $d^2 = 0$.

The following is the analog of the fact, that for $f_i : A_i \to B_i$, $1 \leqslant i \leqslant m$, with $f_1$ null-homotopic $f_1 \otimes f_2 \otimes \ldots \otimes f_m : A_1 \otimes \ldots \otimes A_m \to B_1 \otimes \ldots \otimes B_m$ is also null-homotopic.

**Proposition 2.29.** *Let* $f_i^\bullet \in \mathrm{Com}^{\mathrm{b}}(\mathcal{X}(s_i))(X_i^\bullet, Y_i^\bullet)$ *for* $1 \leqslant i \leqslant m$. *Assume at least one of the* $f_i^\bullet$ *is null-homotopic. Then* $\overline{\beta}(\mathrm{id}_c, f_1^\bullet, ..., f_m^\bullet)$ *is null-homotopic, where* $\mathrm{id}_c \in \overline{\mathrm{C}}(s_0; \varphi(s_1), ..., \varphi(s_m))$.

*Proof.* Without loss of generality, we assume $f_1$ is null-homotopic. Let $h = \{h^i : X_1^i \to Y_1^{i-1}\}$ be a homotopy. It suffices to consider the case $m = 2$. The maps

$$\bigoplus_{a+b=i} \overline{\beta}(\mathrm{id}_c, h^a, f_2^b) : \bigoplus_{a+b=i} \beta(c, X_1^a, X_2^b) \to \bigoplus_{a+b=i} \beta(c, X'^{a-1}_1, X'^b_2)$$

define a homotopy for $\overline{\beta}(\mathrm{id}_c, f_1^\bullet, f_2^\bullet)$. □

As expected, we use the last proposition to define the homotopy category.

**Definition 2.30.** The *homotopy* C*-category of* $\mathrm{Com}^{\mathrm{b}}(\mathcal{X})$, denoted by

$$\mathrm{K}^{\mathrm{b}}(\mathcal{X}) := \{\mathrm{K}^{\mathrm{b}}(\mathcal{X}(s)))\}_{s \in S},$$

is defined as follows:
- the object set $\mathrm{Ob}(\mathrm{K}^{\mathrm{b}}(\mathcal{X}(s)))$ is the same as $\mathrm{Ob}(\mathrm{Com}^{\mathrm{b}}(\mathcal{X}(s)))$.
- the morphisms in $\mathrm{K}^{\mathrm{b}}$ are the morphisms in $\mathrm{Com}^{\mathrm{b}}$ up to null-homotopic, i.e.

$$\mathrm{K}^{\mathrm{b}}(\mathcal{X}(s))(X^\bullet, X'^\bullet) := \mathrm{Com}^{\mathrm{b}}(\mathcal{X}(s))(X^\bullet, X'^\bullet)/\{\text{null-homotopic morphisms}\}.$$

**Remark 2.31.** The constructions of C-categories $\mathrm{Com}^{\mathrm{b}}$ and $\mathrm{K}^{\mathrm{b}}$ work analogously for bounded above or bounded below complexes.

## 3. Extention of category with colored operad action

In this section we introduce C-sets generated by colored elements and C-categories extended by colored morphisms.

3.1. **C-set generated by colored elements.** Let C be an $(S, \varphi)$-colored operad. We introduce C-sets generated by an element $X_s$ with color $s \in S$.

Given a formal symbol $X_s$ with color $s \in S$, we define for $c \in \mathrm{C}(\underline{s})$, where $\underline{s} = (s_0, \underbrace{\varphi(s), ..., \varphi(s)}_{m})$, the formal elements

$$\hat{\alpha}(c; X_s^m) := \hat{\alpha}_{\underline{s}}(c; \underbrace{X_s, ..., X_s}_{m}). \tag{3.1}$$

We set $\hat{\alpha}(c; X_s^0) := c$ for $c \in \mathrm{C}(s_0)$.

**Definition 3.1.** The C-set $(X_s)_{\mathrm{C}}$ *generated* by $X_s$ is the $S$-colored set $\{(X_s)_{\mathrm{C}}(s_0)\}_{s_0 \in S}$, where

$$(X_s)_{\mathrm{C}}(s_0) := \{\hat{\alpha}(c, X_s^m) \mid m \geqslant 0, c \in \mathrm{C}(s_0; \underbrace{\varphi(s), ..., \varphi(s)}_{m})\}, \tag{3.2}$$

together with the structure maps (for non-empty sets)

$$\beta_{\underline{s}} : \mathrm{C}(\underline{s}) \times (X_s)_{\mathrm{C}}(s_1) \times \cdots \times (X_s)_{\mathrm{C}}(s_m) \to (X_s)_{\mathrm{C}}(s_0), \tag{3.3}$$

for $\underline{s} = (s_0; \varphi(s_1), ..., \varphi(s_m))$, defined as

$$\beta_{\underline{s}}(c; x_1, ..., x_m) := \hat{\alpha}(\alpha_{\underline{s}}(c, c_1, ..., c_m); X_s^n), \tag{3.4}$$

where $x_i = \hat{\alpha}(c_i; X_s^{n_i})$ and $n = n_1 + \cdots + n_m$, and for the unit $I_s \in \mathrm{C}(s; \varphi(s))$

$$\beta_{(s;\varphi(s))}(I_s; x_s) = x_s. \tag{3.5}$$

We omit the straightforward check that these structure maps satisfy the required conditions of a C-set which follows directly from the conditions for $\alpha$.

**Example 3.2.** Let $X_+, X_-$ be two formal symbols with the same color $s \in S$. The C-set $(X_+, X_-)_{\mathrm{C}}$ is the $S$-colored set $\{(X_+, X_-)_{\mathrm{C}}(s_0)\}_{s_0 \in S}$, where

$$(X_+, X_-)_{\mathrm{C}}(s_0) := \{\alpha(c, X_{\epsilon_1}, \ldots, X_{\epsilon_m}) \mid m \geqslant 0, \epsilon_i \in \{+, -\}, c \in \mathrm{C}(s_0; \underbrace{\varphi(s), ..., \varphi(s)}_{m})\} \tag{3.6}$$

and the structure maps are defined as before.

**Remark 3.3.** The definition extends to arbitrary finite sets $X = \{X_i \mid i \in I\}$ with coloring map $X \to S$, but we won't need this here.

3.2. **Category extended by morphisms.** Let $\mathcal{A}$ be a category. We consider additional morphisms of the category $\mathcal{A}$ which are morphisms between objects of $\mathcal{A}$ but not included in the morphism set of $\mathcal{A}$. A set $F$ of *additional morphisms* is a disjoint union

$$F = \bigsqcup_{X,Y \in \mathrm{Ob}(\mathcal{A})} F(X, Y),$$

of possibly empty sets. The subset $F(X, Y)$ is the set of *additional morphisms* from $X$ to $Y$ (possibly empty if there are no morphisms from $X$ to $Y$ in $F$).

**Definition 3.4** (Category extended by morphisms)**.** Let $\mathcal{A}$ be a category and $F$ be a set of additional morphisms of $\mathcal{A}$. For any two objects $X, Y$ in $\mathcal{A}$ set

$$\mathcal{A}'(X, Y) = \mathcal{A}(X, Y) \cup F(X, Y).$$

The associated extended category $\mathcal{A}[F]$ is defined by $\mathrm{Ob}(\mathcal{A}[F]) = \mathrm{Ob}(\mathcal{A})$ and

$$\mathcal{A}[F](X, Y) = \{f_r \ldots f_2 f_1 | r \in \mathbb{Z}_{>0}, f_i \in \mathcal{A}'(X_i, X_{i+1}), X_1 = X, X_{r+1} = Y\}$$

where we factor out associativity and identity: $f \circ \mathrm{id}_X = f = \mathrm{id}_Y \circ f$ for all $f \in \mathcal{A}[F](X, Y)$.

3.3. C**-category extended by colored morphisms.** Let $(\mathcal{X}, \beta)$ be a $C$-category. By *additional morphisms* for C we mean an $S$-colored set $G = \{G(s)\}_{s \in S}$ of additional morphisms $G(s)$ for $\mathcal{X}(s)$. For $c \in \mathrm{C}(s_0; s_1, \ldots, s_m)$ and $\mathbf{g} = (g_i)_{1 \leqslant i \leqslant m}$ with $g_i \in \mathcal{X}(s_i)[G(s_i)](A_i, B_i)$ for $A_i, B_i \in \mathcal{X}(s_i)$ we extend $\overline{\beta}$ by adding formal symbols

$$\overline{\beta}(\mathrm{id}_c; \mathbf{g}) := \overline{\beta}(\mathrm{id}_c; g_1, \ldots, g_m).$$

in case not all $g_i$ are in $\mathcal{X}(s_i)$. Let $G^{\mathrm{C}}(s_0)$ be the equivalence classes of all $\overline{\beta}(\mathrm{id}_c; \mathbf{g})$ (for varying $\mathbf{g}$) modulo the relation (2.6) and let $G^{\mathrm{C}} = \bigsqcup_{s_0 \in S} G^{\mathrm{C}}(s_0)$. In the following we do not distinguish in the notation between $\overline{\beta}(\mathrm{id}_c; \mathbf{g})$ and its equivalence class. Given this data we have the extended $S$-colored category

$$\mathcal{X}[G] = \{\mathcal{X}[G](s)\}_{s \in S} := \{\mathcal{X}(s)[G^{\mathrm{C}}(s)]\}_{s \in S}. \tag{3.7}$$

Note that the objects are the same as in $\mathcal{X}$, but the morphisms are extended. We define now the corresponding extended $C$-category.

**Definition 3.5** (Extended C-category)**.** Let $(\mathcal{X}, \beta)$ be a $C$-category and $G$ be a set of additional $S$-colored morphisms of $\mathcal{X}$. The associated *extended* C*-category* $(\mathcal{X}[G], \beta^{\mathcal{X}[G]})$ is the $S$-colored category $\mathcal{X}[G]$ as in (3.7) together with the structure maps $\beta^{\mathcal{X}[G]} := \beta$ and $\overline{\beta}^{\mathcal{X}[G]}$ given by extending the $\overline{\beta}$'s linearly in the sense of (2.7).

Although the above definition looks rather technical it should be seen just as the $C$-category generated by the set of additional morphisms $G$.

## 4. Colored operad P

In this section we define oriented planar arc diagrams and an operad structure on them. In the next section, we introduce then a new family of cobordism categories for Khovanov–Rozansky homology associated to matrix factorization with the potentials $\frac{n}{n+1}x^{n+1}$ and realize them as a category over the operad of arc diagrams.

4.1. **The $(S_\pm, \iota)$-colored operad** P**.** Let $D_0$ be a closed unit disk in $\mathbb{R}^2$ with center the origin $(0,0)$ and let $D(k)$ $(k \in \mathbb{Z}_{\geqslant 0})$ be the complement in $D_0$ of the finite union of $k$ open disks $D_i$ $(i = 1, ..., k)$ with diameter $\frac{1}{k}$ and center $(\frac{2i}{k} - \frac{1}{k} - 1, 0)$

$$D(k) = D_0 \backslash \bigcup_{i=1}^{k} D_i,$$

together with a fixed base point on each circle. For instance, we fix a base point $*_j$ at the angle $\frac{3\pi}{2}$ of each boundary $\partial D_j$ $(j = 0, ..., k)$, for instance

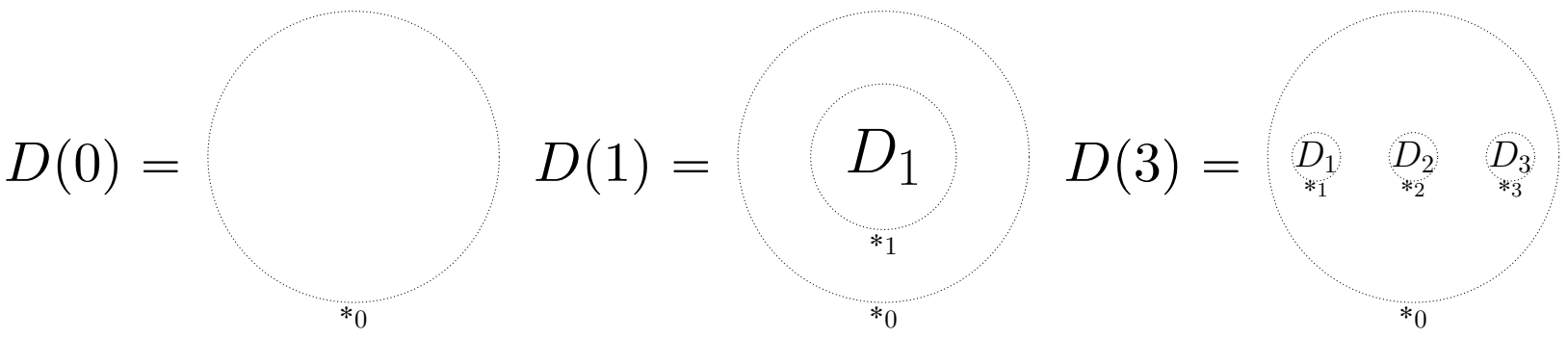

**Definition 4.1.** Let $I = [0,1]$ be the unit interval oriented from 0 to 1 and $S^1$ the unit circle in $\mathbb{R}^2$ oriented counter-clockwise.

An *oriented $(l_1, l_2)$-arc diagram of type $k$* is a finite union $T$ of embeddings $f_i : I \to D(k)$, $1 \leqslant i \leqslant l_1$, and $g_j : S^1 \to D(k)$, $1 \leqslant j \leqslant l_2$, such that $f_i(0), f_i(1) \in \partial(D(k)) \backslash \{*_0, ..., *_k\}$ and the images are pairwise disjoint. We usually identify $T$ withe the image $\bigcup_{a=1}^{l_1} f_a \cup \bigcup_{b=1}^{l_2} g_b$ of the embeddings.

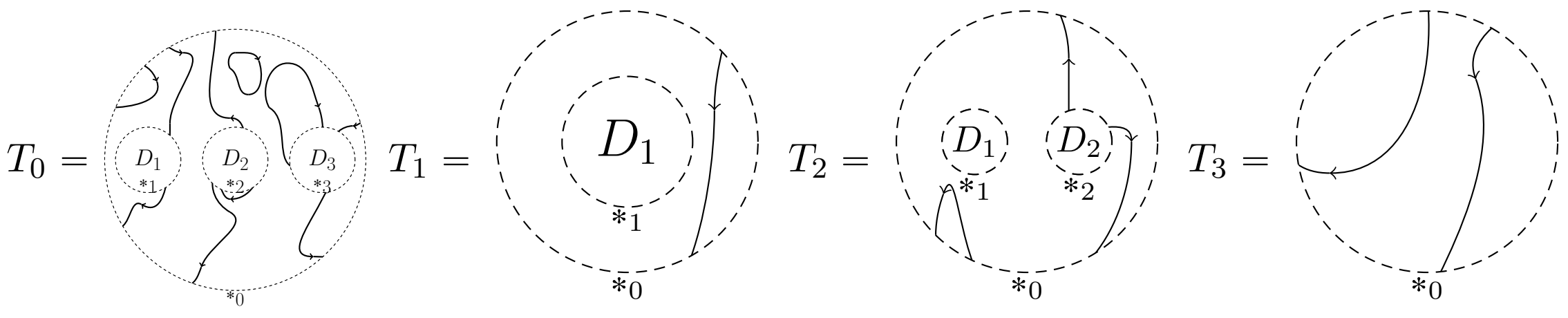

Figure 2. Planar arc diagrams

**Definition 4.2** (Isotopy of arc diagrams)**.** Let $T$ and $T'$ be oriented $(l_1, l_2)$-arc diagrams of type $k$. $T$ and $T'$ are *isotopic* if there exists isotopy, i.e. a continuous map

$$\varphi : \quad \begin{array}{ccc} D(k) \times [0,1] & \longrightarrow & D(k) \\ \Downarrow & & \Downarrow \\ (d,t) & \longmapsto & \varphi_t(d) \end{array}$$

such that each $\varphi_t$ is a homeomorphism of $D(k)$, with $\varphi_t(*_j) = *_j$ $(j = 0, ..., k)$, $\varphi_0 = \mathrm{id}_{D(k)}$ and $\varphi_1(T) = T'$. Write $T \sim T'$ if $T$ and $T'$ are isotopic.

For a given $(l_1, l_2)$-arc diagram $T = \bigsqcup_{a=1}^{l_1} f_a \cup \bigsqcup_{b=1}^{l_2} g_b$ of type $k$, we define a *coloring* $\underline{s}(T)$ as follows. Put a $+$-sign at $f_i(1)$ and a $-$-sign at $f_i(0)$ $(i = 1, ..., l_1)$ and define $s_j(T)$ $(j = 0, ..., k)$ by a sequence of $\pm$-signs on $\partial D_j$ reading signs counter-clockwise from the base point $*_j$. We set $s_j(T) = \varnothing$ if there is no $\pm$-sign on $\partial D_j$. Finally we define a coloring of $T$ by

$$\underline{s} = \underline{s}(T) = (s_0(T); s_1(T), ..., s_k(T)).$$

Denote by $\#_+(\underline{s})$ the total number of $+$-sign in $\underline{s}$ and by $\#_-(\underline{s})$ the total number of $-$-sign in $\underline{s}$. By definition we find the following lemma. Hence by definition

**Lemma 4.3.** *For any arc diagrams $T$ we have $\#_+(\underline{s}(T)) = \#_-(\underline{s}(T))$.*

Going back to the arc diagrams of Figure 2, we find that the coloring $\underline{s}(T_0)$ consists of $s_0(T_0) = (+, -, -, -, +, +, +)$, $s_1(T_0) = (-, +)$, $s_2(T_0) = (-, -, +)$, $s_3(T_0) = (-, +, +, -)$, $\underline{s}(T_1)$ consists of $s_0(T_1) = (+, -)$, $s_1(T_1) = \varnothing$, the coloring $\underline{s}(T_2)$ consists of $s_0(T_2) = (+, +, -)$, $s_1(T_2) = \varnothing$, $s_2(T_2) = (-)$ and the coloring $\underline{s}(T_3)$ consists of $s_0(T_3) = (+, -, -, +)$.

Let $S_\pm$ be the set $\bigsqcup_{k \geqslant 1} \{+, -\}^k \cup \{\varnothing\}$. We define an involution $\iota : S_\pm^{k+1} \to S_\pm^{k+1}$ by $\iota(\varnothing) = \varnothing$, $\iota(+) = -$, $\iota(-) = +$, $\iota(\epsilon_1, ..., \epsilon_n) = (\iota(\epsilon_1), ..., \iota(\epsilon_n))$ for $(\epsilon_1, ..., \epsilon_n) \in S_\pm$ and $\iota(s_0; s_1, ..., s_k) = (\iota(s_0); \iota(s_1), ..., \iota(s_k))$ for $(s_0; s_1, ..., s_k) \in S_\pm^{k+1}$.

**Definition 4.4.** For a given $\underline{s} = (s_0; \iota(s_1), ..., \iota(s_k)) \in S_\pm^{k+1}$, we consider the isotopy classes of arc diagrams $\mathrm{P}(\underline{s}) = \{T \mid \underline{s}(T) = \underline{s}\} / \sim$ . We define the $\bigsqcup_{k \geqslant 0} S_\pm^{k+1}$-colored set $\mathrm{P}$ as

$$\mathrm{P} = \{\mathrm{P}(\underline{s})\}_{\underline{s} \in \bigsqcup_{k \geqslant 0} S_\pm^{k+1}}.$$

Obviously, by Lemma 4.3 we have $\mathrm{P}(\underline{s}) = \varnothing$ if $\#_+(\underline{s}) \neq \#_-(\underline{s})$.

4.2. **Colored operad structure of** P**.** In this section we define an $(S_\pm, \iota)$-operad structure on P. We start by illustrating the operad product before giving precise definitions.

For arc diagrams $T_0 \in \mathrm{P}(s_0; \iota(s_0^{(1)}), ..., \iota(s_0^{(m)}))$ and $T_i \in \mathrm{P}(s_0^{(i)}; \iota(s_1^{(i)}), ..., \iota(s_{k_i}^{(i)}))$ $(i = 1, ..., m)$, the operad product $\alpha$ first takes $T_i$ and puts it into the $i$-th hole of $T_0$ by matching the base points and connecting the arcs at the boundaries, and then relabels and rescales the discs in the obvious way. Using for example the arc diagrams of Figure 2 the product $\alpha(T_0;T_1,T_2,T_3)$ looks as follows:

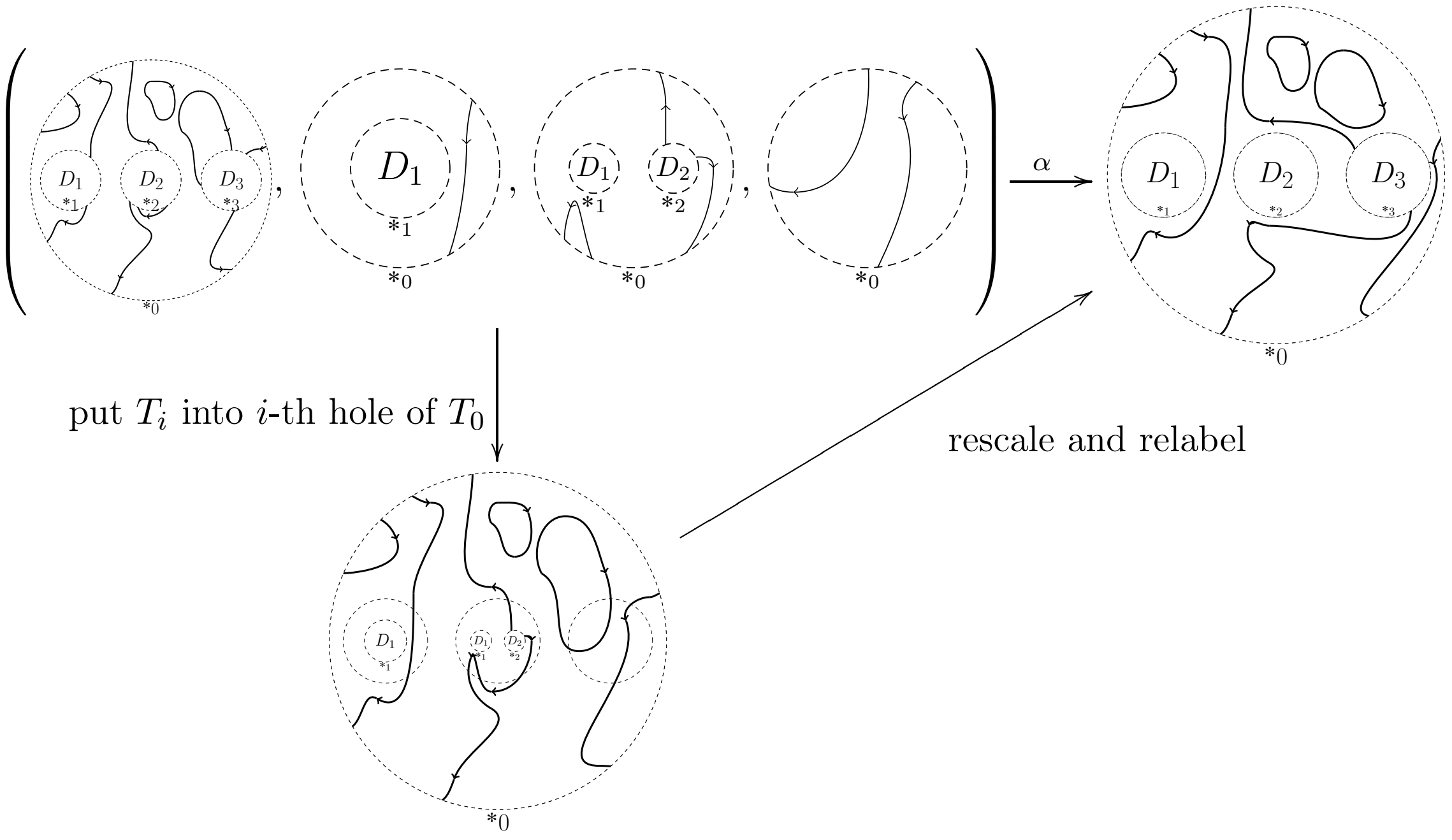

• **The operad product** $\alpha$:
Let

$$T_0 = \bigcup_{a=1}^{l_1} f_a \cup \bigcup_{b=1}^{l_2} g_b$$

be an $(l_1, l_2)$-arc diagram of type $m$ with coloring $\underline{s}(T_0) = (s_0; \iota(s_0^{(1)}), ..., \iota(s_0^{(m)})) \in S_{\pm}^{m+1}$ and let

$$T_i = \bigcup_{a=1}^{l_1^{(i)}} f_a^{(i)} \cup \bigcup_{b=1}^{l_2^{(i)}} g_b^{(i)} \quad (i = 1, ..., m)$$

be an $(l_1^{(i)}, l_2^{(i)})$-arc diagram of type $k_i$ with coloring $\underline{s}(T_i) = (s_0^{(i)}; \iota(s_1^{(i)}), ..., \iota(s_{k_i}^{(i)})) \in S_{\pm}^{k_i+1}$, where $s_0^{(i)}$ is the coloring on $\partial D_i$ of $T_0$. Let $\epsilon_a^{(i)} \in \{+,-\}$ be the sign at $f_{j_a}(t_a) \in \partial D_i$ respectively $f_{j'_a}^{(i)}(1-t_a) \in \partial D_0$ with $1 \leqslant a \leqslant i$, $j_a \in \{1, ..., l_1\}$, $j'_a \in \{1, ..., l_1^{(i)}\}$ and $t_a \in \{0,1\}$.

There exists an embedding $e_i : D_0 \to D_i$ $(i = 1, ..., m)$ such that $e_i$ on

$$D_\theta = \{d \in D_0 | \, |d| < 1 - \theta\} \subset D_0 \quad (\theta \ll 0)$$

equals to the scaling embedding $\frac{1}{k} : D_0 \to D_i; d \mapsto \frac{1}{k}d$, $e_i(*_0) = *_i$ and

$$(f_{j_1}(t_1), ..., f_{j_{n_i}}(t_{n_i})) = (e_i(f_{j'_1}^{(i)}(1-t_1)), ..., e_i(f_{j'_{n_i}}^{(i)}(1-t_{n_i}))).$$

The operad product $\alpha$ on P is defined by assigning to arc diagrams $T_0$,..., $T_m$ of type $k_i$ an arc diagram $\alpha(T_0; T_1, ..., T_m)$ of type $k_1 + \cdots + k_m$ as follows:

(1) putting each arc diagram $T_i$ into the $i$-th hole of $T_0$ by the embedding $e_i$,
(2) connecting the curve $f_{j_a}$ of $T_0$ and the curve $e_i(f_{j'_a}^{(i)})$ of $e_i(T_i)$ at all points $f_{j_a}(t_a)$ of $T_0$ $(a = 1, ..., n_i)$,
(3) rescaling sizes of disks by a homeomorphism $\psi$ of $D_0$ such that

$$\psi(\{(x,0) | x \in [-1,1]\}) = \{(x,0) | x \in [-1,1]\}$$

and $\psi(e_i(D_j \text{ of } D(k_i))) = D_{k_1+\cdots+k_{i-1}+j}$ of $D(k_1 + \cdots + k_m)$,
(4) relabeling names of disks.

The product is well-defined on P, since it does not depend on the choice of embeddings $e_i$ and the homeomorphism $\psi$ up to isotopy.

• **Units:** For each coloring $s = (\epsilon_1, ..., \epsilon_n) \in S_\pm$, the unit $I_s \in \mathrm{P}(s, \iota(s))$ is given by $I_s = \bigcup_{a=1}^{n} f_{\epsilon_a}$, where

$$f_{\epsilon_a}(t) = \begin{cases} \left(\frac{t+1}{2}\cos(\frac{(2a+1)\pi}{n} - \frac{\pi}{2}), \frac{t+1}{2}\sin(\frac{(2a+1)\pi}{n} - \frac{\pi}{2})\right) & \text{if } \epsilon_a = +, \\ \left(\frac{-t+2}{2}\cos(\frac{(2a+1)\pi}{n} - \frac{\pi}{2}), \frac{-t+2}{2}\sin(\frac{(2a+1)\pi}{n} - \frac{\pi}{2})\right) & \text{if } \epsilon_a = -. \end{cases}$$

For example,

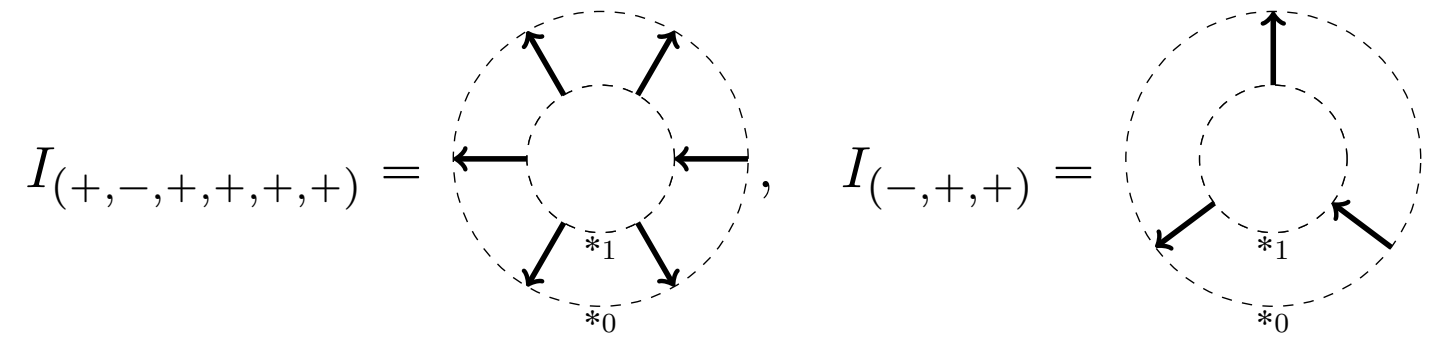

For any arc diagram $T$ of type $k$ with a coloring $(s_0; \iota(s_1), ..., \iota(s_k)) \in S_\pm^{k+1}$, we have
(1) $\alpha_1(I_{s_0}; T) \sim T$,
(2) $\alpha_k(T; I_{s_1}, ..., I_{s_k}) \sim T$.
• **Action of the braid group:** We define an action of the braid group $B_k$ on an arc diagram $T$ of type $k$ as follow. We define an action of $b_i \in B_k$ on an arc diagram $T$ of type $k$ by renaming the label $i$ with $i+1$ and the label $i+1$ with $i$ and moving $T$ by the homeomorphism $\mathfrak{b}_i$ of $D_0$ in Figure 3.

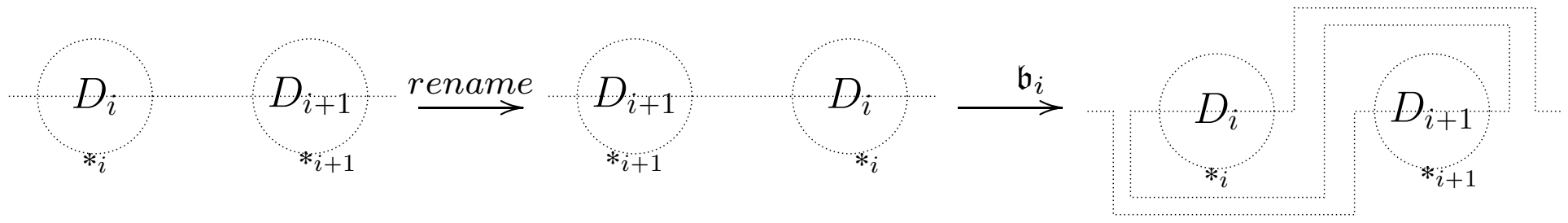

FIGURE 3. Action of braid group on arc diagram

**Lemma 4.5.** *Let $T$ be an oriented arc planar diagram of type $k$. For $i = 1, ..., k-1$, we have*
(1) $b_i b_j(T) \sim b_j b_i T$ *for* $|i - j| > 1$,
(2) $b_{i+1} b_i b_{i+1}(T) \sim b_i b_{i+1} b_i(T)$ *for* $i = 1, ..., k-2$.

*Proof.* This can be directly deduced from the definitions. □

**Theorem 4.6.** P *is a braided $(S_\pm, \iota)$-colored operad.*

*Proof.* An operad structure $\alpha$ with units $I_s$, $s \in S_\pm$, is given above. By Lemma 4.5, P is braided and it is easy to see that the other axioms of the colored operad is also satisfied. □

4.3. **The identities operad $\overline{\mathrm{P}}$.** We now introduce a cobordism presentation of the identities operad.

**Definition 4.7.** Let $T$ be an oriented $(l_1, l_2)$-arc diagram of type $k$. We define an identity cobordism of $T$, denoted by $\mathrm{id}_T$, is the image of the embedding

$$\begin{array}{rccc} \mathrm{id}_T : & \left(\bigcup_{a=1}^{j_1}[0,1] \cup \bigcup_{b=1}^{l_2} S^1\right) \times [0,1] & \longrightarrow & D(k) \times [0,1] \\ & \rotatebox{90}{$\in$} & & \rotatebox{90}{$\in$} \\ & (u,t) & \longmapsto & (T(u), t) \end{array}$$

**Definition 4.8.** Let $T$ and $T'$ be oriented $(l_1, l_2)$-arc diagrams of type $k$. Two identity cobordisms $\mathrm{id}_T$ and $\mathrm{id}_{T'}$ are called *cobordism isotopic* if there exists a cobordism isotopy, i.e. a continuous map

$$\begin{array}{rccc} \Phi: & D(k)\times[0,1]\times[0,1] & \longrightarrow & D(k)\times[0,1] \\ & \Cup & & \Cup \\ & (d,t_1,t_2) & \longmapsto & \Phi_{t_2}(d,t_1) \end{array}$$

where $\Phi_{t_2}$ for each $t_2 \in [0,1]$ is a homeomorphism of $D(k) \times [0,1]$ satisfying $\Phi_0 = \mathrm{id}_{D(k)\times[0,1]}$, $\Phi_1(\mathrm{id}_T) = \mathrm{id}_{T'}$, $\Phi_{t_2}(\{*_i\} \times [0,1]) = \{*_i\} \times [0,1]$ for $0 \leqslant i \leqslant k$.

We write $\mathrm{id}_T \overset{c}{\sim} \mathrm{id}_{T'}$ if $\mathrm{id}_T$ and $\mathrm{id}_{T'}$ are cobordism isotopic.

**Definition 4.9.** For a given $\underline{s} = (s_0; \iota(s_1), ..., \iota(s_k)) \in S_\pm^{k+1}$, we consider the set

$$\overline{\mathrm{P}}(\underline{s}) \quad = \quad \{\mathrm{id}_T \,|T \in \mathrm{P}(\underline{s})\}/ \overset{c}{\sim}$$

We define the $\bigsqcup_{k\geqslant 0} S_\pm^{k+1}$-colored set $\overline{\mathrm{P}}$ as

$$\overline{\mathrm{P}} = \{\overline{\mathrm{P}}(\underline{s})\}_{\underline{s}\in\bigsqcup_{k\geqslant 0} S_\pm^{k+1}}.$$

By definition of cobordism isotopies we have the following lemma.

**Lemma 4.10.** *Oriented $(l_1, l_2)$-arc diagrams $T$ and $T'$ are isotopic if and only these identity cobordisms $\mathrm{id}_T$ and $\mathrm{id}_{T'}$ are cobordism isotopic.*

**Proposition 4.11.** *We have an isomorphism as a set*

$$\overline{\mathrm{P}}(\underline{s}) \simeq \mathrm{P}(\underline{s}).$$

*In particular, $\overline{\mathrm{P}}$ inherits a symmetric $(S_\pm, \iota)$-colored operad structure via Theorem 4.6.*

## 5. Graded cobordism P-category $\mathcal{COB}^{gr}_{G_n}$ and $\mathcal{COB}^{gr}_{G_n/L}$

5.1. **The cobordism** P**-category** $\widehat{\mathcal{COB}}$**.** We consider embeddings $f^{(i)} : [0,1] \to D_0$ $(i = 0, 1)$ defined by

$$t \longmapsto \begin{cases} \left(\frac{1-3t}{\sqrt{2}} + \frac{3t}{2}, \frac{(-1)^i(1-3t)}{\sqrt{2}}\right) & (0 \leqslant t \leqslant \frac{1}{3}) \\ \left(\frac{3}{2} - 3t, 0\right) & (\frac{1}{3} < t \leqslant \frac{2}{3}) \\ \left(-\frac{3t-3}{2} - \frac{3t-2}{\sqrt{2}}, \frac{(-1)^i(3t-2)}{\sqrt{2}}\right) & (\frac{2}{3} < t \leqslant 1). \end{cases}$$

Denote by $X$ the union $f^{(0)} \cup f^{(1)}$ which is a trivalent diagram (See Figure 4). The coloring of $X$ is defined by $s_X = (-, -, +, +)$.

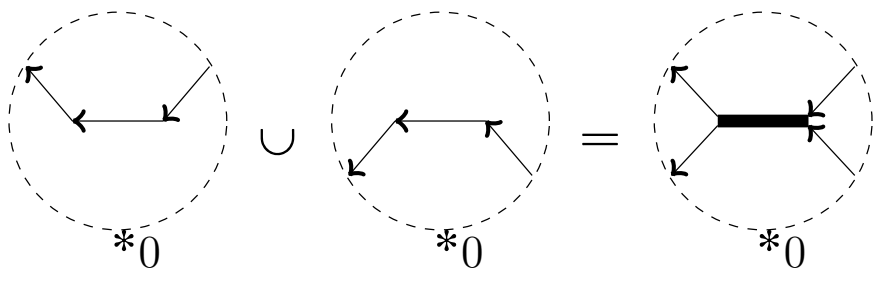

Figure 4. $f^{(0)} \cup f^{(1)} = X$

We consider the P-set $(X)_\mathrm{P}$ generated by the colored element $X$ defined as below (also see Definition 3.1).

**Definition 5.1.** We define the colored set $(X)_\mathrm{P}(s)$, $s \in S_\pm$, as

$$(X)_\mathrm{P}(s) := \{\alpha(T; X^m) \,|\, m \geqslant 0, T \in P(s; \underbrace{\iota(s_X), ..., \iota(s_X)}_{m})\}/\sim$$

and the P-set $(X)_\mathrm{P}$ as

$$(X)_\mathrm{P} := \{(X)_\mathrm{P}(s)\}_{s\in S_\pm}.$$

The formal symbols $\alpha(T; X^m)$ have the following geometric presentation.

**Example 5.2** (Geometric presentation of formal symbols in the P-set $(X)_{\mathrm{P}}$).

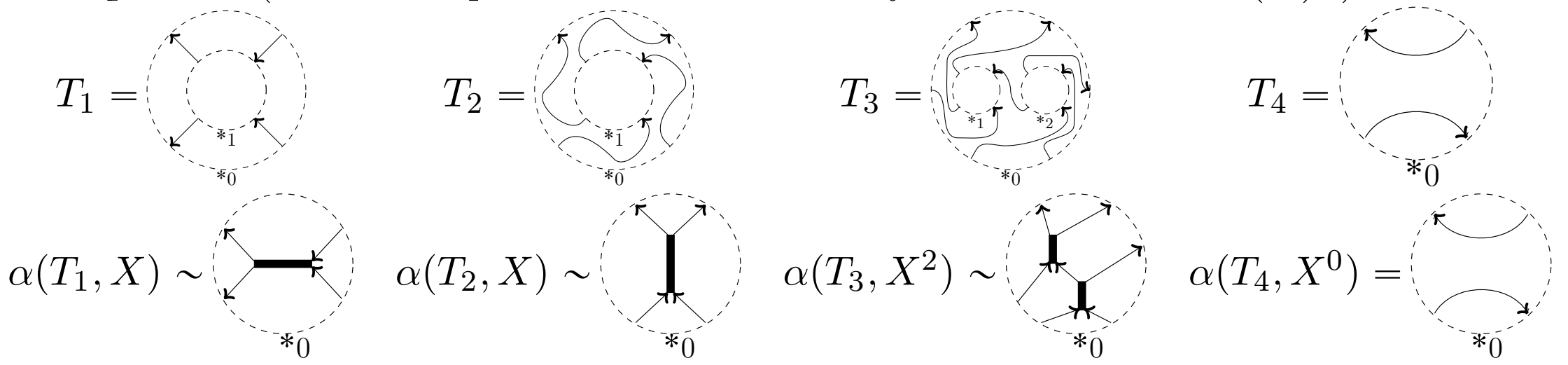

We find that $\alpha(T_1, X) \in (X)_{\mathrm{P}}(-,-,+,+)$, $\alpha(T_2, X) \in (X)_{\mathrm{P}}(-,+,+,-)$, $\alpha(T_3, X^2) \in (X)_{\mathrm{P}}(-,+,+,+,-,-)$, $\alpha(T_4, X) \in (X)_{\mathrm{P}}(+,-,+,-)$.

Let $\overline{X}$ be a cobordism of the embedding

$$\begin{array}{rccc} \mathrm{id}_X : & ([0,1] \cup [0,1]) \times [0,1] & \longrightarrow & D(k) \times [0,1] \\ & \Psi & & \Psi \\ & (u,t) & \longmapsto & (X(u), t). \end{array}$$

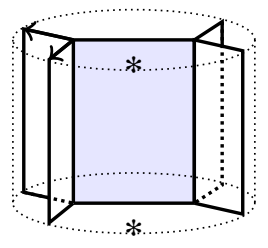

FIGURE 5. Cobordism $\mathrm{id}_X$

Define the $S_\pm$-colored set

$$(\overline{X})_{\overline{\mathrm{P}}}(s) := \{\overline{\alpha}(\mathrm{id}_T; \overline{X}^m) \mid \mathrm{id}_T \in \overline{P}(s; \underbrace{\iota(s_X), ..., \iota(s_X)}_{m}), m \geqslant 0\}/ \overset{c}{\sim}$$

and a $\overline{\mathrm{P}}$-set $(\overline{X})_{\overline{\mathrm{P}}} := \{(\overline{X})_{\overline{\mathrm{P}}}(s)\}_{s \in S_\pm}$.

**Definition 5.3.** We define $\widehat{\mathcal{COB}}$ as the P-category whose object set is $(X)_{\mathrm{P}}$ and, for objects $x$, $x' \in (X)_{\mathrm{P}}(s)$, the morphism set is

$$\widehat{\mathcal{COB}}(x, x') := \begin{cases} \{\mathrm{id}_x\} \subset (\overline{X})_{\overline{\mathrm{P}}}(s) & \text{if } x \sim x' \\ \varnothing & \text{otherwise.} \end{cases}$$

5.2. **Graded cobordism P-category $\mathcal{COB}^{gr}_{G_n}$.** We introduce now a $\mathbb{Z} \oplus \mathbb{Z}/2$-graded P-category extended by a set $G_n$ of colored morphisms. This set should be viewed as a set of foams or cobordisms, but we treat them formally.

**Definition 5.4.** The P-category $\widehat{\mathcal{COB}}[G_n]'$ is defined as the P-category $\widehat{\mathcal{COB}}$ extended by the set $G_n$ which consists of the following additional morphisms between objects of $\widehat{\mathcal{COB}}$:

- Morphisms with coloring $(+,-,+,-)$:

  $f_{s1}=$ 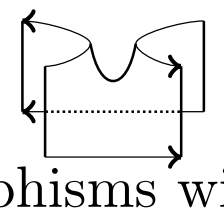 $f_{s2}=$ 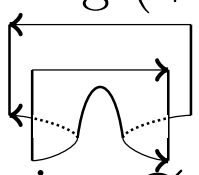
- Morphisms with coloring $\varnothing$:

  $f_{u1} =$ 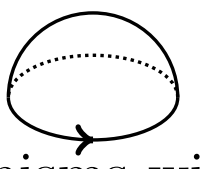 $f_{u2} =$ 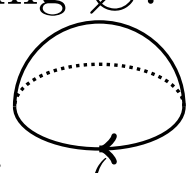 $f_{cu1} =$ 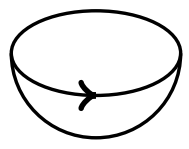 $f_{cu2} =$ 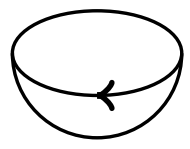
- Morphisms with coloring $(-,-,+,+)$:

  $f_z =$ 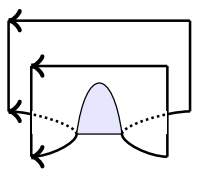 $f_{uz} =$ 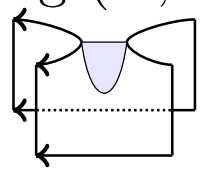

- Morphisms with coloring $(-,-,+,+)$:

  $f_{in} =$ 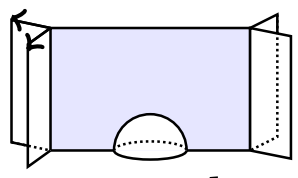 $f_{pr} =$ 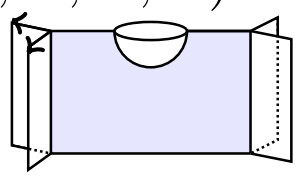

- Morphisms with coloring $(-,-,-,+,+,+)$:

  $f_{tri} =$ 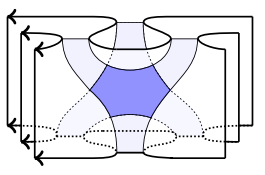

Here, the additional morphisms $f_{s1}$, $f_{s2}$, $f_{u1}$, $f_{u2}$, $f_{cu1}$ and $f_{cu2}$ describe the generators for usual cobordisms, while $f_z$, $f_{uz}$, $f_{in}$ and $f_{pr}$ are analogously for cobordisms with sheets. We need those, because we need a map from the generating object $X$ to two parallel lines. The additional morphism $f_{tri}$ is a generator composed of the sheets which are overlapped with two and three sheets; we need to include it to have maps which we need to prove the invariance under Reidemeister move 3.

**Remark 5.5.** (1) Note that the set $G_n$ does not depend on $n$. But for $n = 2$ we later set the generator $f_{tri}$ to zero, so we could also define the generating set in this case without $f_{tri}$. Our approach makes the similarities in the relations for $n = 2$ and $n \geqslant 3$ clear.
(2) The reason why we call it $G_n$ is that cobordism generators have a $\mathbb{Z} \oplus \mathbb{Z}/2\mathbb{Z}$-grading depending on $n$. See Definition 5.8.

We define cobordism isotopic on morphisms of $\widehat{\mathcal{COB}}[G_n]'$ using a cobordism isotopy which is same to Definition 4.8 in the case of $k = 0$.

**Definition 5.6** (Cobordism isotopy)**.** Let $f$ and $f'$ be cobordisms in $\widehat{\mathcal{COB}}[G_n]'$. The cobordisms $f$ and $f'$ are *cobordism isotopic* if if there exists a cobordism isotopy, i.e. a continuous map

$$\Phi : \quad \begin{array}{ccc} D_0 \times [0,1] \times [0,1] & \longrightarrow & D_0 \times [0,1] \\ \Cup & & \Cup \\ (d, t_1, t_2) & \longmapsto & \Phi_{t_2}(d, t_1) \end{array}$$

where $\Phi_{t_2}$ for each $t_2 \in [0,1]$ is a homeomorphism of $D_0 \times [0,1]$ satisfying $\Phi_0 = \mathrm{id}_{D_0 \times [0,1]}$, $\Phi_1(f) = f'$, $\Phi_{t_2}(\{*_0\} \times [0,1]) = \{*_0\} \times [0,1]$. Write $f \overset{c}{\simeq} f'$ if $f$ and $f'$ are cobordism isotopic.

**Definition 5.7.** We define $\widehat{\mathcal{COB}}[G_n]$ as the P-category consisting of the same objects to the category $\widehat{\mathcal{COB}}[G_n]'$ and morphism sets modulo cobordism isotopies.

We give the following list of some cobordism isotopies. We find that there exist homeomorphisms giving such a cobordism isotopy.

(I1) 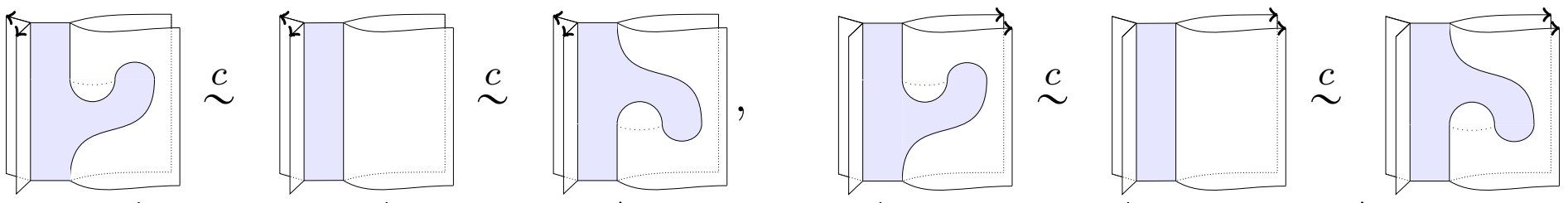

(I2) 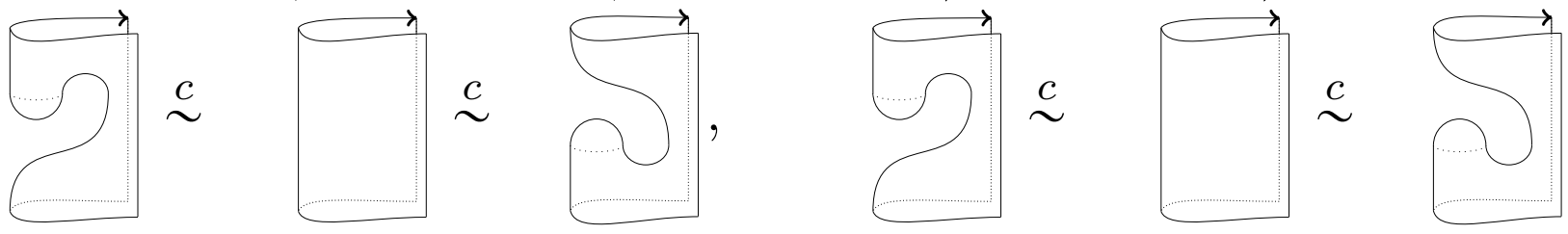

(I3) 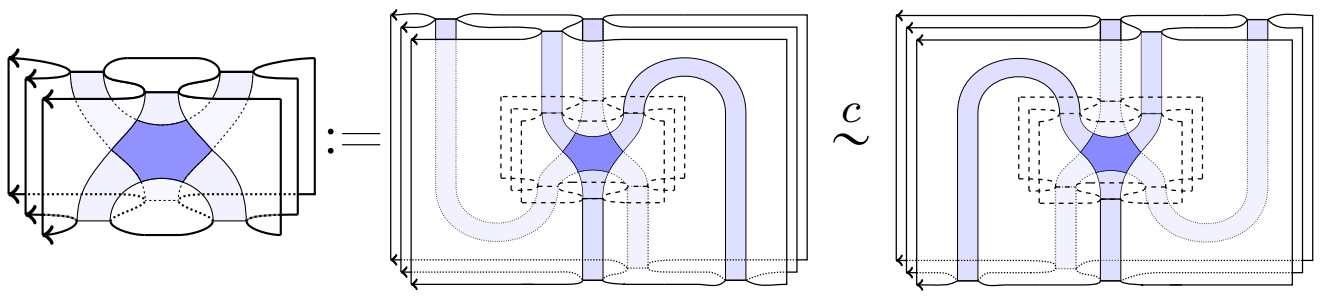

Using cobordism isotopy, the cobordism 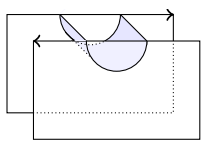 is obtained as the following composition of additional morphisms

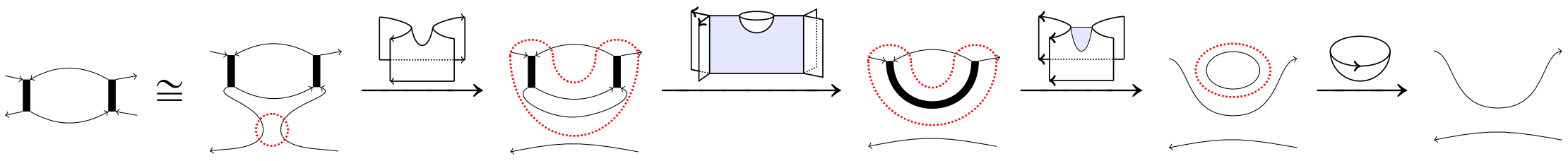

where we always apply the indicated map to the part indicated in the red circle.

Similarly, the cobordism 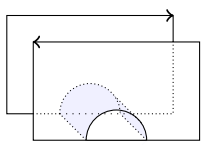 is composed from the generators of the above composition turned upside down and the cobordisms 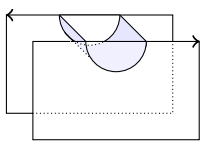 and 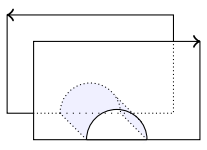 are the compositions of the above generators with reversed orientation.

By cobordism isotopy, the following different compositions of additional morphisms are isotopic.

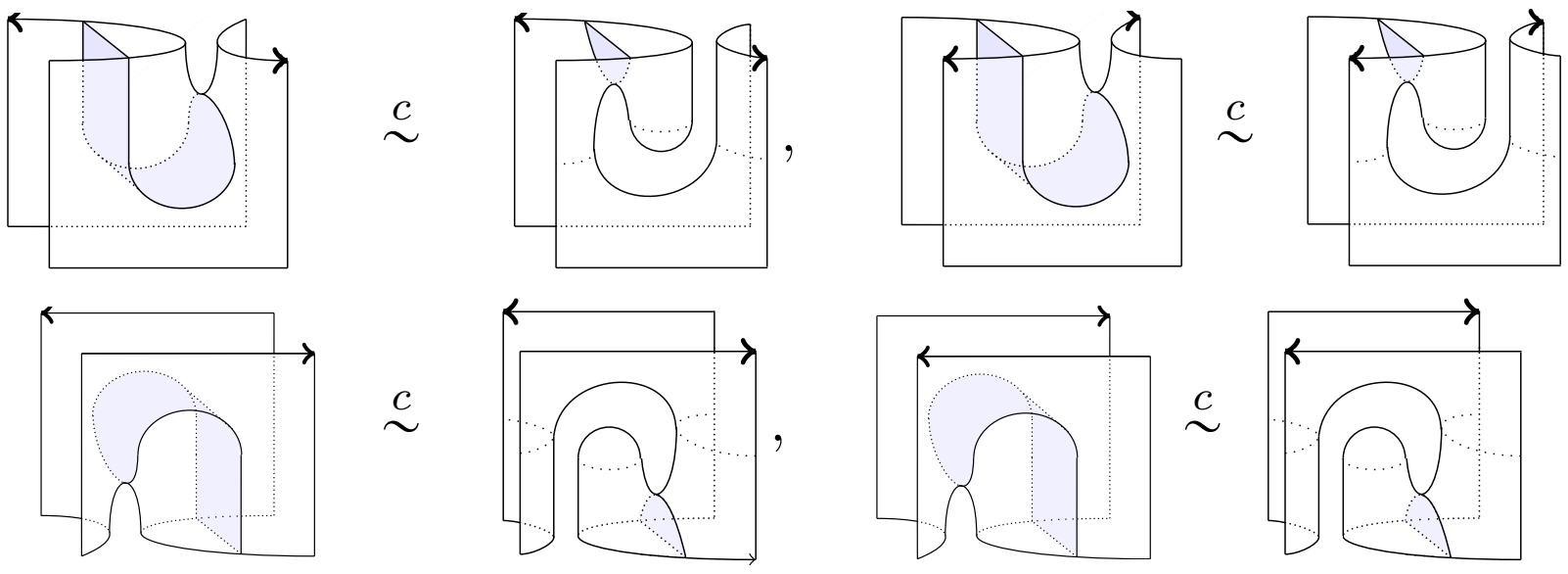

We construct the $\mathbb{Z} \oplus \mathbb{Z}/2$-graded category $\widehat{\mathcal{COB}}[G_n]^{gr}$ as follows.

**Definition 5.8.** The P-category $\widehat{\mathcal{COB}}[G_n]^{gr}$ is the $\mathbb{Z} \oplus \mathbb{Z}/2$-graded category $\widehat{\mathcal{COB}}[G_n]$ whose objects consist of $\mathbb{Z} \oplus \mathbb{Z}/2$-graded objects $x\{t, u\}$, where $x \in ob(\widehat{\mathcal{COB}}[G_n])$, $t \in \mathbb{Z}$ and $u \in \mathbb{Z}/2$.

In the graded category $\widehat{\mathcal{COB}}[G_n]$, the identity morphisms in $\widehat{\mathcal{COB}}$ and $f_{tri}$ in $G_n$ are of degree (0,0), $f_{s1}$ and $f_{s2}$ in $G_n$ are of degree $(n-1, 1)$, $f_{u1}$, $f_{u2}$, $f_{cu1}$ and $f_{cu2}$ in $G_n$ are of degree $(1-n, 1)$, $f_z$ and $f_{uz}$ in $G_n$ are of degree $(1, 0)$, and $f_{in}$ and $f_{pr}$ in $G_n$ are of degree $(-1, 0)$.

The composition of morphisms and the $\overline{\mathrm{P}}$-action on morphisms are preserving degrees, that is,

$$\deg(f \circ g) = \deg(f) + \deg(g), \quad \deg(\overline{\beta}(T, f_1, ..., f_m)) = \deg(f_1) + \cdots + \deg(f_m),$$

where $f \in \mathrm{Hom}_{\widehat{\mathcal{COB}}[G_n]^{gr}}(x', x'')(s)$, $g \in \mathrm{Hom}_{\widehat{\mathcal{COB}}[G_n]^{gr}}(x, x')(s)$, $T \in \overline{\mathrm{P}}(s_0; \iota(s_1), ..., \iota(s_k))$ and $f_i$ $(i = 1, ..., m)$ is a morphism of $\widehat{\mathcal{COB}}[G_n]^{gr}(s_i)$.

In this situation the morphism set from $x$ to $x'$ in $\widehat{\mathcal{COB}}[G_n]^{gr}(s)$ is

$$\mathrm{Hom}_{\widehat{\mathcal{COB}}[G_n]^{gr}(s)}(x, x') = \bigoplus_{z \in \mathbb{Z} \oplus \mathbb{Z}/2} \mathrm{Hom}^z(x, x'),$$

where $\mathrm{Hom}^z(x, x')$ is the set of morphisms from $x$ to $x'$ of degree $z \in \mathbb{Z} \oplus \mathbb{Z}/2$ modulo cobordism isotopies.

**Remark 5.9.** The above degree of the additional morphism of $G_n$ is derived from the degree of the corresponding morphisms of matrix factorizations. We discuss the correspondence in Section 8.

**Definition 5.10.** For any $n \geqslant 2$ the graded category $\mathcal{COB}^{gr}_{G_n}$ is defined as the scalar extended preadditive P-category $(\widehat{\mathcal{COB}}[G_n]^{gr})^{p.add}_{\mathbb{Z}[\frac{1}{n}]}$.

**Remark 5.11.** The definition of $\mathcal{COB}^{gr}_{G_n}$ involves the extension of scalars from $\mathbb{Z}$ to $\mathbb{Z}[\frac{1}{n}]$. This is done because the following Lemma 5.17 and Lemma 5.20 require that n is invertible.

Denote by 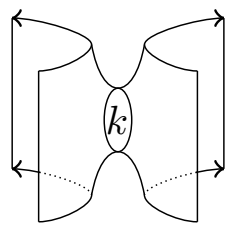 the cobordism with degree $(2k, 0)$ associated to the morphism $(f_{zip} \cdot f_{unzip})^{k+1}$, which is the cobordism with $k$ genus and $k+1$ sheets 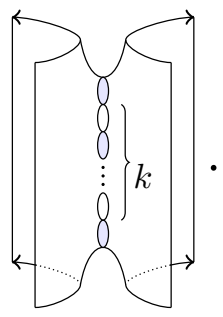. Moreover, we denote by 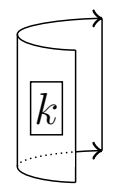 the cobordism of foams 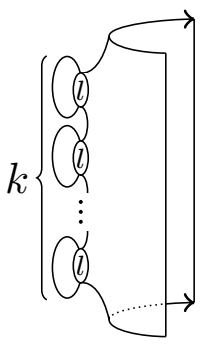 with degree $(2k, 0)$, where $l = n - 1$.

We do not include the dot on cobordism in the generating set $G_n$. The dot on a cobordism is represented by the cobordism 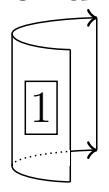.

**Remark 5.12.** We canonically regard a cobordism with a number in the box as a cobordism with foams using the orientation of the cobordism. In Figure 6, note that we find that $f_{a_1}$ and $f_{b_1}$ are isotopic, whereas $f_{a_1}$ and $f_{c_1}$ are not isotopic since rotating by 180 degrees does not preserve the base point. Therefore the cobordism $f_{b_1}$ is regarded as the cobordism $f_{a_2}$ but not to $f_{b_2}$. The cobordism $f_{c_1}$ is regarded as the cobordism $f_{c_2}$.

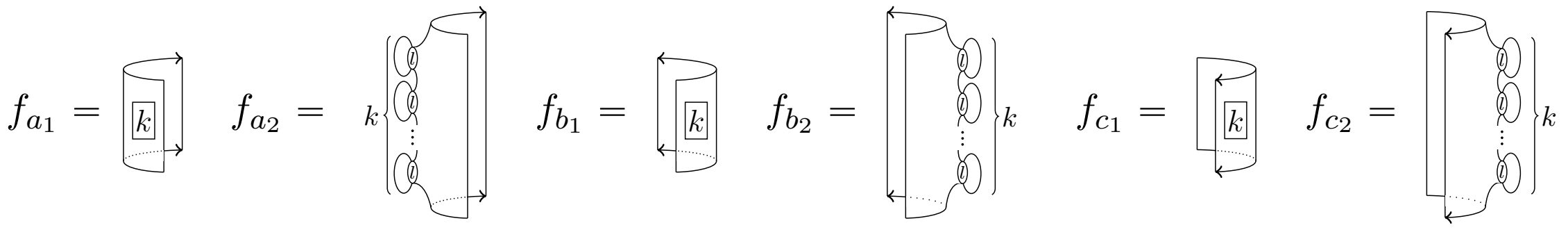

FIGURE 6. Orientation and cobordism of foams

5.3. **Graded cobordism P-category $\mathcal{COB}^{gr}_{G_n/L}$.** We define a cobordism category $\mathcal{COB}^{gr}_{G_n/L}$ associated to Khovanov–Rozansky's $sl_n$ link homology [KR08] which is analogous to the Bar-Natan's cobordism category $\mathcal{C}ob^3_l$ associated to Khovanov's $sl_2$ link homology [BN05, Kho00].

**Definition 5.13.** We define the cobordism category $\mathcal{COB}^{gr}_{G_n/L}$ as the quotient category of the preadditive category $\mathcal{COB}^{gr}_{G_n}$ consisting of the same objects, that is, $\mathrm{Ob}(\mathcal{COB}^{gr}_{G_n/L}(s)) = \mathrm{Ob}(\mathcal{COB}^{gr}_{G_n}(s))$ and the morphism sets of $\mathcal{COB}^{gr}_{G_n}$ modulo the local relations $L$ listed in Definition 5.14 below.

We define a set of local relations, denoted by $L$, in the the preadditive category $\mathcal{COB}^{gr}_{G_n}$.

**Definition 5.14.** Let $L$ be a set of the following relations in $\mathcal{COB}^{gr}_{G_n}$.

(R1) a) 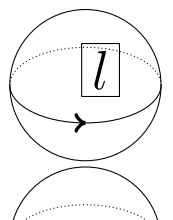 $= 0$ b) 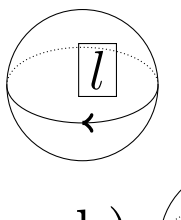 $= 0$ $(1 \leqslant l \leqslant n-2)$

(R2) a) 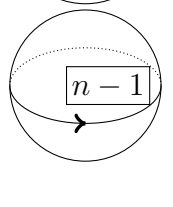 $= \frac{(-1)^{n-1}}{n}$ b) 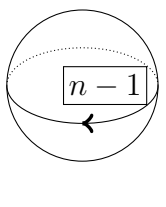 $= \frac{(-1)^{n-1}}{n}$

(R3) a) 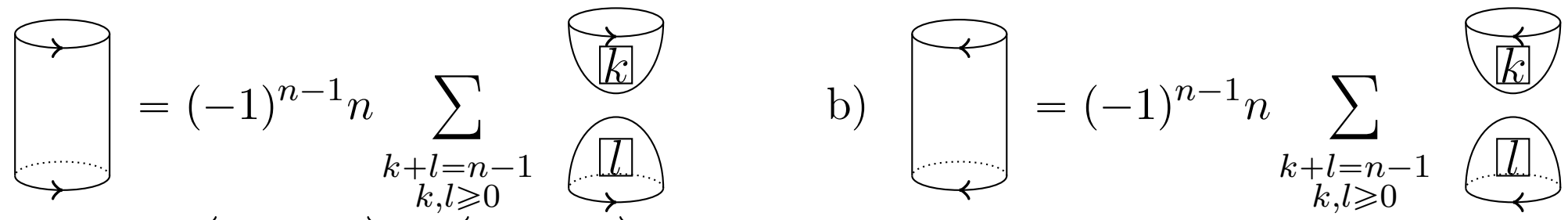
a) $= (-1)^{n-1} n \sum_{\substack{k+l=n-1\\ k,l\geqslant 0}}$  b) $= (-1)^{n-1} n \sum_{\substack{k+l=n-1\\ k,l\geqslant 0}}$

(R4) 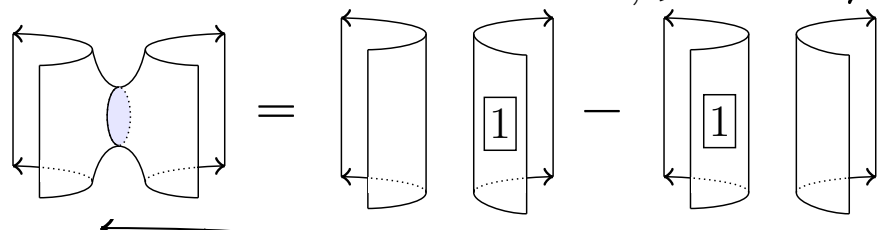
$= \quad - $

(R5) a) 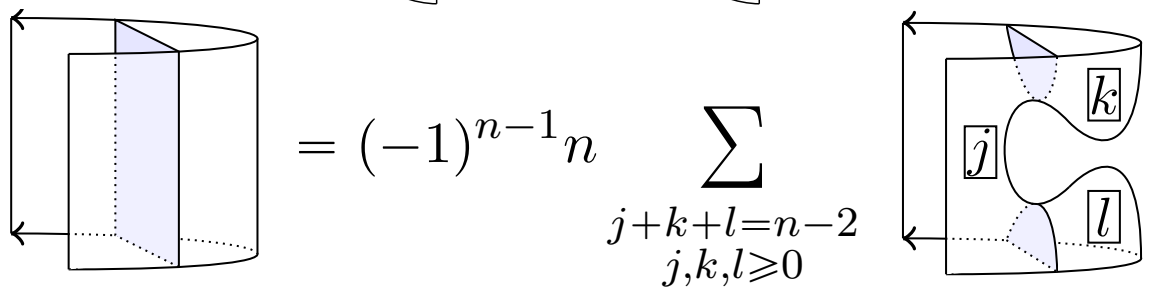
$= (-1)^{n-1} n \sum_{\substack{j+k+l=n-2\\ j,k,l\geqslant 0}}$

b) 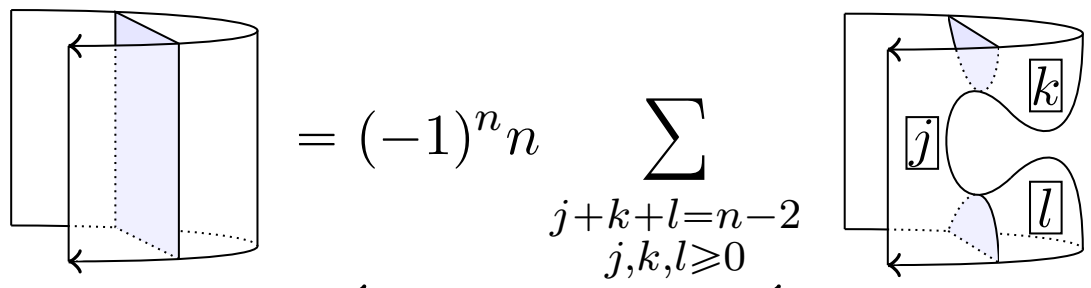
$= (-1)^{n} n \sum_{\substack{j+k+l=n-2\\ j,k,l\geqslant 0}}$

(R6) 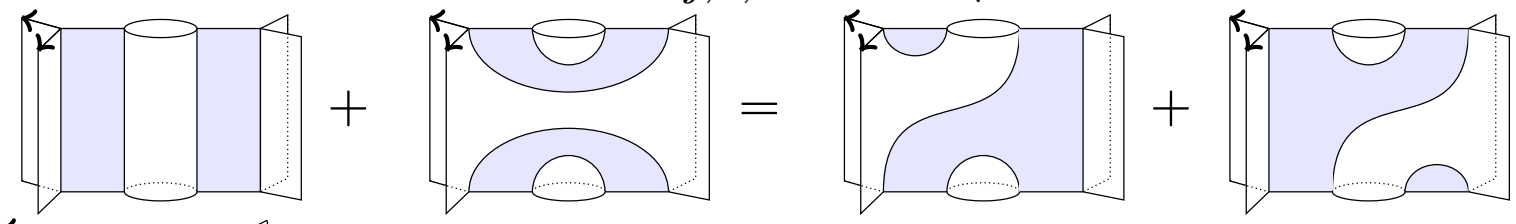
$+ \quad = \quad +$

(R7) 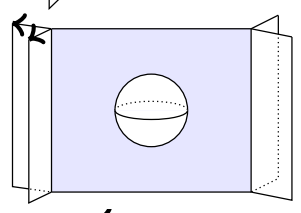 $= 0$

(R8) a) 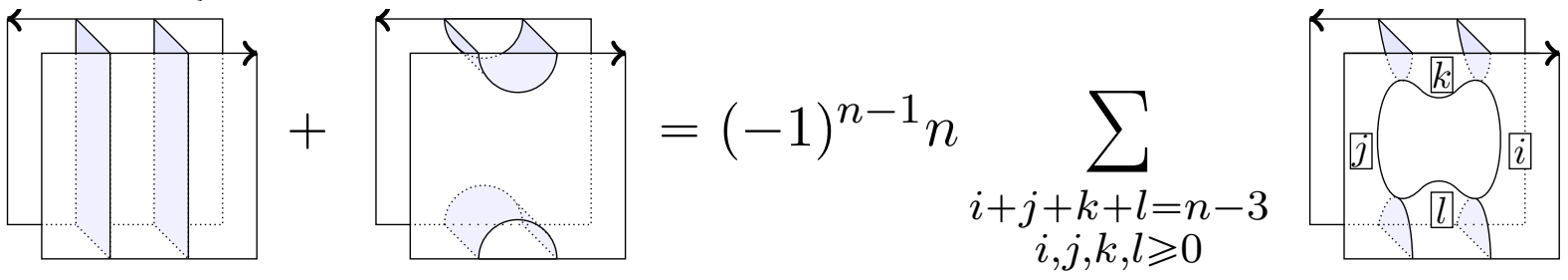
$+ \quad = (-1)^{n-1} n \sum_{\substack{i+j+k+l=n-3\\ i,j,k,l\geqslant 0}}$

b) 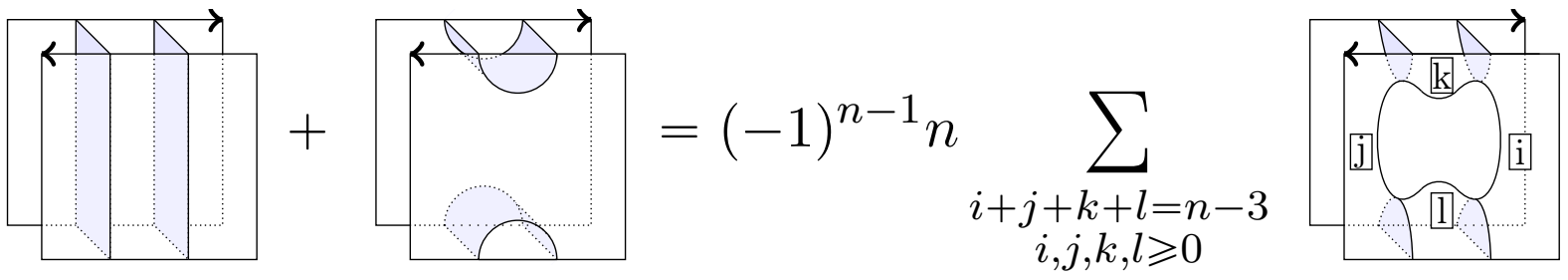
$+ \quad = (-1)^{n-1} n \sum_{\substack{i+j+k+l=n-3\\ i,j,k,l\geqslant 0}}$

(R9) a) 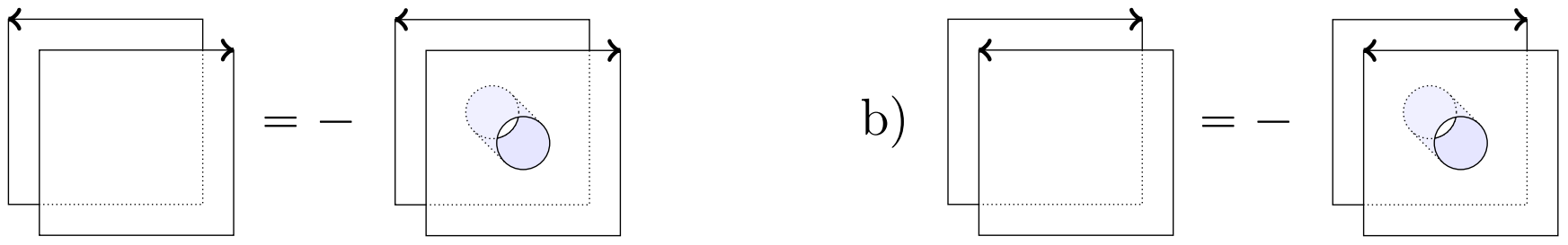
a) $= -$  b) $= -$

(R10) 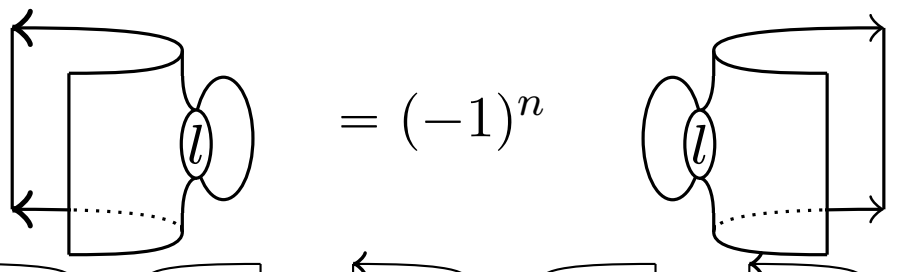
$= (-1)^{n}$ where $l = n-1$

(R11) 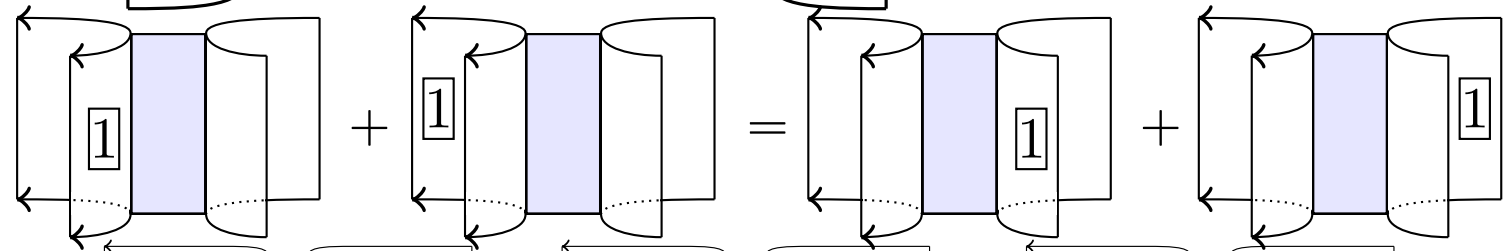
$+ \quad = \quad +$

(R12) a) 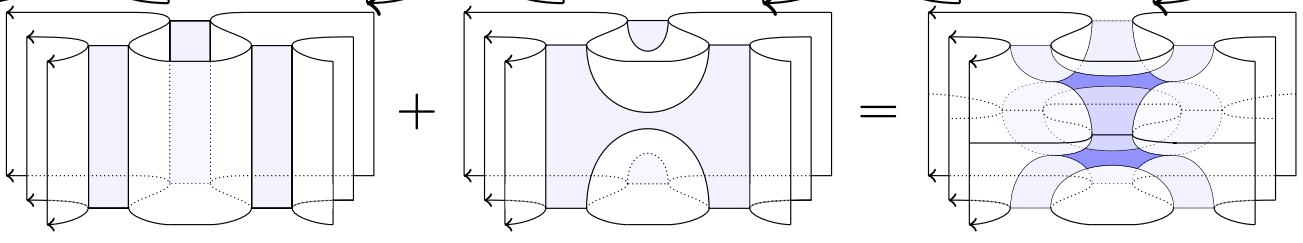
$+ \quad =$

b) 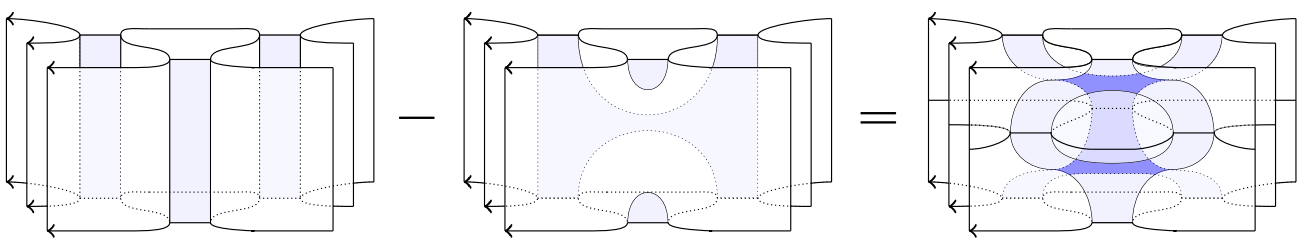
$- \quad =$

(R13) a) 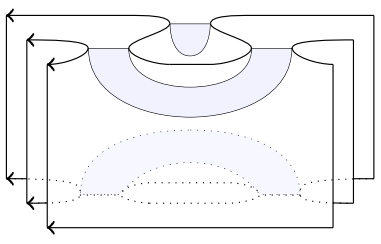 $+$ 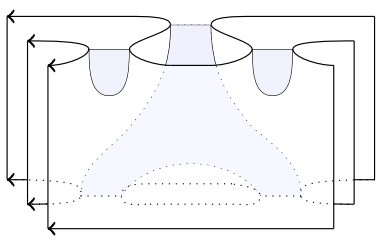 $=$ 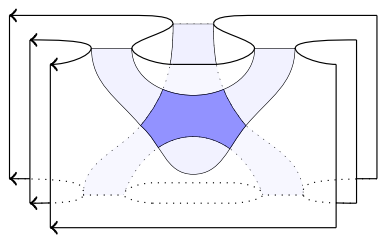

b) 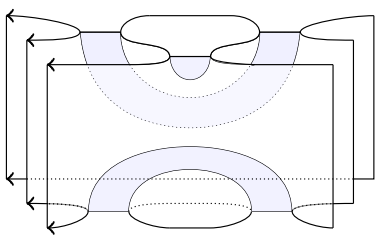 $+$ 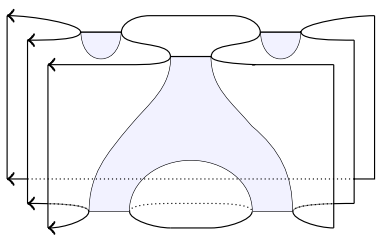 $=$ 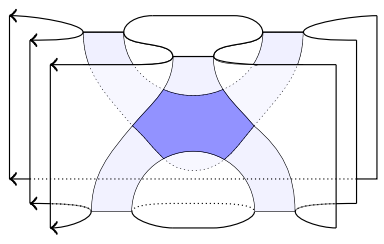

(R14) if $n = 2$: 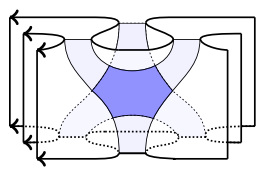 $= 0$

**Conjecture 5.15** (Consistency conjecture)**.** *The* $\mathcal{COB}^{gr}_{G_n/L}$ *is a non trivial* P*-category.*

The category $\mathcal{COB}^{gr}_{G_n}$ is a non trivial P-category but we have not still understood that the quotient category $\mathcal{COB}^{gr}_{G_n/L}$ is non trivial. For defining the relation set $L$ in $\mathcal{COB}^{gr}_{G_n}$ we use a structure of morphisms of matrix factorizations whose indeterminate sign for the saddle cobordisms is fixed. See Remark 8.19. We conjecture that our choice of the sign does not induce that the identity cobordism is equal to 0. We discuss the consistency of this cobordism category $\mathcal{COB}^{gr}_{G_n/L}$ in Section 8.

5.4. **Consequences of the relations.** We first deduce some further relations which hold in $\mathcal{COB}^{gr}_{G_n/L}$.

**Remark 5.16.** By defining 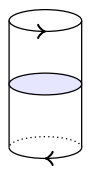 $:=$ 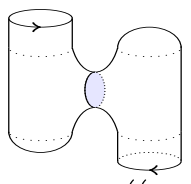 we obtain the relation 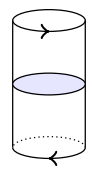 $=$ 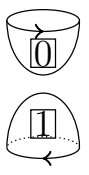 $-$ 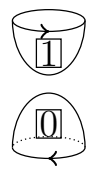 directly from (R4). Hence, we can always "cut" at a sheet independent of whether it lies horizontally or vertically.

**Lemma 5.17.** *We have the following local relations in the category* $\mathcal{COB}^{gr}_{G_n/L}$*.*

(Ra) 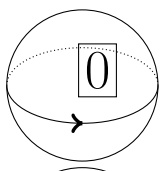 $= 0$ 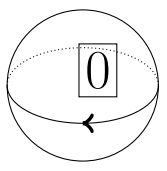 $= 0$

(Rb) 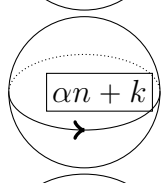 $= 0$, 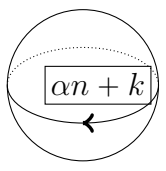 $= 0$ $(\alpha \geqslant 1,\ 0 \leqslant k \leqslant n-2)$

(Rc) 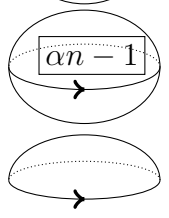 $= \frac{(-1)^{n-1}}{n}$ 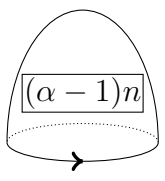 $(\alpha \geqslant 2)$.

(Rd) 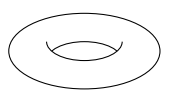 $= n$

(Re) a) $\sum\limits_{\substack{j+k=n-1\\ j,k\geqslant 0}}$ 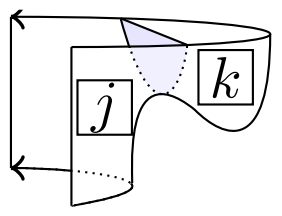 $= 0$, b) $\sum\limits_{\substack{j+k=n-1\\ j,k\geqslant 0}}$ 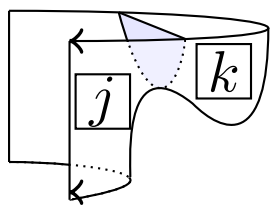 $= 0$

(Rf) a) $\sum\limits_{\substack{j+l=n-1\\ j,l\geqslant 0}}$ 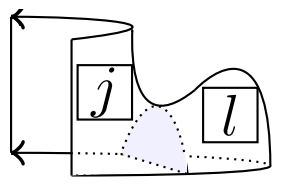 $= 0$, b) $\sum\limits_{\substack{j+l=n-1\\ j,l\geqslant 0}}$ 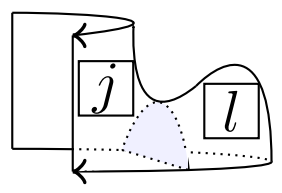 $= 0$

(Rg) 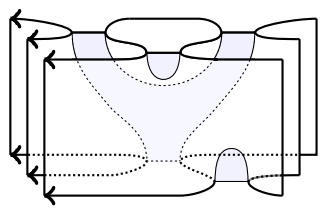 $+$ 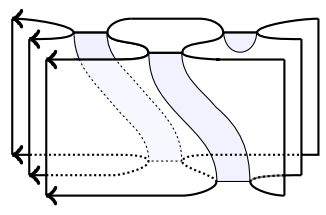 $= \begin{cases}$ 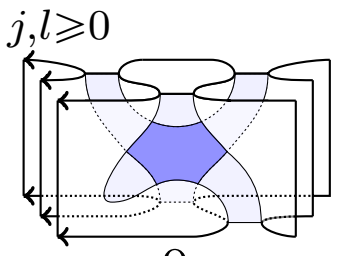 $& \text{if } n \geqslant 3 \\ 0 & \text{if } n = 2 \end{cases}$

(Rh) 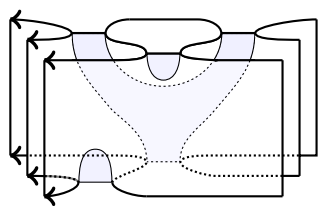 $+$ 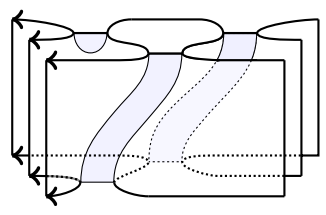 $= \begin{cases}$ 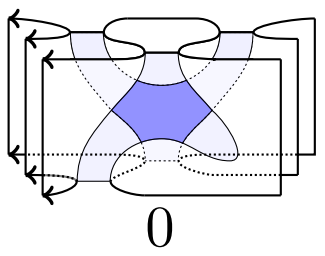 $& \text{if } n \geqslant 3 \\ 0 & \text{if } n = 2 \end{cases}$

*Proof.* Relation (Ra): We apply the cobordism 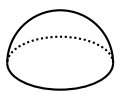 to the top of (R3) and its mirror image to the bottom. Applying (R1) only two summands remain, each of which contain one undecorated sphere and one sphere with $\boxed{n-1}$. Using (R2) we obtain that the sphere equals twice the sphere and hence must be zero.

Relation (Rb): Let $y = \alpha n + k$. Assume first that $\alpha = 1$. We apply the cobordism 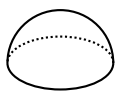 decorated with $\boxed{n-1}$ and its mirror image decorated with $\boxed{k+1}$ to the top respectively bottom of Relation (R3). If $k = 0$ the only two summands which can be nonzero contain a sphere with $\boxed{n-1}$ at the bottom, and a sphere with $\boxed{n}$ at the top, and the other way around. Using (R2) we obtain that the sphere decorated with $y$ equals twice this sphere and hence must be zero. If $k > 0$ we can argue exactly in the same way assuming the claim for smaller $k$. Now assume $\alpha \geqslant 2$ and the claim holds for smaller $y$. Then apply the cobordism 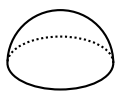 decorated with $\boxed{(\alpha-1)n+k}$ and its mirror image decorated with $\boxed{n}$ to the top respectively bottom of Relation (R3). The only summand which can be nonzero has a sphere decorated with $\boxed{2n-1}$ at the bottom, but then we have $\boxed{(\alpha-1)n+k}$ at the top, and hence we obtain zero by induction.

Relation (Rc): We now apply the cobordism 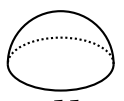 decorated with $\boxed{n}$ to the top of Relation (R3). Then from what we proved already all summands vanish except of the one containing a sphere decorated with $\boxed{2n-1}$ at the top and no decorations at the bottom. This is however equivalent to the assertion.

Relation (Rd): We obtain the torus 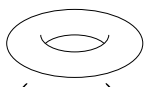, if we connect the two boundaries of the cylinder on the left-hand side of Relation (R3). By Relation (R2) we obtain the right-hand side equals $n$.

Relation (Re): We compose Relation (R5) with the cobordism for $x = n-1$ and then obtain Relation (Re) from Relation (R4). Relation (Rf) is deduced analogously.

Relation (Rg) and (Rh): By putting Relation (R13) a) into $f_{rot_1}$ in Figure 7 and using a cobordism isotopy, we obtain Relation (Rg).

Relation (Rh) is also obtained by putting Relation (R13) a) into $f_{rot_2}$ in Figure 7.

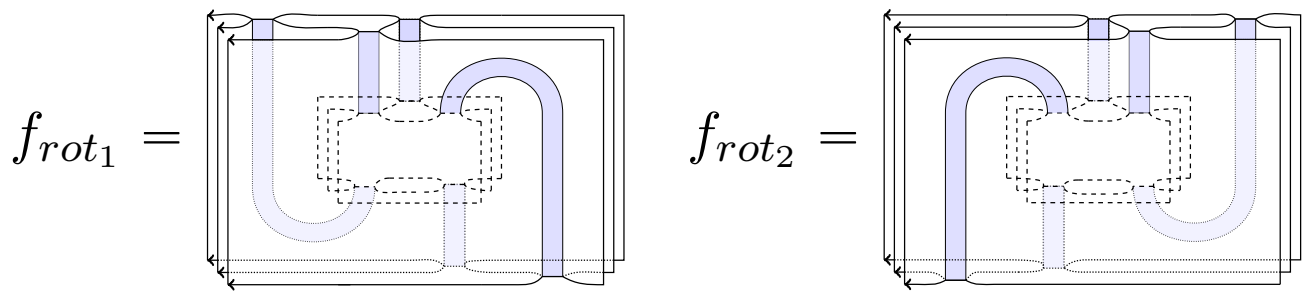

FIGURE 7. Rotation cobordisms

□

**Lemma 5.18.** *From Relation (R5) we obtain the following equalities in the category* $\mathcal{COB}^{gr}_{G_n/L}$.

(a) 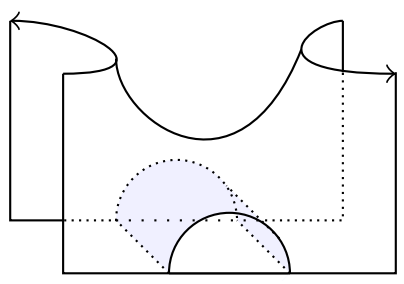 $= (-1)^n n \sum_{\substack{j+k+l=n-2 \\ j,k,l \geqslant 0}}$ 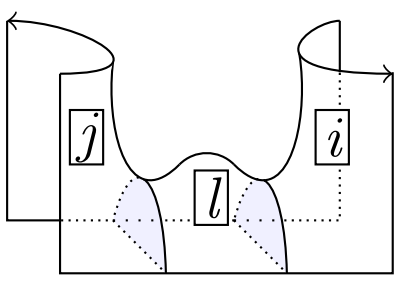

(b) 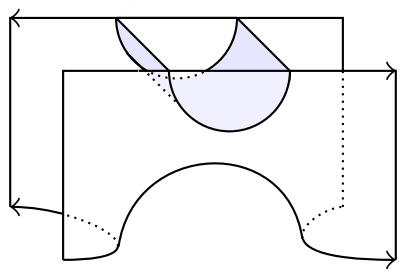 $= (-1)^n n \sum_{\substack{j+k+l=n-2 \\ j,k,l \geqslant 0}}$ 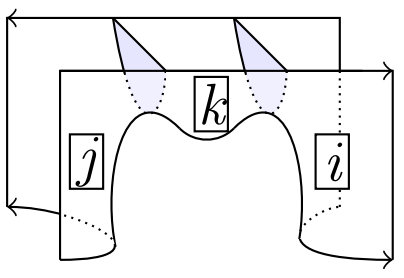

$(c)$ 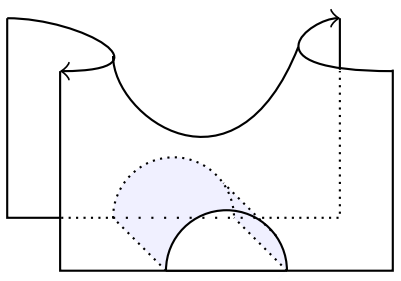 $= (-1)^{n-1} n \sum_{\substack{j+k+l=n-2 \\ j,k,l \geqslant 0}}$ 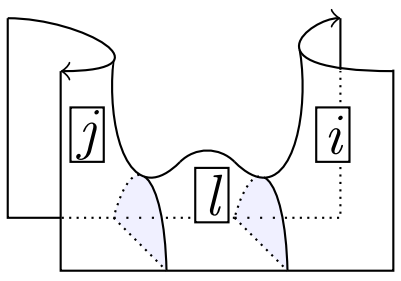

$(d)$ 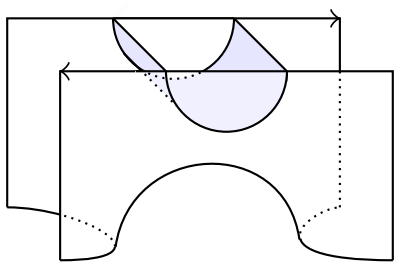 $= (-1)^{n-1} n \sum_{\substack{j+k+l=n-2 \\ j,k,l \geqslant 0}}$ 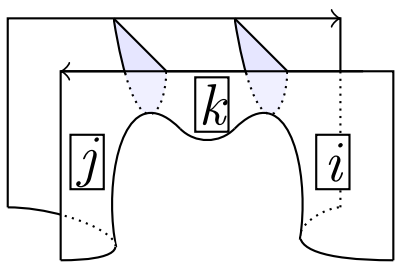

*Proof.* Put (R5) into 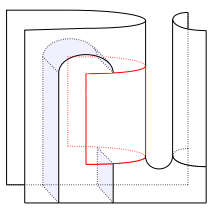 resp. its mirror image through the x-axis. □

**Lemma 5.19.** *We have the following local relations in the category $\mathcal{COB}^{gr}_{G_n/L}$.*

$(a)$ 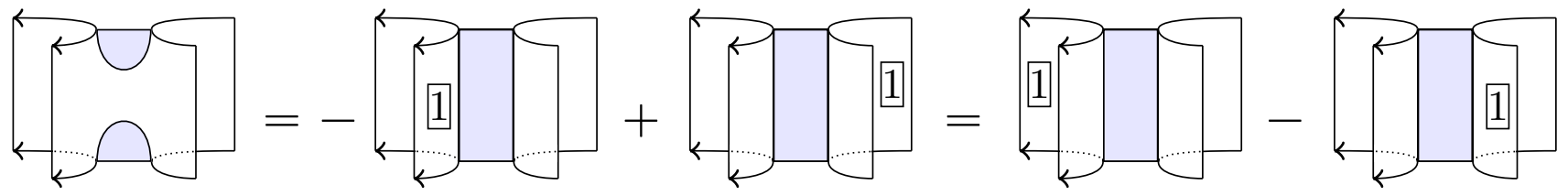

$(b)$ 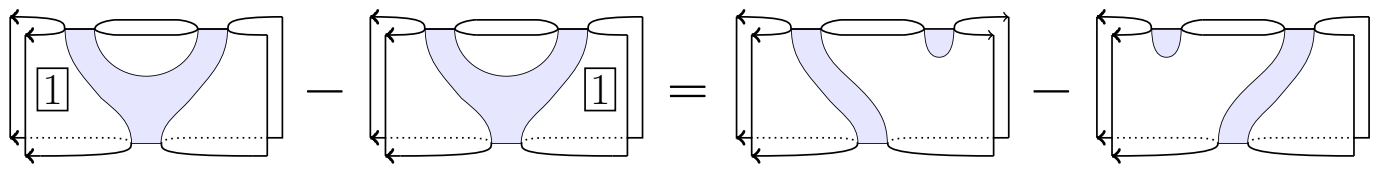

$(c)$ 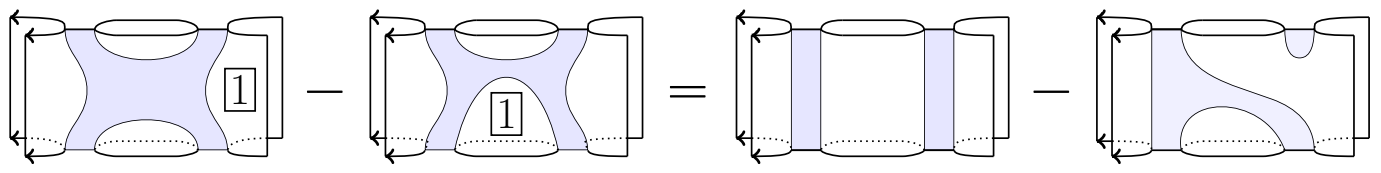

*Proof.* Putting Relation (R6) into 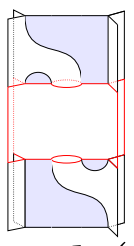 and applying Relation (R4) and (R11), we obtain the relation $(a)$. The relations $(b)$ and $(c)$ are applications of the relation $(a)$. □

Immediately, we have the following local relations by Relation (R4) and Lemma 5.19 $(b)$ and $(c)$

**Lemma 5.20.** *The following equalities hold in the category $\mathcal{COB}^{gr}_{G_n/L}$.*

$(a)$ 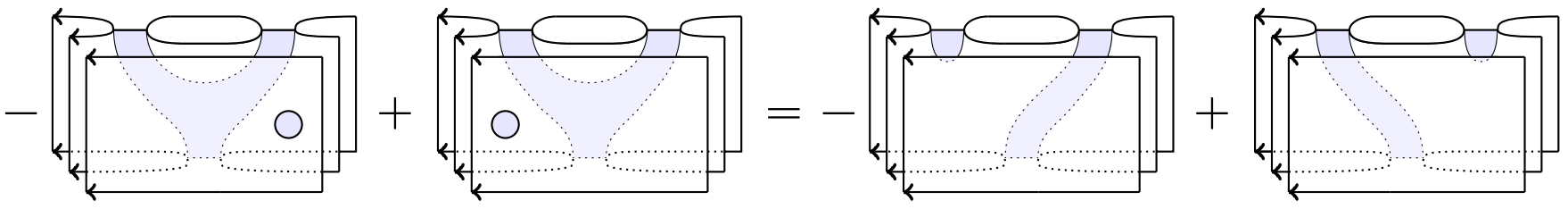

$(b)$ 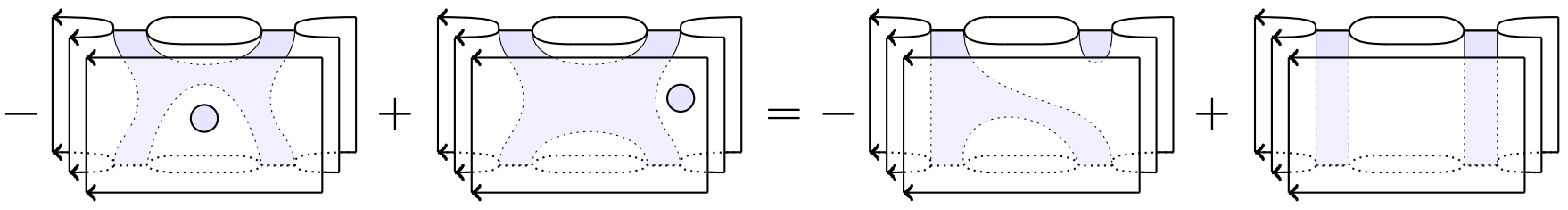

## 6. Chain complex of cobordisms for tangle

Bar-Natan defined an oriented tangle invariant with values in a cobordism category associated to Khovanov homology [BN05]. In this section, we introduce a complex of the category $\mathrm{K}^{\mathrm{b}}(\mathrm{Mat}(\mathcal{COB}^{gr}_{G_n/L}))$ associated to Khovanov–Rozansky homology for an oriented tangle. Assuming Conjecture 5.15 is true implies that the complex is an oriented tangle invariant. We discuss the invariance under the Reidemeister moves in Appendix A.

6.1. **Oriented tangle diagram and colored operad** P**.** We consider the positive and negative crossing with color $s_{cr} = (-,-,+,+)$, denoted by $Cr_+$ and $Cr_-$ (see Figure 8) and put $Cr = \{Cr_+, Cr_-\}$.

For a given oriented planar arc diagram $T \in \mathrm{P}(s; \overbrace{\iota(s_{cr}), \ldots, \iota(s_{cr})}^{m})$, $s \in S_{\pm}$, we obtain an oriented tangle diagram with $m$ crossings by putting the positive or negative crossings $Cr_+$,

$Cr_+ =$ 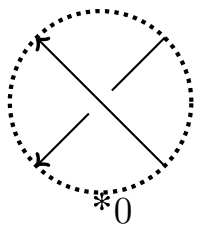

$Cr_- =$ 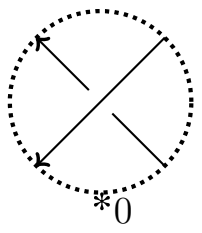

FIGURE 8. Crossings with base point and color $s_{cr} = (-,-,+,+)$

$Cr_-$ into each hole of $T$ using an embedding $e_i : D_0 \to D_i$ $(i = 1, ..., m)$. This construction is analogous to the construction of $\alpha(T; X^m)$ in Section 5. Denote by $\alpha(T; Cr_+^{m+}, Cr_-^{m-})$, $m_+ + m_- = m$, the oriented tangle diagram obtained by putting $m_+$ positive crossing into the holes $D_1$ to $D_{m_+}$ of $T$ and $m_-$ negative crossing into the holes $D_{m_++1}$ to $D_m$.

**Definition 6.1.** For $s \in S_\pm$ we define the $S_\pm$-colored set of oriented tangle diagrams

$$(Cr)_{\mathrm{P}}(s) = \{\alpha(T; Cr_+^{m+}, Cr_-^{m-}) | m_+, m_- \geqslant 0, T \in \mathrm{P}(s; \overbrace{\iota(s_{cr}), ..., \iota(s_{cr})}^{m_+ + m_-})\}/\sim,$$

where $\sim$ is the isotopic relation which is the same to Definition 4.2 in the case $k = 0$, and the P-set of oriented tangle diagrams $(Cr)_{\mathrm{P}} = \{(Cr)_{\mathrm{P}}(s)\}_{s \in S_\pm}$.

**Example 6.2.** This example is Example 3.2 in the case that $C = \mathrm{P}$, $X_+ = Cr_+$, and $X_- = Cr_-$.

$T_1 =$ 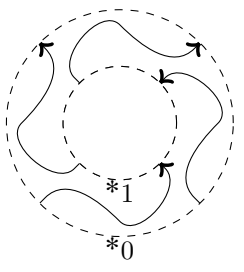

$T_2 =$ 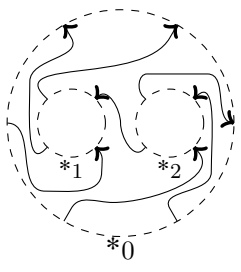

$\alpha(T_1; Cr_+^0, Cr_-) \sim$ 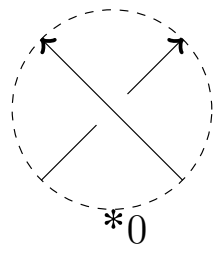

$\alpha(T_2; Cr_+^0, Cr_-^2) \sim$ 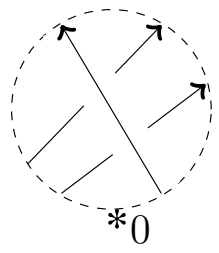

$\alpha(T_2; Cr_+, Cr_-) \sim$ 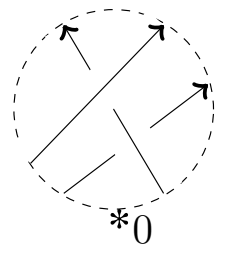

$\alpha(T_1; Cr_+^0, Cr_-)$ is an element of $(Cr)_{\mathrm{P}}(-,+,+,-)$. $\alpha(T_2; Cr_+^0, Cr_-^2)$ and $\alpha(T_2; Cr_+, Cr_-)$ are elements of $(Cr)_{\mathrm{P}}(-,+,+,+,-,-)$.

For a given oriented tangle diagram $D$ with color $s \in S_\pm$, there exists a planar arc diagram $T \in \mathrm{P}(s; \underbrace{\iota(s_{cr}), ..., \iota(s_{cr})}_{m_+ + m_-})$ such that $\alpha(T; Cr_+^{m+}, Cr_-^{m-})$ is isotopic to $D$. Conversely, for a given planar arc diagram $T \in \mathrm{P}(s; \underbrace{\iota(s_{cr}), ..., \iota(s_{cr})}_{m})$, $\alpha(T; Cr_+^{m+}, Cr_-^{m-})$, $m_+ + m_- = m$, is an oriented tangle diagram. Therefore we have the following Lemma.

**Lemma 6.3.** *There is a natural bijection between $(Cr)_{\mathrm{P}}(s)$ and the set of oriented tangle diagrams with coloring $s$.*

6.2. **From crossings and tangle diagrams to complexes of cobordisms.** We define an object in the category $\mathrm{K}^{\mathrm{b}}(\mathrm{Mat}(\mathcal{COB}^{gr}_{G_n/L}))$ for an oriented tangle diagram. We start by assigning complexes to positive and negative crossings.

**Definition 6.4.** To the positive respectively negative crossing we assign complexes in $\mathrm{K}^{\mathrm{b}}(\mathrm{Mat}(\mathcal{COB}^{gr}_{G_n/L}))(s_{cr})$ as follows:

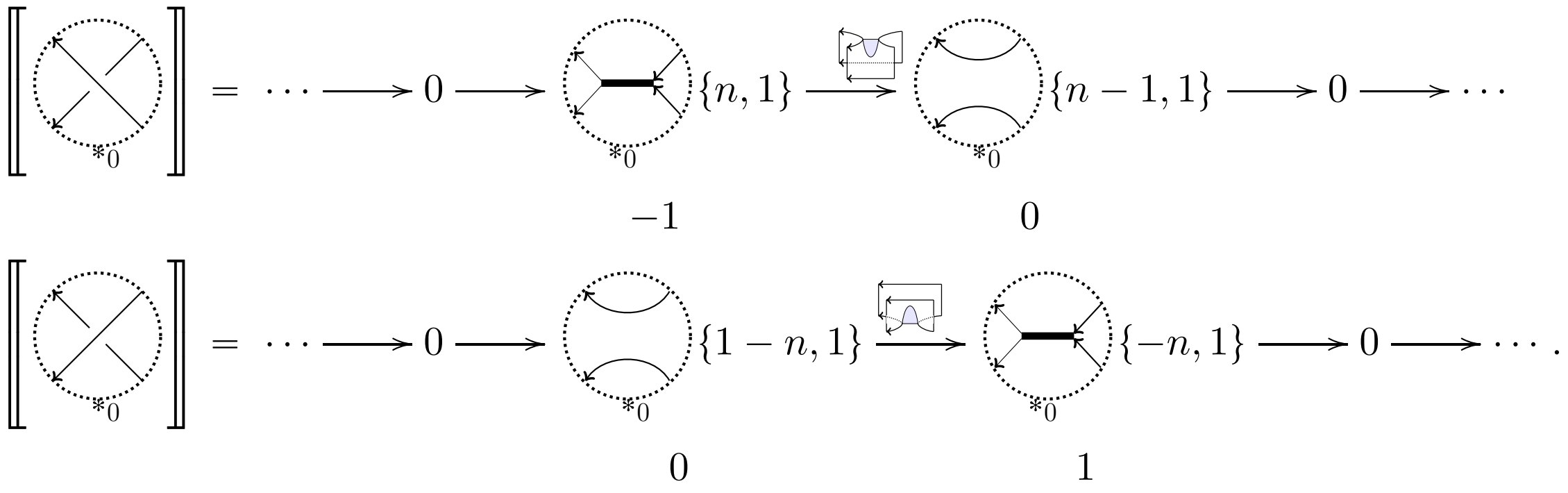

(with 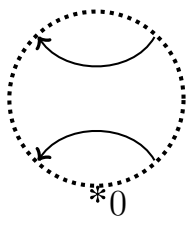 always in homological degree zero)

**Definition 6.5.** Let $T$ be a planar arc diagram in $\mathrm{P}(s; \overbrace{s_{cr}, ..., s_{cr}}^{m})$ and let $D$ be an oriented tangle diagram which is isotopic to $\alpha(T; {Cr_+}^{m_+}, {Cr_-}^{m_-})$, $m_+ + m_- = m$. To the planar diagram $D$ we define a complex $[\![D]\!]$ of $\mathrm{K}^{\mathrm{b}}(\mathrm{Mat}(\mathcal{COB}^{gr}_{G_n/L}))(s)$ as follows:

$$\beta(T, [\![Cr_+]\!]^{k_+}, [\![Cr_-]\!]^{k_-}) = \beta(T, \overbrace{[\![Cr_+]\!], ..., [\![Cr_+]\!]}^{k_+}, \overbrace{[\![Cr_-]\!], ..., [\![Cr_-]\!]}^{k_-}),$$

where $\beta$ is the structure map for complexes (see Definition 2.27).

Apriori, this definition depends on the choice of the planar arc diagram, but the complex $[\![D]\!]$ is in fact well-defined:

**Lemma 6.6.** *Let $D$ be an oriented tangle diagram. We assume that $D$ is isomorphic to $\alpha(T_1; {Cr_+}^{m_+}, {Cr_-}^{m_-})$ and $\alpha(T_2; {Cr_+}^{m_+}, {Cr_-}^{m_-})$ for $T_1$ and $T_2$ planar arc diagrams in $\mathrm{P}(s; \overbrace{s_{cr}, ..., s_{cr}}^{m})$. Then we have an isomorphism in $\mathrm{K}^{\mathrm{b}}(\mathrm{Mat}(\mathcal{COB}^{gr}_{G_n/L}))(s)$ of*

$$\beta(T_1, [\![Cr_+]\!]^{k_+}, [\![Cr_-]\!]^{k_-}) \simeq \beta(T_2, [\![Cr_+]\!]^{k_+}, [\![Cr_-]\!]^{k_-}).$$

*Proof.* We can prove this proposition by a similar way to [BN05, 2.7]. This follows by exactly the same arguments [BN05, 2.7] as in Bar Natan's original construction. □

$[\![-]\!]$ is a map from the P-sets of oriented tangle diagrams $(Cr)_{\mathrm{P}}$ to the P-set of complexes of planar arc diagrams $\mathrm{K}^{\mathrm{b}}(\mathrm{Mat}(\mathcal{COB}^{gr}_{G_n/L}))$.

We show that this map $[\![-]\!] : D \mapsto [\![D]\!]$ is invariant under the Reidemeister moves in Appendix A.

**Definition 6.7.** Let $\mathcal{A}$ be an additive category. The split Grothendieck group $K_0(\mathrm{K}^{\mathrm{b}}(\mathcal{A}))$ of the bounded homotopy category $\mathrm{K}^{\mathrm{b}}(\mathcal{A})$ is the following quotient of the free abelian group on the isomorphism classes of objects in $\mathrm{K}^{\mathrm{b}}(\mathcal{A})$:

$$K_0(\mathrm{K}^{\mathrm{b}}(\mathcal{A})) = \mathbb{Z}\langle \mathrm{Iso}(\mathrm{K}^{\mathrm{b}}(\mathcal{A}))\rangle \Big/ \left( [A_\bullet] = \sum_{i=-\infty}^{\infty} (-1)^i [A_i], [A \oplus B] = [A] + [B] \right), \quad (6.1)$$

where $\mathrm{Iso}(\mathrm{K}^{\mathrm{b}}(\mathcal{A}))$ is the isomorphism class of $\mathrm{K}^{\mathrm{b}}(\mathcal{A})$.

If in addition $\mathcal{A}$ is $\mathbb{Z} \oplus \mathbb{Z}/2\mathbb{Z}$-graded, we set $q^i[A] = [A\{i, 0\}]$ and $s[A] = [A\{0, 1\}]$, which makes $K_0(K(\mathcal{A}))$ a $\mathbb{Z}[q, q^{-1}, s]/(s^2 - 1)$-module.

From this we get immediately from Definition 6.4:

**Proposition 6.8.** *In the split Grothendieck group of the homotopy category of we have the following equality*

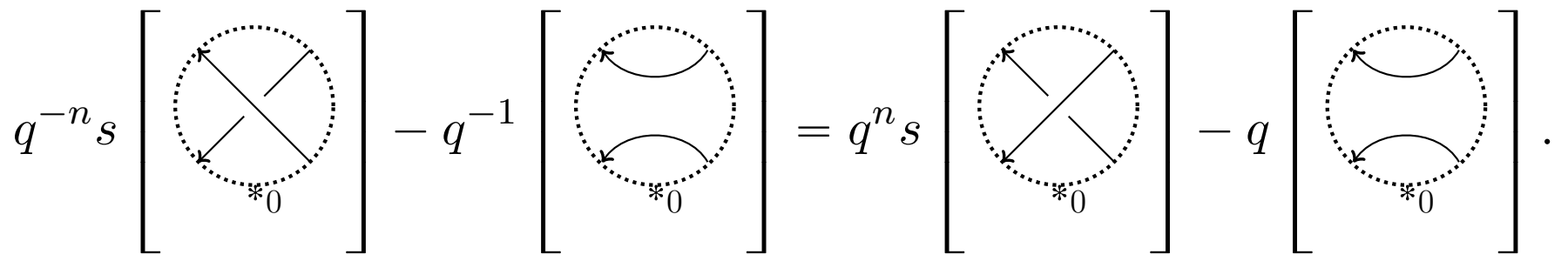

This equation in the case $s=1$ is the same to the skein relation of the polynomial $sl_n$ link invariant (specialized HOMFLY-PT polynomial). Together with the last theorem this means that we get a categorification of the $sl_n$ link invariant.

**Remark 6.9.** There are categorifications of the $(sl_n, V_n)$ link invariant which is associated to the quantum group $U_q(sl_n)$ and its vector representation $V_n$. They were obtained using the category of matrix factorizations [KR08], Lie theoretic category $\mathcal{O}$ [MS09] and the category of coherent sheaves of Grassmannian [CK08], Mackaay–Stosic–Vaz's foam category [MSV09] and, moreover, using the category of Soergel bimodule we obtain a categorification of HOMFLY-PT polynomial[Kho07]. Mackaay–Webster prove that there exists an equivalence of these categorifications [MW18]. We conjecture that our categorification using the cobordism category has the flavor of "universal" property of a categorification of the $(sl_n, V_n)$ invariant. In other words, for any categorification of $(sl_n, V_n)$ invariant

$$K_0 : \mathcal{D} \to \{(sl_n, V_n) \text{ link invariant}\}$$

there exists a full factor $\mathcal{G}$ from our categorification $\mathrm{K^b}(\mathrm{Mat}(\mathcal{COB}^{gr}_{G_n/L}))$ to the categorification $\mathcal{D}$ such that the following diagram commutes

$$\begin{array}{ccc} \mathrm{K^b}(\mathrm{Mat}(\mathcal{COB}^{gr}_{G_n/L})) & \xrightarrow{K_0} & \{(sl_n, V_n) \text{ link invariant}\} \\ & \searrow^{\mathcal{G}} \quad \mathcal{D} \quad \nearrow_{K_0} & \end{array}$$

In the rest of paper we discuss relation the cobordism category and the category of matrix factorizations.

## 7. Operad morphism from P to $\mathrm{P_{HMF}}$

In this section we define an operad morphism from the symmetric $(S_\pm, \iota)$-colored operad P to a symmetric $(S_\pm, \iota)$-colored operad $\mathrm{P_{HMF}}$ where the object sets are objects in the homotopy category of matrix factorizations HMF. The operad structure of $\mathrm{P_{HMF}}$ has been discussed by Webster [Web07].

7.1. **Matrix factorization.** We recall the notion of a matrix factorization.

Let $\Bbbk$ be a field of characteristic 0 and let $R = \Bbbk[x_1, ..., x_l]$ for $l \geqslant 1$ be the graded polynomial ring equipped with an even positive grading $R = \oplus_{i \in 2\mathbb{Z}_{\geqslant 0}} R_i$, i.e. $\deg(x_i) = 2a_i$ for some $a_i \in \mathbb{Z}_{\geqslant 1}$. A graded $R$-module $M$ is a $\mathbb{Z}$-graded vector space $\oplus_{i\in\mathbb{Z}} M^i$ with an $R$-module structure $R^j M^i \subset M^{i+j}$. For $m \in \mathbb{Z}$ define the functor $\{m\}$ which shifts the grading up by $m$, i.e.

$$(M\{m\})^i = M^{i-m}.$$

For a Laurent polynomial $f(q) = \sum a_i q^i \in \mathbb{N}_{\geqslant 0}[q, q^{-1}]$ we set

$$M^{\oplus f(q)} = \bigoplus_i (M\{i\})^{\oplus a_i}.$$

For graded $R$-modules $M$ and $N$, we denote by $\mathrm{Hom}_{R-\mathrm{gr}}(M, N)$ the vector space of grading-preserving $R$-module morphisms inside all $R$-module homomorphisms $\mathrm{Hom}_R(M, N)$ and denote by $\mathrm{Hom}_R(M, N)_i$ the vector space of homogeneous morphisms of degree $i \in \mathbb{Z}$ so that

$$\mathrm{Hom}_R(M, N)_i = \mathrm{Hom}_{R-\mathrm{gr}}(M\{i\}, N) = \mathrm{Hom}_{R-\mathrm{gr}}(M, N\{-i\}).$$

It is obvious that there is a canonical isomorphism between $\mathrm{Hom}_R(M\{j\}, N\{j\})_i$ and $\mathrm{Hom}_R(M, N)_i$ for any $j \in \mathbb{Z}$. Define $\mathrm{HOM}_R(M, N) := \oplus_{i\in\mathbb{Z}}\mathrm{Hom}_R(M, N)_i$.

**Definition 7.1.** Let $\omega$ be a homogeneous polynomial of even degree. A *(graded) matrix factorization with potential* $\omega$ is a 4-tuple $\widehat{M} = (M_0, M_1, d_{M_0}, d_{M_1})$, where $M_0$ and $M_1$ are free $R$-modules (possibly of infinite rank) and $d_{M_i} : M_i \to M_{i+1}$ is an $R$-module morphism of degree $\frac{1}{2}\deg(\omega)$ such that $d_{M_{i+1}} d_{M_i} = \omega \mathrm{id}_{M_i}$ for $i \in \mathbb{Z}/2\mathbb{Z}$.

A *grading-preserving morphism* between matrix factorizations $\widehat{M}$ and $\widehat{N}$ is defined by a pair of morphisms $f_i : M_i \to N_i \in \mathrm{Hom}_{R-\mathrm{gr}}(M_i, N_i)$, $i \in \mathbb{Z}/2\mathbb{Z}$, such that

$$d_{N_0} f_0 = f_1 d_{M_0}, \quad d_{N_1} f_1 = f_0 d_{M_1}.$$

We denote by $\mathrm{Hom}_{R-\mathrm{gr}}(\widehat{M}, \widehat{N})$ the vector space of grading-preserving morphisms from $\widehat{M}$ to $\widehat{N}$.

A *morphism of degree* $i \in \mathbb{Z}$ between $\widehat{M}$ and $\widehat{N}$ is an element in $\mathrm{Hom}_{R-\mathrm{gr}}(\widehat{M}\{i\}, \widehat{N}) = \mathrm{Hom}_{R-\mathrm{gr}}(\widehat{M}, \widehat{N}\{-i\})$, where $\widehat{M}\{i\} = (M_0\{i\}, M_1\{i\}, d_{M_0}, d_{M_1})$. Denote the $R$-module of all morphisms from $\widehat{M}$ to $\widehat{N}$ by

$$\mathrm{HOM}_{\mathrm{MF}}(\widehat{M}, \widehat{N}) = \oplus_{i\in\mathbb{Z}}\mathrm{Hom}_{R-\mathrm{gr}}(\widehat{M}\{i\}, \widehat{N}).$$

The $R$-action on this set $\mathrm{HOM}_{\mathrm{MF}}(\widehat{M}, \widehat{N})$ is $r(f_0, f_1) = (rf_0, rf_1)$ for $r \in R$ and $(f_0, f_1) \in \mathrm{HOM}_{\mathrm{MF}}(\widehat{M}, \widehat{N})$.

**Remark 7.2.** We think of a matrix factorization $\widehat{M} = (M_0, M_1, d_{M_0}, d_{M_1})$ as a 2-periodic chain of $R$-modules

$$\cdots \xrightarrow{d_{M_1}} M_0 \xrightarrow{d_{M_0}} M_1 \xrightarrow{d_{M_1}} M_0 \xrightarrow{d_{M_0}} \cdots .$$

Let $\mathrm{MF}_R^{gr}(\omega)$ be the additive category with objects all matrix factorizations with potential $\omega$ and morphism sets $\mathrm{HOM}_{\mathrm{MF}}(\widehat{M}, \widehat{N})$ for matrix factorizations $\widehat{M}$ and $\widehat{N}$.

Define a translation $\langle 1\rangle$ on a matrix factorization $\widehat{M} = (M_0, M_1, d_{M_0}, d_{M_1})$ by

$$\widehat{M}\langle 1\rangle \;:=\; (M_1, M_0, -d_{M_1}, -d_{M_0})$$

and on a morphism $\widehat{f} = (f_0, f_1)$ by $\widehat{f}\langle 1\rangle := (f_1, f_0)$.

A morphism $\widehat{f} \in \mathrm{HOM}_{\mathrm{MF}}(\widehat{M}, \widehat{N})$ is *null-homotopic* if there exists a pair of $R$-module morphisms $h_0 : M_0 \to N_1$ and $h_1 : M_1 \to N_0$ such that

$$f_0 = d_{N_1} h_0 + h_1 d_{M_0}, \quad f_1 = d_{N_0} h_1 + h_0 d_{M_1}.$$

Let $\mathrm{HMF}_R^{gr}(\omega)$ be the homotopy category of matrix factorizations, i.e. this category consists of the object set $\mathrm{Ob}(\mathrm{HMF}_R^{gr}(\omega)) = \mathrm{Ob}(\mathrm{MF}_R^{gr}(\omega))$ and the morphism set $\mathrm{HOM}_{\mathrm{HMF}}(\widehat{M}, \widehat{N})$ for $\widehat{M}$ and $\widehat{N} \in \mathrm{Ob}(\mathrm{HMF}_R^{gr}(\omega))$ is $\mathrm{HOM}_{\mathrm{MF}}(\widehat{M}, \widehat{N})$ modulo null-homotopic morphisms.

Let $\mathbb{X}_1$ and $\mathbb{X}_2$ be sets of variables with the set of common variables $\mathbb{Y}$. We consider the graded polynomial rings $R_1 = \Bbbk[\mathbb{X}_1]$, $R_2 = \Bbbk[\mathbb{X}_2]$ and $S = \Bbbk[\mathbb{Y}]$. For $\widehat{M} = (M_0, M_1, d_{M_0}, d_{M_1})$ in $\mathrm{HMF}_{R_1}^{gr}(\omega_1)$ and $\widehat{N} = (N_0, N_1, d_{N_0}, d_{N_1})$ in $\mathrm{HMF}_{R_2}^{gr}(\omega_2)$, we define a tensor product $\widehat{M} \otimes_S \widehat{N}$ in $\mathrm{HMF}_{R_1\otimes_S R_2}^{gr}(\omega_1 + \omega_2)$ by

$$\left( \begin{pmatrix} M_0 \otimes_S N_0 \\ \oplus \\ M_1 \otimes_S N_1 \end{pmatrix}, \begin{pmatrix} M_1 \otimes_S N_0 \\ \oplus \\ M_0 \otimes_S N_1 \end{pmatrix}, \begin{pmatrix} d_{M_0} & -d_{N_1} \\ d_{N_0} & d_{M_1} \end{pmatrix}, \begin{pmatrix} d_{M_1} & d_{N_1} \\ -d_{N_0} & d_{M_0} \end{pmatrix} \right).$$

For homogeneous polynomials $p$, $q \in R$ and an $R$-module $M$, the Koszul matrix factorization $K(p; q)_M$ with potential $pq$ is defined by

$$K(p; q)_M \;:=\; (M, M\{\frac{1}{2}(\deg(q) - \deg(p))\}, p, q).$$

For sequences $\mathbf{p} = (p_1, p_2, ..., p_k)$, $\mathbf{q} = (q_1, q_2, ..., q_k)$ of homogeneous polynomials in $R$ and an $R$-module $M$, a matrix factorization $K(\mathbf{p};\mathbf{q})_M$ with potential $\sum_{i=1}^{k} p_i q_i$ is defined by

$$K(\mathbf{p};\mathbf{q})_M = \bigotimes_{i=1}^{k} K(p_i; q_i)_R \otimes_R (M, 0, 0, 0).$$

7.2. **Matrix factorization of an arc diagram.** Fix $n \geqslant 2$. Khovanov–Rozansky introduced an $sl_n$ link homology whose Euler characteristic is the $sl_n$ link invariant [KR08]. We recall matrix factorizations of planar arc diagrams in Figure 9. The variables $x_1$ and $x_2$ of $T_f$ in Figure 9 have degree 2 and the variable $y$ of $T_g$ also has degree 2.

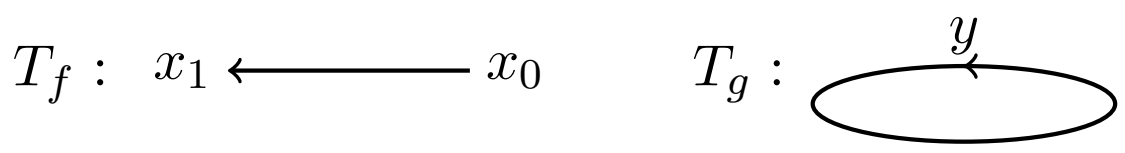

FIGURE 9. Arc diagrams $T_f$ and $T_g$

Let $h^k(x_1, x_2, ..., x_l)$ be the $k$-th complete symmetric function in variables $x_1, x_2, ..., x_l$ and let $p(e_1, e_2)$ be the polynomial associated to a representation of $x^{n+1} + y^{n+1}$ by $e_1 = x + y$ and $e_2 = xy$. In other words, $p(x + y, xy) = x^{n+1} + y^{n+1}$. For convenience, put $\overline{h}^k(x_1, x_2, ..., x_l) = \frac{n}{n+1} h^k(x_1, x_2, ..., x_l)$.

**Definition 7.3.** For the empty diagram $\varnothing$ we define the matrix factorization by

$$\widehat{\varnothing} := (\Bbbk, 0, 0, 0). \tag{7.1}$$

For the oriented diagram $T_f$ we define the matrix factorization by

$$\widehat{T_f}(x_1, x_2) := (\Bbbk[x_1, x_2], \Bbbk[x_1, x_2]\{1 - n\}, \overline{h}^n(x_1, x_2), x_1 - x_2). \tag{7.2}$$

For the oriented diagram $T_g$ we define the matrix factorization by

$$\widehat{T_g}(y) \quad = \quad (0, \Bbbk[y]/(y^n)\{1 - n\}, 0, 0). \tag{7.3}$$

7.3. **Operad morphism from P to $\mathrm{P}_{\mathrm{HMF}}$.** In this section, we define a matrix factorization for an oriented planar arc diagram $T$ and show the set of these matrix factorizations has an $(S_\pm, \iota)$-operad structure.

Let $T$ be an oriented planar arc diagram of type $k$, i.e. $T = \bigsqcup_{i=0}^{l_1} f_i \sqcup \bigsqcup_{j=1}^{l_2} g_j$ : $[0,1]^{l_1} \times (S^1)^{l_2} \to D(k)$. The color $\underline{s}(T) = (s_0(T); s_1(T), ..., s_k(T)) \in S_\pm^{k+1}$ consists of $s_a(T) = (s^{(a,1)}, ..., s^{(a,m_a)})$ for $a = 0, ..., k$.

We define the map $\varphi$ from the set of boundary points of $T$, denoted by $\partial T = \{f_i(j) | 0 \leqslant i \leqslant l_1, j = 0, 1\}$ to the set of colors of $T$, denoted by $s(T) = \{s^{(a,b)} | 0 \leqslant a \leqslant k, 1 \leqslant b \leqslant m_a\}$ by a boundary point $f_i(j)$ mapping to the color of $T$ at $f_i(j)$.

We assign two kind of variables at each boundary of an oriented planar arc diagram $T$. We assign the variable $x_{(i,j)}$ to $f_i(j)$ in $T$, where $j = 1, ..., l_1$ and $j = 0, 1$, and assign the variable $x^{(a,b)}$ to $f_i(j)$ such that $\varphi(f_i(j)) = s^{(a,b)}$, where $0 \leqslant a \leqslant k, 1 \leqslant b \leqslant m_a$. We suppose that the variables $x^{(a,b)}$ and $x_{(i,j)}$ at $f_i(j) \in \partial T$ are equal.

We use the variables $x_{(i,j)}$ for defining a matrix factorization of an oriented planar arc diagram $T$ and we use the variables $x^{(a,b)}$ for defining an operad structure on the set of matrix factorizations of arc diagrams.

For an oriented planar arc diagram $T = \bigsqcup_{i=0}^{l_1} f_i \sqcup \bigsqcup_{j=1}^{l_2} g_j$, we assign the variable $x_{(i,j)} = x^{(a,b)}$ to all boundary points of $T$ and assign the variable $y_j$ to $g_j(1, 0)$, where $(1, 0) \in S^1$.

Let $x_a(T)$ be the set of variables $x^{(a,1)}$, ..., $x^{(a,m_a)}$ on the boundary $D_a$.

**Example 7.4.** We consider the oriented planar arc diagram in Figure 10. We have

$$\begin{aligned} x_0(T) &= \{x^{(0,1)}, x^{(0,2)}, x^{(0,3)}, x^{(0,4)}, x^{(0,5)}\} = \{x_{(4,1)}, x_{(2,0)}, x_{(1,0)}, x_{(1,1)}, x_{(3,1)}\}, \\ x_1(T) &= \{x^{(1,1)}, x^{(1,2)}\} = \{x_{(3,0)}, x_{(2,1)}\}, \\ x_2(T) &= \{x^{(2,1)}, x^{(2,2)}, x^{(2,3)}\} = \{x_{(4,0)}, x_{(5,1)}, x_{(5,0)}\}. \end{aligned}$$

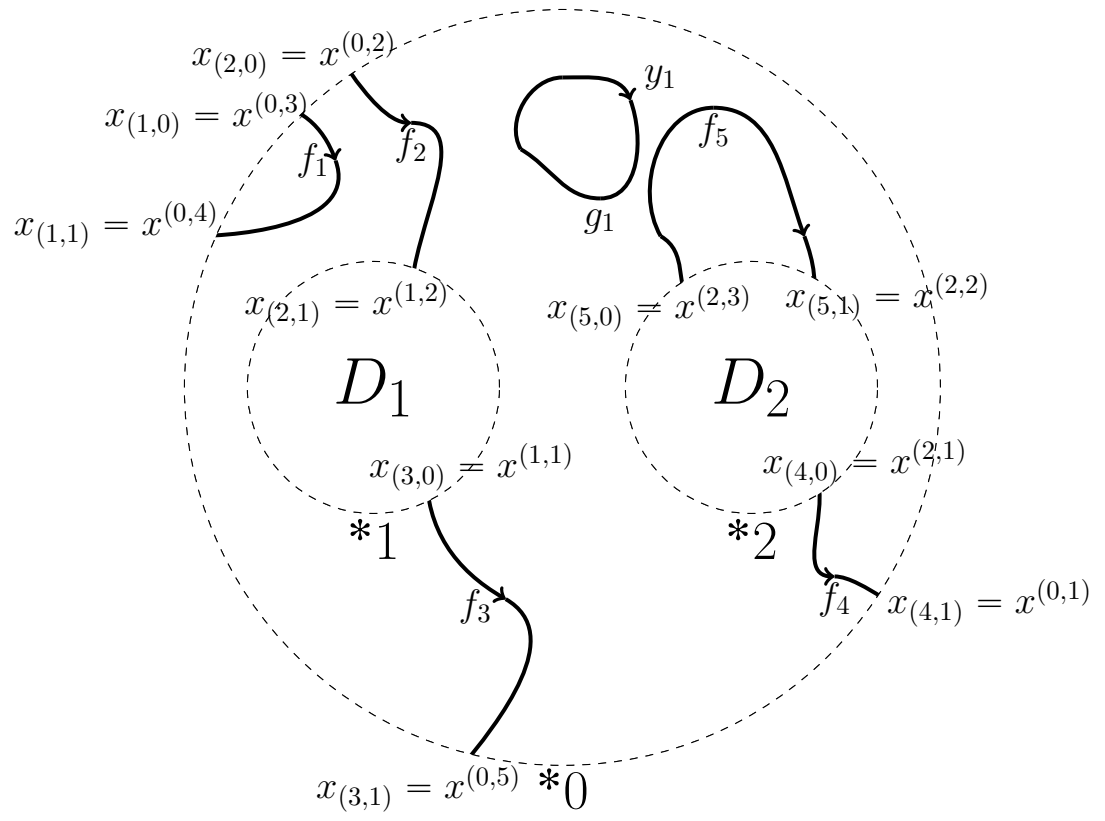

FIGURE 10. Variables associated to oriented planar arc diagram $T$

**Definition 7.5.** For an oriented planar arc diagram $T = \bigsqcup_{i=1}^{l_1} f_i \sqcup \bigsqcup_{j=1}^{l_2} g_j$ we define the matrix factorization by

$$\widehat{T} = \widehat{\bigsqcup_{i=1}^{l_1} f_i \sqcup \bigsqcup_{j=1}^{l_2} g_j} := \bigotimes_{i=1}^{l_1} \widehat{T}_{f_i}(x_{(i,1)}, x_{(i,0)}) \otimes \bigotimes_{j=1}^{l_2} \widehat{T}_{g_j}(y_j), \tag{7.4}$$

where these tensor products are defined over $\Bbbk$.

**Example 7.6.** The matrix factorization of the oriented planar arc diagram $T$ in Figure 10 is

$$\bigotimes_{i=1}^{5} \widehat{T}_{f_i}(x_{(i,1)}, x_{(i,0)}) \otimes \widehat{T}_{g_1}(y_1).$$

**Definition 7.7** (Operad $\mathrm{P}_{\mathrm{HMF}}$)**.** We define a set of matrix factorizations with color $\underline{s} \in S_{\pm}^{k+1}$ by

$$\mathrm{P}_{\mathrm{HMF}}(\underline{s}) := \{\widehat{T} | T \in \mathrm{P}(\underline{s})\} / \overset{h}{\sim}, \tag{7.5}$$

where $\overset{h}{\sim}$ is the homotopy equivalence between matrix factorizations, and a set

$$\mathrm{P}_{\mathrm{HMF}} := \bigsqcup_{k \geqslant 0} \bigsqcup_{\underline{s} \in S_{\pm}^{k+1}} \mathrm{P}_{\mathrm{HMF}}(\underline{s}). \tag{7.6}$$

We define an operad multiplication $\alpha^{\mathrm{HMF}}$ on $\mathrm{P}_{\mathrm{HMF}}$ as follows: Let $T$ be an oriented planar arc diagram $\bigsqcup_{i=1}^{l_1} f_i \sqcup \bigsqcup_{i=1}^{l_2} g_i \in \mathrm{P}(s_0, \iota(s_0^{(1)}), \dots, \iota(s_0^{(m)}))$ and let $T_j$ be an oriented planar arc diagram $\bigsqcup_{i=1}^{l_1^{(j)}} f_i^{(j)} \sqcup \bigsqcup_{i=1}^{l_2^{(j)}} g_i^{(j)} \in \mathrm{P}(s_0^{(j)}, s_1^{(j)}, \dots, s_{k_j}^{(j)})$, $j = 1, ..., m$. As before, we have matrix factorizations $\widehat{T} = \bigotimes_{i=1}^{l_1} \widehat{T}_{f_i}(x_{(i,1)}, x_{(i,0)}) \otimes \bigotimes_{i=1}^{l_2} \widehat{T}_{g_i}(y_i)$ and $\widehat{T}_j = \bigotimes_{i=1}^{l_1^{(j)}} \widehat{T}_{f_i^{(j)}}(x_{(i,1)}^{(j)}, x_{(i,0)}^{(j)}) \otimes \bigotimes_{i=1}^{l_2^{(j)}} \widehat{T}_{g_i^{(j)}}(y_i^{(j)})$.

Let $x_a(T)$ be the set of variables $x^{(a,1)}$, ..., $x^{(a,m_a)}$ on the boundary $D_a$. We define a multiplication $\alpha^{\mathrm{HMF}}_{\underline{s}}(\widehat{T};\widehat{T}_1,\ldots,\widehat{T}_m)$ by

$$\widehat{T}\otimes_{\Bbbk[x_1(T)]} e_1(\widehat{T_1})\otimes_{\Bbbk[x_2(T)]} e_2(\widehat{T_2})\otimes\ldots\otimes_{\Bbbk[x_m(T)]} e_m(\widehat{T_m}),$$

where $e_j(\widehat{T}_j)$ is the matrix factorization obtained by replacing each variable in $x_0(T_j)$ of $\widehat{T}_j$ by variables in $x_j(T)$.

**Theorem 7.8.** *The multiplication $\alpha^{\mathrm{HMF}}$ defines an operad structure on $\mathrm{P}_{\mathrm{HMF}}$ up to isomorphism of the homotopy category of matrix factorizations.*

*Proof.* The associativity of the operad structure $\alpha^{\mathrm{HMF}}$ follows from the associativity and commutativity of tensor products of matrix factorizations. The $\widehat{I}_{\underline{s}}$ has a unit structure using a homotopy equivalence (see [KR08][Proposition 15, 17]). □

The following theorem is obtained by this proposition.

**Corollary 7.9.** $\mathrm{P}_{\mathrm{HMF}}$ *is an $(S_\pm,\iota)$-colored operad.*

**Definition 7.10.** Let $T$ be an oriented arc planar diagram of type $k$. For the generator $s_a\in\mathfrak{S}_k$, $a=1,...,k-1$, we define the matrix factorization $s_a(\widehat{T})$ by replacing the variables $x^{(a,b)}$ in $\widehat{T}$ by $x^{(a+1,b)}$, $1\leqslant b\leqslant m_a$, and $x^{(a+1,b')}$ by $x^{(a,b')}$, $1\leqslant b'\leqslant m_{a+1}$.

**Theorem 7.11.** *The map $\Phi$ from $\mathrm{P}$ to $\mathrm{P}_{\mathrm{HMF}}$ defined by $T\mapsto\widehat{T}$ is a braided operad morphism.*

*Proof.* Since there is an isomorphism $\alpha^{\mathrm{HMF}}_{\underline{s}}(\widehat{T};\widehat{T}_1,\ldots,\widehat{T}_m)\simeq\widehat{\alpha_{\underline{s}}(T;T_1,\ldots,T_m)}$, the map $\Phi$ is compatible with the operad multiplications and the braid group action. □

## 8. Relation between $\mathcal{COB}^{gr}_{G_n/L}$ and $\mathrm{HMF}^{KR_n}$

Fix $n\geqslant 2$. In this section, we discuss a relation of objects of the cobordism category $\mathcal{COB}^{gr}_{G_n/L}$ and matrix factorizations and the sign problem of morphisms of matrix factorizations. We conjecture that the category of matrix factorizations $\mathrm{HMF}^{KR_n}$ is a P-category and there exists a functor of P-categories from $\mathcal{COB}^{gr}_{G_n/L}$ to $\mathrm{HMF}^{KR_n}$.

8.1. **Matrix factorizations for objects of $\mathcal{COB}^{gr}_{G_n/L}$.** We recall the matrix factorization of a trivalent diagram defined by Khovanov–Rozansky. We define matrix factorizations associated to objects of $\mathcal{COB}^{gr}_{G_n/L}$.

The object set $(X)_{\mathrm{P}}$ of $\mathcal{COB}^{gr}_{G_n/L}$ is generated by the trivalent diagram $X$ using the operad action of P. For the trivalent diagram $X$ in Figure 4, we assign the variable $x^{(0,1)}$ to $f^{(1)}(0)$, $x^{(0,2)}$ to $f^{(0)}(0)$, $x^{(0,3)}$ to $f^{(0)}(1)$ and $x^{(0,4)}$ to $f^{(1)}(1)$ (See Figure 11). We denote by $X$ this assigned trivalent diagram of $X$. The four variables $x^{(0,i)}$ $(i=1,2,3,4)$ have degree 2.

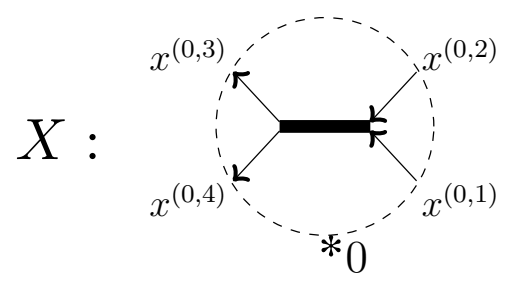

Figure 11. Assigned trivalent diagram $X$

We define the matrix factorization of the assigned trivalent diagram $X$ by

$$\begin{aligned}
&\widehat{X}(x^{(0,1)}, x^{(0,2)}, x^{(0,3)}, x^{(0,4)}) \\
:=\ &(R_{1234}, R_{1234}\{1-n\}, u_{1234}, x^{(0,3)} + x^{(0,4)} - x^{(0,1)} - x^{(0,2)}) \\
&\otimes_{R_{1234}} (R_{1234}, R_{1234}\{3-n\}, v_{1234}, x^{(0,3)}x^{(0,4)} - x^{(0,1)}x^{(0,2)})\{-1\},
\end{aligned} \tag{8.1}$$

where

$$\begin{aligned}
R_{1234} &= \Bbbk[x^{(0,1)}, x^{(0,2)}, x^{(0,3)}, x^{(0,4)}], \\
u_{1234} &= \frac{n}{n+1} \frac{p(x^{(0,3)} + x^{(0,4)}, x^{(0,3)}x^{(0,4)}) - p(x^{(0,1)} + x^{(0,2)}, x^{(0,3)}x^{(0,4)})}{x^{(0,3)} + x^{(0,4)} - x^{(0,1)} - x^{(0,2)}}, \\
v_{1234} &= \frac{n}{n+1} \frac{p(x^{(0,1)} + x^{(0,2)}, x^{(0,3)}x^{(0,4)}) - p(x^{(0,1)} + x^{(0,2)}, x^{(0,1)}x^{(0,2)})}{x^{(0,3)}x^{(0,4)} - x^{(0,1)}x^{(0,2)}}.
\end{aligned}$$

**Definition 8.1.** We define a map $\mathcal{F}$ from the set $(X)_{\mathrm{P}}$ to the set of matrix factorizations as follows. The image of the trivalent diagram $X$ with coloring $s_X$ is defined by

$$\mathcal{F}(X) := \widehat{T}_X(x^{(0,1)}, x^{(0,2)}, x^{(0,3)}, x^{(0,4)}).$$

Let $T$ be an oriented planar arc diagram of $\mathrm{P}(\underline{s})$, $\underline{s} = (s; \overbrace{\iota(s_X), \dots, \iota(s_X)}^{m})$. The matrix factorization of the formal symbol $\alpha_{\underline{s}}(T, X^m)$ is defined by

$$\mathcal{F}(\alpha_{\underline{s}}(T; X^m)) := \widehat{T} \otimes_{\Bbbk[x_1(T)]} e_1(\mathcal{F}(X)) \otimes_{\Bbbk[x_2(T)]} \cdots \otimes_{\Bbbk[x_m(T)]} e_m(\mathcal{F}(X)), \tag{8.2}$$

where $x_a(T) = \{x^{(a,1)}, x^{(a,2)}, x^{(a,3)}, x^{(a,4)}\}$ and $e_a(\mathcal{F}(X)) := \widehat{T}_X(x^{(a,1)}, x^{(a,2)}, x^{(a,3)}, x^{(a,4)})$ for $a = 1, ..., m$.

For $s = (\epsilon_1, ..., \epsilon_{2k}) \in S_\pm$, let $x(s)$ be the variables $(x^{(0,1)}, ..., x^{(0,2k)})$ and let $\omega(s)$ be the power sum $\frac{n}{n+1}(\epsilon_1(x^{(0,1)})^{n+1} + \cdots + \epsilon_{2k}(x^{(0,2k)})^{n+1})$. Let $\mathrm{HMF}^{KR_n}(s)$ for $s \in S_\pm$ by the category of matrix factorizations whose objects are composed of matrix factorizations $\mathcal{F}(Y)$ for $Y \in (X)_{\mathrm{P}}(s)$. Denote by $\mathrm{HMF}^{KR_n}$ a collection of the categories $\{\mathrm{HMF}^{KR_n}(s)\}_{s \in R}$.

**Problem 8.2.** Can we construct a $\mathrm{P}_{\mathrm{HMF}}$-category structure on $\mathrm{HMF}^{KR_n}$?

On matrix factorizations of formal symbols, we define an action of $\mathrm{P}_{\mathrm{HMF}}$ as follows.

**Definition 8.3.** The multiplication $\beta_{\underline{s}}^{\mathrm{HMF}}$ is defined by

$$\beta_{\underline{s}(T)}(\widehat{T}, \widehat{M}_1, ..., \widehat{M}_m) := \widehat{T} \otimes_{\Bbbk[x_1(T)]} e_1(\widehat{M}_1) \otimes_{\Bbbk[x_2(T)]} \cdots \otimes_{\Bbbk[x_m(T)]} e_m(\widehat{M}_m),$$

where $\widehat{T}$ is the matrix factorization associated with the oriented planar arc diagram $T$ with color $\underline{s}(T) = (s; \iota(s_1), \dots, \iota(s_m))$, $\widehat{M_i}$ is a matrix factorization in $\mathrm{Ob}(\mathrm{HMF}^{KR_n}(s_i))$ and $e_i$ is replacing the variables $x(s_i)$ by $x_i(T)$.

**Theorem 8.4.** *The collection of object sets* $\{\mathrm{Ob}(\mathrm{HMF}^{KR_n}(s_i))\}$ *is a* $\mathrm{P}_{\mathrm{HMF}}$*-set.*

*Proof.* Using homogopy equivalence of matrix factorizations, we find the multiplication $\beta_{\underline{s}}^{\mathrm{HMF}}$ as above induces a $\mathrm{P}_{\mathrm{HMF}}$-set. □

For defining a $\mathrm{P}_{\mathrm{HMF}}$-category structure on $\mathrm{HMF}^{KR_n}$, we need to construct a multiplication $\overline{\beta}_{\underline{s}}^{\mathrm{HMF}}(\widehat{T}, f_1, ..., f_k)$ on morphisms of matrix factorizations, where $\widehat{T}$ is the matrix factorization associated with the diagram $T$ with color $\underline{s}(T) = (s; \iota(s_1), \dots, \iota(s_m))$, $f_i$ is a morphism in $\mathrm{HMF}^{KR_n}(s_i)(\widehat{M_i}, \widehat{N_i})$ and $e_i$ is replacing the variables $x(s_i)$ by $x_i(T)$. This $\overline{\beta}_{\underline{s}}^{\mathrm{HMF}}$ needs to satisfy the axioms of Definition 2.20. However, it is hard to define such an action by the sign ambiguity in morphisms of matrix factorizations. We expect that this problem is solved by an operad action. We discuss this problem in Section 8.3.

8.2. **Cohomology of 2-cyclic complex.** The $\mathbb{Z}/2\mathbb{Z}$-grading on the cobordism of $\mathcal{COB}^{gr}_{G_n/L}$ comes from the $\mathbb{Z}/2\mathbb{Z}$-grading of a morphism of matrix factorizations. This structure can be seen in the cohomology of 2-cyclic complex.

We define a 2-cyclic complex $\mathrm{HOM}_R(\widehat{M},\widehat{N})$ by

$$\mathrm{HOM}^0_R(\widehat{M},\widehat{N}) \xrightarrow{d_0} \mathrm{HOM}^1_R(\widehat{M},\widehat{N}) \xrightarrow{d_1} \mathrm{HOM}^0_R(\widehat{M},\widehat{N}),$$

where

$$\begin{aligned}\mathrm{HOM}^0_R(\widehat{M},\widehat{N}) &= \mathrm{HOM}_R(M_0,N_0)\oplus\mathrm{HOM}_R(M_1,N_1),\\ \mathrm{HOM}^1_R(\widehat{M},\widehat{N}) &= \mathrm{HOM}_R(M_0,N_1)\oplus\mathrm{HOM}_R(M_1,N_0),\end{aligned}$$

and

$$d_i(f) = d_N\, f + (-1)^i f\, d_M \quad (i=0,1).$$

The cohomology of this complex is denoted by

$$\mathrm{EXT}(\widehat{M},\widehat{N}) = \mathrm{EXT}^0(\widehat{M},\widehat{N})\oplus\mathrm{EXT}^1(\widehat{M},\widehat{N}).$$

By definition, we have the following.

**Lemma 8.5.** *We have*

$$\mathrm{EXT}^0(\widehat{M},\widehat{N}) \simeq \mathrm{HOM}_{\mathrm{HMF}}(\widehat{M},\widehat{N}), \tag{8.3}$$
$$\mathrm{EXT}^1(\widehat{M},\widehat{N}) \simeq \mathrm{HOM}_{\mathrm{HMF}}(\widehat{M},\widehat{N}\langle 1\rangle). \tag{8.4}$$

We define the composition on $\mathrm{EXT}(\widehat{M},\widehat{N})$ and $\mathrm{EXT}(\widehat{L},\widehat{M})$ as follows.

**Definition 8.6.** Let $(f_0,f_1)_i$ be an element of $\mathrm{EXT}^i(\widehat{M},\widehat{N})$, $i=0,1$, and let $(g_0,g_1)_j$ be an element of $\mathrm{EXT}^j(\widehat{L},\widehat{M})$, $j=0,1$. We define a composition of $(f_0,f_1)_i$ and $(g_0,g_1)_j$ by

$$(f_0,f_1)_i(g_0,g_1)_j = \begin{cases}(f_0g_0,f_1g_1)_i & \text{if } j=0\\ (f_1g_0,f_0g_1)_{i+1} & \text{if } j=1.\end{cases} \tag{8.5}$$

By this composition, $\mathrm{EXT}(\widehat{M},\widehat{N})$ is a $\mathbb{Z}/2\mathbb{Z}$-graded $R$-module.

8.3. **Map from cobordisms of $\mathcal{COB}^{gr}_{G_n/L}$ to morphisms of** HMF**.** In Section 8.1, we defined the map $\mathcal{F}$ from objects of $\mathcal{COB}^{gr}_{G_n/L}$ to matrix factorizations. Here, we discuss a relation of cobordism generators in Section 5.4 and morphisms of matrix factorizations.

For an identity cobordism $\mathrm{id}_Y$, $Y\in(X)_{\mathrm{P}}$ we define a morphism of matrix factorization

$$\mathcal{F}(\mathrm{id}_Y) = \mathrm{id}_{\widehat{Y}}.$$

Denote by $R_{1234}$ the polynomial ring $\Bbbk[x^{(0,1)},x^{(0,2)},x^{(0,3)},x^{(0,4)}]$. We assign variables to the cobordisms $f_z$ and $f_{uz}$ as follows:

$f_z =$ 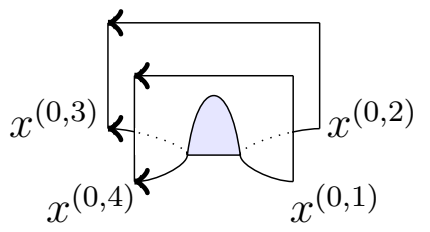

$f_{uz} =$ 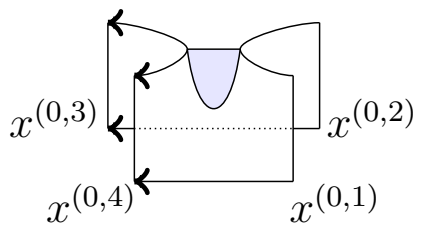

Let $\widehat{M_z}$ be the matrix factorization of the diagram on the top of $f_z$ (the bottom of $f_{uz}$):

$$\widehat{M_z} = \widehat{T_f}(x^{(0,4)},x^{(0,1)})\otimes_{\Bbbk}\widehat{T_f}(x^{(0,3)},x^{(0,2)})$$

and let $\widehat{N}_z$ be the matrix factorization of the diagram on the cobordism of $f_z$ (the top of $f_{uz}$):

$$\widehat{X}(x^{(0,1)},x^{(0,2)},x^{(0,3)},x^{(0,4)}).$$

The following lemma follows by direct calculation.

**Lemma 8.7** ([KR08]). *We have*

$$\operatorname{rank}_{R_{1234}} \mathrm{EXT}^0(\widehat{M}_z, \widehat{N}_z) = 1, \quad \operatorname{rank}_{R_{1234}} \mathrm{EXT}^1(\widehat{M}_z, \widehat{N}_z) = 0 \tag{8.6}$$

$$\operatorname{rank}_{R_{1234}} \mathrm{EXT}^0(\widehat{N}_z, \widehat{M}_z) = 1, \quad \operatorname{rank}_{R_{1234}} \mathrm{EXT}^1(\widehat{N}_z, \widehat{M}_z) = 0. \tag{8.7}$$

The image of the cobordisms $f_z$ and $f_{uz}$ by the map $\mathcal{F}$ is defined as follows.

**Definition 8.8.** Morphisms of matrix factorization $\mathcal{F}(f_z) \in \mathrm{EXT}^0(\widehat{M}_z, \widehat{N}_z)$ and $\mathcal{F}(f_{uz}) \in \mathrm{EXT}^0(\widehat{N}_z, \widehat{M}_z)$ are defined by

$$\mathcal{F}(f_z) = \left( \begin{pmatrix} -x^{(0,3)} + x^{(0,1)} & 0 \\ \alpha^{(n-1)} & 1 \end{pmatrix}, \begin{pmatrix} x^{(0,1)} & -x^{(0,3)} \\ -1 & 1 \end{pmatrix} \right) \tag{8.8}$$

$$\mathcal{F}(f_{uz}) = \left( \begin{pmatrix} 1 & 0 \\ -\alpha^{(n-1)} & -x^{(0,3)} + x^{(0,1)} \end{pmatrix}, \begin{pmatrix} 1 & x^{(0,3)} \\ 1 & x^{(0,1)} \end{pmatrix} \right), \tag{8.9}$$

where $\alpha^{(n-1)} = \frac{(x^{(0,3)} - x^{(0,1)})v - \overline{h}^n(x^{(0,4)}, x^{(0,1)}) + \overline{h}^n(x^{(0,3)}, x^{(0,2)})}{x^{(0,4)} + x^{(0,3)} - x^{(0,1)} - x^{(0,2)}}$.
The $\mathbb{Z} \oplus \mathbb{Z}/2\mathbb{Z}$-grading of these morphisms is $(1, 0)$.

We assign the variable $y$ to the cobordisms $f_{u1}$, $f_{u2}$, $f_{cu1}$ and $f_{cu2}$ as follows:

$f_{u1} =$ 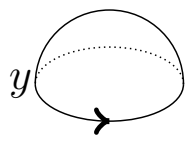

$f_{u2} =$ 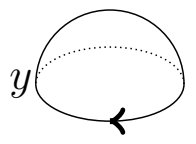

$f_{cu1} =$ 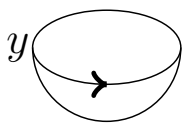

$f_{cu2} =$ 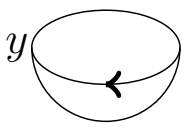

The following lemma follows by direct calculation.

**Lemma 8.9** ([KR08]). *We have*

$$\dim_{\Bbbk} \mathrm{EXT}^0(\widehat{\varnothing}, \widehat{T}_g(y)) = 0, \quad \dim_{\Bbbk} \mathrm{EXT}^1(\widehat{\varnothing}, \widehat{T}_g(y)) = n, \tag{8.10}$$

$$\dim_{\Bbbk} \mathrm{EXT}^0(\widehat{T}_g(y), \widehat{\varnothing}) = 0, \quad \dim_{\Bbbk} \mathrm{EXT}^1(\widehat{T}_g(y), \widehat{\varnothing}) = n. \tag{8.11}$$

*Proof.* This lemma follows by the following equalities

$$\widehat{\varnothing} = (\Bbbk, 0, 0, 0), \quad \widehat{T_g}(y) = (0, \Bbbk[y]/(y^n)\{1-n\}, 0, 0).$$

□

**Definition 8.10.** Morphisms of matrix factorization $\mathcal{F}(f_{u1}), \mathcal{F}(f_{u2}) \in \mathrm{EXT}^1(\widehat{\varnothing}, \widehat{T}_g(y))$ and $\mathcal{F}(f_{cu1}), \mathcal{F}(f_{cu2}) \in \mathrm{EXT}^1(\widehat{T}_g(y), \widehat{\varnothing})$ are defined by

$$\mathcal{F}(f_{u1}) = \mathcal{F}(f_{u2}) = (1, 0), \tag{8.12}$$

$$\mathcal{F}(f_{cu1}) = \mathcal{F}(f_{cu2}) = \left(0, \frac{1}{n!}\left(\frac{d}{dy}\right)^{n-1}\right). \tag{8.13}$$

The $\mathbb{Z} \oplus \mathbb{Z}/2\mathbb{Z}$-grading of these morphisms is $(1-n, 1)$.

We assign variables to the cobordisms $f_{in}$ and $f_{pr}$ and variables as follows:

$f_{in} =$ 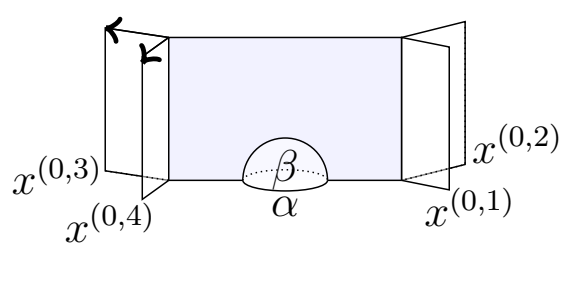

$f_{pr} =$ 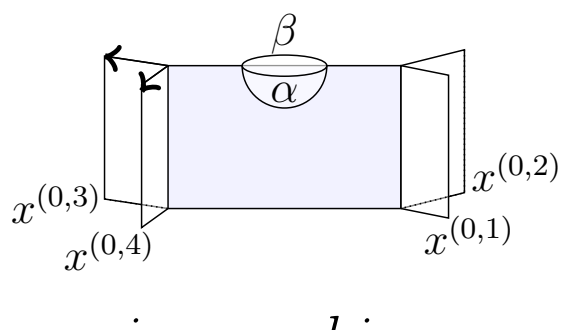

**Proposition 8.11** ([KR08, Proposition 30]). *We have an isomorphism*

$\mathcal{F}($ 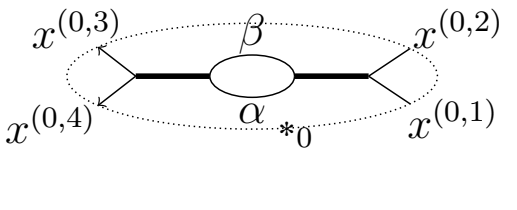

$)$

$$= \widehat{X}(\alpha, \beta, x^{(0,3)}, x^{(0,4)}) \otimes_{\Bbbk[\alpha,\beta]} \widehat{X}(x^{(0,1)}, x^{(0,2)}, \beta, \alpha)$$

$$\simeq \widehat{X}(x^{(0,1)}, x^{(0,2)}, x^{(0,3)}, x^{(0,4)})\{-1\} \oplus \alpha \widehat{X}(x^{(0,1)}, x^{(0,2)}, x^{(0,3)}, x^{(0,4)})\{-1\} \tag{8.14}$$

$\simeq \mathcal{F}($ 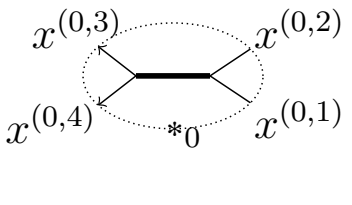

$)\{-1\} \oplus \mathcal{F}($ 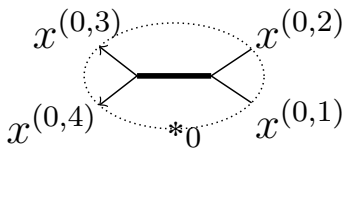

$)\{1\}.$

**Definition 8.12.** The morphism $\mathcal{F}(f_{in}) \in \mathrm{EXT}^0(\mathcal{F}($ 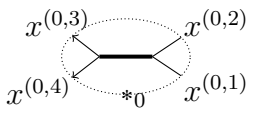 $), \mathcal{F}($ 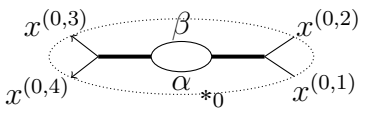 $))$ is defined by the natural inclusion

$\mathcal{F}(f_{in}) : \mathcal{F}($ 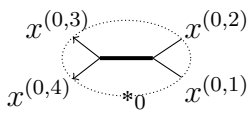 $) \longrightarrow \widehat{X}(x^{(0,1)}, x^{(0,2)}, x^{(0,3)}, x^{(0,4)})\{-1\} \subset \mathcal{F}($ 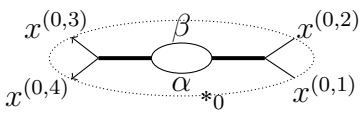 $).$

The morphism $\mathcal{F}(f_{pr}) \in \mathrm{EXT}^0(\mathcal{F}($ 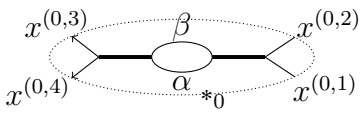 $), \mathcal{F}($ 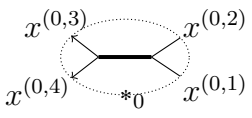 $))$ is defined by

$$\mathcal{F}(f_{pr}) = (\partial_{\alpha,\beta}, \partial_{\alpha,\beta}), \tag{8.15}$$

where $\partial_{\alpha,\beta}$ is the $R_{1234}$-linear map given by $\partial_{\alpha,\beta}(a(\alpha,\beta)) = \frac{a(\alpha,\beta)-a(\beta,\alpha)}{\alpha-\beta}$ for a polynomial $a(\alpha,\beta) \in \Bbbk[x^{(0,1)}, x^{(0,2)}, x^{(0,3)}, x^{(0,4)}, \alpha, \beta]$.

The $\mathbb{Z} \oplus \mathbb{Z}/2\mathbb{Z}$-grading of these morphisms is $(-1, 0)$.

In other words, this morphism $\mathcal{F}(f_{pr})$ is the following projection

$\mathcal{F}($ 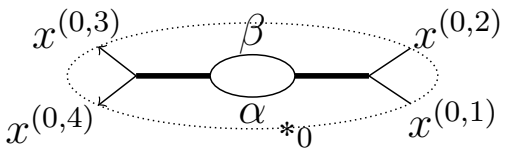 $) \supset \widehat{X}(x^{(0,1)}, x^{(0,2)}, x^{(0,3)}, x^{(0,4)})\{1\} \longrightarrow \mathcal{F}($ 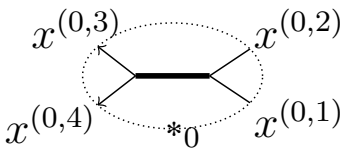 $)$

$\mathcal{F}($ 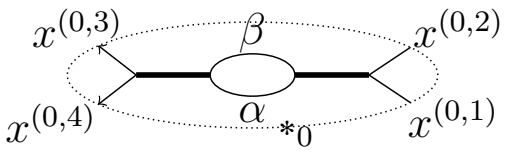 $) \supset \widehat{X}(x^{(0,1)}, x^{(0,2)}, x^{(0,3)}, x^{(0,4)})\{-1\} \longrightarrow 0.$

We assign variables to the cobordism $f_{tri}$ as follows:

$f_{tri} =$ 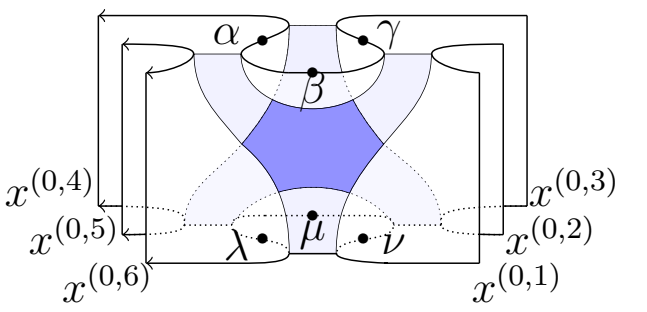

**Proposition 8.13** ([KR08], See also [Yon11, Proposition 4.12.][Wu14]). *We have isomorphisms*

$\mathcal{F}($ 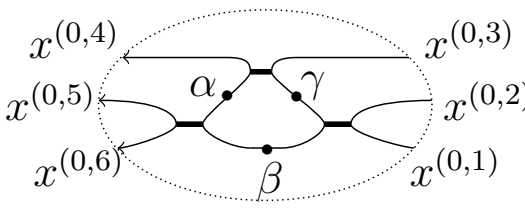 $)$

$$\begin{aligned}
&= (\widehat{X}(\beta, \alpha, x^{(0,5)}, x^{(0,6)}) \otimes_{\Bbbk[\alpha]} \widehat{X}(\gamma, x^{(0,3)}, x^{(0,4)}, \alpha)) \otimes_{\Bbbk[\beta,\gamma]} \widehat{X}(x^{(0,1)}, x^{(0,2)}, \gamma, \beta) \\
&\simeq \begin{cases} \widehat{X}(x^{(0,1)}, x^{(0,2)}, x^{(0,5)}, x^{(0,6)}) \otimes (\gamma - x^{(0,4)})\widehat{T}_f(x^{(0,4)}, x^{(0,3)})\{-2\} \oplus \widehat{Y} & \text{if } n \geqslant 3 \\ \widehat{X}(x^{(0,1)}, x^{(0,2)}, x^{(0,5)}, x^{(0,6)}) \otimes (\gamma - x^{(0,4)})\widehat{T}_f(x^{(0,4)}, x^{(0,3)})\{-2\} & \text{if } n = 2 \end{cases} \\
&\simeq \begin{cases} \widehat{X}(x^{(0,1)}, x^{(0,2)}, x^{(0,5)}, x^{(0,6)}) \otimes \widehat{T}_f(x^{(0,4)}, x^{(0,3)}) \oplus \widehat{Y} & \text{if } n \geqslant 3 \\ \widehat{X}(x^{(0,1)}, x^{(0,2)}, x^{(0,5)}, x^{(0,6)}) \otimes \widehat{T}_f(x^{(0,4)}, x^{(0,3)}) & \text{if } n = 2. \end{cases}
\end{aligned} \tag{8.16}$$

$\mathcal{F}($ 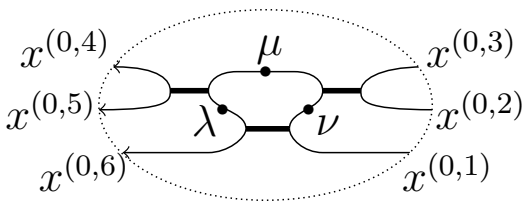 $)$

$$\begin{aligned}
&= (\widehat{X}(\lambda, \mu, x^{(0,4)}, x^{(0,5)}) \otimes_{\Bbbk[\lambda]} \widehat{X}(x^{(0,1)}, \nu, \lambda, x^{(0,6)})) \otimes_{\Bbbk[\mu,\nu]} \widehat{X}(x^{(0,2)}, x^{(0,3)}, \mu, \nu) \\
&\simeq \begin{cases} \widehat{X}(x^{(0,2)}, x^{(0,3)}, x^{(0,4)}, x^{(0,5)}) \otimes (\nu - x^{(0,6)})\widehat{T}_f(x^{(0,6)}, x^{(0,1)})\{-2\} \oplus \widehat{Y} & \text{if } n \geqslant 3 \\ \widehat{X}(x^{(0,2)}, x^{(0,3)}, x^{(0,4)}, x^{(0,5)}) \otimes (\nu - x^{(0,6)})\widehat{T}_f(x^{(0,6)}, x^{(0,1)})\{-2\} & \text{if } n = 2 \end{cases} \\
&\simeq \begin{cases} \widehat{X}(x^{(0,2)}, x^{(0,3)}, x^{(0,4)}, x^{(0,5)}) \otimes (\nu - x^{(0,6)})\widehat{T}_f(x^{(0,6)}, x^{(0,1)}) \oplus \widehat{Y} & \text{if } n \geqslant 3 \\ \widehat{X}(x^{(0,2)}, x^{(0,3)}, x^{(0,4)}, x^{(0,5)}) \otimes (\nu - x^{(0,6)})\widehat{T}_f(x^{(0,6)}, x^{(0,1)}) & \text{if } n = 2. \end{cases}
\end{aligned} \tag{8.17}$$

*where $\widehat{Y}$ is an indecomposable matrix factorization of a planar diagram with a triple edge (see Proposition 33* [KR08]*).*

**Corollary 8.14.** *We have*

$\mathrm{rank}_{R_{1,2,3,4,5,6}}\mathrm{EXT}^0(\mathcal{F}($ 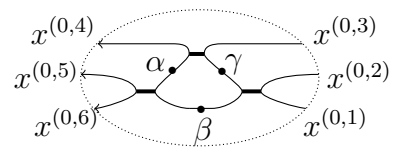 $),\mathcal{F}($ 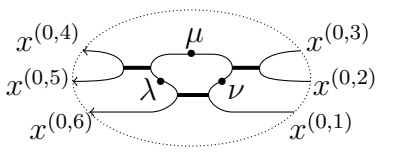 $)) = \begin{cases} 1 & \text{if } n \geqslant 3 \\ 0 & \text{if } n = 2, \end{cases}$

$\mathrm{rank}_{R_{1,2,3,4,5,6}}\mathrm{EXT}^1(\mathcal{F}($ 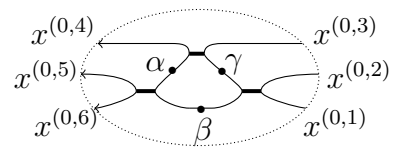 $),\mathcal{F}($ 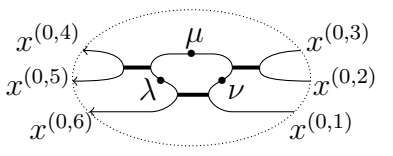 $)) = 0,$

*where* $R_{1,2,3,4,5,6} = \Bbbk[x^{(0,1)}, x^{(0,2)}, x^{(0,3)}, x^{(0,4)}, x^{(0,5)}, x^{(0,6)}]$.

When $n \geqslant 3$, the morphism $\mathcal{F}(f_{tri})$ is defined by a morphism which factors through $\widehat{Y}$:

$\mathcal{F}(f_{tri}) : \mathcal{F}($ 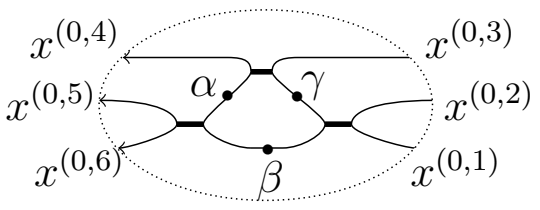 $) \longrightarrow \widehat{Y} \longrightarrow \mathcal{F}($ 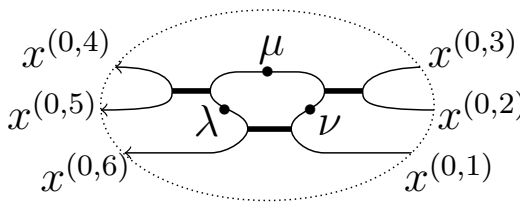 $).$

**Definition 8.15.** The morphism $\mathcal{F}(f_{tri}) \in \mathrm{EXT}^0$ is defined by

$$\mathcal{F}(f_{tri}) = \begin{cases} (\partial_{\beta,\gamma}(\beta - x^{(0,4)}), \partial_{\beta,\gamma}(\beta - x^{(0,4)})) & \text{if } n \geqslant 3 \\ 0 & \text{if } n = 2. \end{cases} \tag{8.18}$$

The $\mathbb{Z} \oplus \mathbb{Z}/2\mathbb{Z}$-grading of the morphism is $(n-1, 1)$.

**Lemma 8.16.** *When $n \geqslant 3$, the morphism $\mathcal{F}(f_{tri})$ is homotopic to*

$$(\partial_{\alpha,\beta}(\beta - x^{(0,3)}), \partial_{\alpha,\beta}(\beta - x^{(0,3)})).$$

We assign variables to the cobordisms $f_{s_1}$ and $f_{s_2}$ as follows:

$f_{s_1} =$ 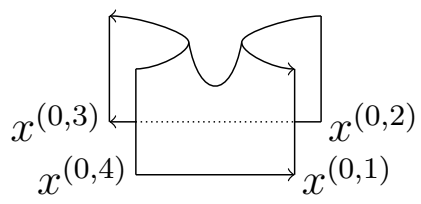 $f_{s_2} =$ 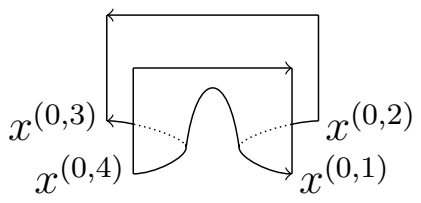

The matrix factorization of the top of the cobordism $f_{s_1}$ (the bottom of $f_{s_2}$) is

$\mathcal{F}($ 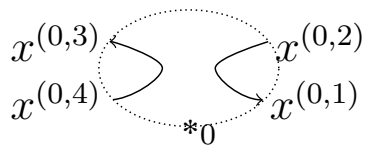 $) = \widehat{T}_f(x^{(0,3)}, x^{(0,4)}) \otimes_{\Bbbk} \widehat{T}_f(x^{(0,1)}, x^{(0,2)})$

and the matrix factorization of the bottom of the cobordism $f_{s_1}$ (the top of $f_{s_2}$) is

$\mathcal{F}($ 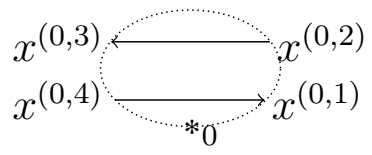 $) = \widehat{T}_f(x^{(0,3)}, x^{(0,2)}) \otimes_{\Bbbk} \widehat{T}_f(x^{(0,1)}, x^{(0,4)}).$

**Proposition 8.17.** *We have*

$\mathrm{rank}_{R_{1234}}\mathrm{EXT}^0(\mathcal{F}($ 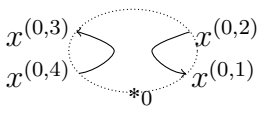 $),\mathcal{F}($ 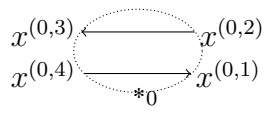 $)) = 0,$

$\mathrm{rank}_{R_{1234}}\mathrm{EXT}^1(\mathcal{F}($ 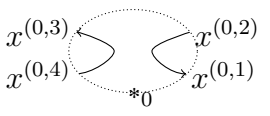 $),\mathcal{F}($ 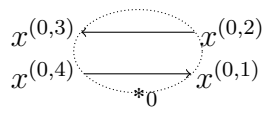 $)) = 1,$

$\mathrm{rank}_{R_{1234}}\mathrm{EXT}^0(\mathcal{F}($ 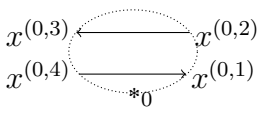 $),\mathcal{F}($ 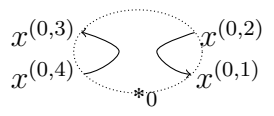 $)) = 0,$

$\mathrm{rank}_{R_{1234}}\mathrm{EXT}^1(\mathcal{F}($ 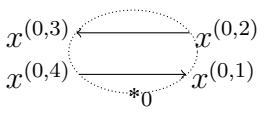 $),\mathcal{F}($ 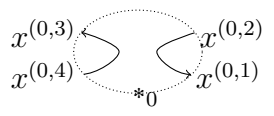 $)) = 1.$

By definition, we have the matrix factorization of the top diagram of $f_{s_1}$ (the bottom of $f_{s_2}$) is

$$\left( \begin{matrix} R_{1234} \\ \oplus \\ R_{1234}\{2-2n\} \end{matrix}, \begin{matrix} R_{1234}\{1-n\} \\ \oplus \\ R_{1234}\{1-n\} \end{matrix}, \begin{pmatrix} \overline{h}^n_{34} & x^{(0,2)} - x^{(0,1)} \\ \overline{h}^n_{12} & x^{(0,3)} - x^{(0,4)} \end{pmatrix}, \begin{pmatrix} x^{(0,3)} - x^{(0,4)} & x^{(0,1)} - x^{(0,2)} \\ -\overline{h}^n_{12} & \overline{h}^n_{34} \end{pmatrix} \right)$$

and the matrix factorization of the bottom diagram of $f_{s_1}$ (the top of $f_{s_2}$) is

$$\left( \begin{array}{c} R_{1234} \\ \oplus \\ R_{1234}\{2-2n\} \end{array}, \begin{array}{c} R_{1234}\{1-n\} \\ \oplus \\ R_{1234}\{1-n\} \end{array}, \begin{pmatrix} \overline{h}^n_{23} & x^{(0,4)} - x^{(0,1)} \\ \overline{h}^n_{14} & x^{(0,3)} - x^{(0,2)} \end{pmatrix}, \begin{pmatrix} x^{(0,3)} - x^{(0,2)} & x^{(0,1)} - x^{(0,4)} \\ -\overline{h}^n_{14} & \overline{h}^n_{23} \end{pmatrix} \right),$$

where $\overline{h}^n_{ij}$ is the $n$-th complete symmetric function $\overline{h}^n(x^{(0,i)}, x^{(0,j)})$. We consider the following morphisms of matrix factorizations for the saddle cobordisms.

**Definition 8.18.** We define the morphisms $\mathcal{F}(f_{s_1}) \in \mathrm{EXT}^1$ and $\mathcal{F}(f_{s_2}) \in \mathrm{EXT}^1$ by

$$\mathcal{F}(f_{s_1}) = \left( \begin{pmatrix} \overline{h}^{n-1}_{123} & -1 \\ -\overline{h}^{n-1}_{134} & -1 \end{pmatrix}, \begin{pmatrix} -1 & 1 \\ \overline{h}^{n-1}_{123} & \overline{h}^{n-1}_{134} \end{pmatrix} \right), \tag{8.19}$$

$$\mathcal{F}(f_{s_2}) = \left( \begin{pmatrix} \overline{h}^{n-1}_{134} & -1 \\ -\overline{h}^{n-1}_{123} & -1 \end{pmatrix}, \begin{pmatrix} -1 & 1 \\ \overline{h}^{n-1}_{134} & \overline{h}^{n-1}_{123} \end{pmatrix} \right), \tag{8.20}$$

where $\overline{h}^{n-1}_{ijk}$ is the $(n-1)$-th complete symmetric function $\overline{h}^{n-1}(x^{(0,i)}, x^{(0,j)}, x^{(0,k)})$.

**Remark 8.19.** The cobordisms $f_{s_1}$ and $f_{s_2}$ have the rotation symmetry by 180 degrees which swaps the two variables $x^{(0,1)}$ and $x^{(0,3)}$ and swaps the two variables $x^{(0,2)}$ and $x^{(0,4)}$. If we suppose the morphism $\overline{\mathcal{F}(f_{s_1})}$, $i = 1, 2$, is symmetric to the morphism $\mathcal{F}(f_{s_i})$ under the rotation of 180 degrees, we have $\overline{\mathcal{F}(f_{s_i})} = -\mathcal{F}(f_{s_i})$. Therefore $\mathcal{F}(f_{s_i})$ can be defined only up to a sign (See also Section 9 [KR08]). However, in our construction, we can regard the rotation as a change of the base point by an operad action. We expect that the sign problem can be solved using the operad structure.

**Problem 8.20.** Can we define a morphism of matrix factorizations of $\beta_{s_0}(\mathrm{id}_T, f_1, ..., f_m)$, where $f_i$ $(i = 1, ..., m)$ is a cobordism of $\mathcal{COB}^{gr}(s_i)$ and $T$ is an element of $\mathrm{P}(\underline{s})$, $\underline{s} = (s_0; s_1, ..., s_m)$?

As we state Remark 8.19, we have the sign problem of morphisms of matrix factorizations corresponding to cobordisms. We have to solve the above sign problem for getting a functor $\mathcal{COB}^{gr}_{G_n/L}$ to $\mathrm{HMF}^{KR_n}$.

We conjecture that there exists an integer $\#(T)$ for $T$ defining the morphism of matrix factorizations corresponding to a cobordism $\beta_{s_0}(\mathrm{id}_T, f_1, ..., f_m)$ by

$$\mathcal{F}(\beta_{s_0}(\mathrm{id}_T, f_1, ..., f_m)) = (-1)^{\#(T)} \mathrm{id}_{\widehat{T}} \otimes_{x_1(T)} \mathcal{F}(e_1(f_1)) \otimes_{x_2(T)} \cdots \otimes_{x_m(T)} \mathcal{F}(e_m(f_m)).$$

Subsequently, we need to show the assignment $\mathcal{F}$ preserves the local relations in $L$ of $\mathcal{COB}^{gr}_{G_n/L}$ in the category $\mathrm{HMF}^{KR_n}$. However, the problem is that we do not know generators of the cobordism isotopy relations.

**Problem 8.21.** Can we find a generating set of all cobordism isotopy of $\mathcal{COB}^{gr}_{G_n/L}$?

If we have generators spanning all isotopy relations by the operad action, it enough to show that the assignment $\mathcal{F}$ preserves generators in $\mathrm{HMF}^{KR_n}$ for getting a functor.

We expect that we solve Problem 8.20 and 8.21 and we construct a functor from $\mathcal{COB}^{gr}_{G_n/L}$ to $\mathrm{HMF}^{KR_n}$. Then, the consistency of $\mathcal{COB}^{gr}_{G_n/L}$ follows from existence of a functor.

# Appendix A. Invariance under Reidemeister moves

Two oriented link diagrams are in the same isotopy class of link diagrams if these are related by a finite sequence of oriented Reidemeister moves. In [Pol10] it was shown that in fact the following moves already generate all oriented Reidemeister moves.

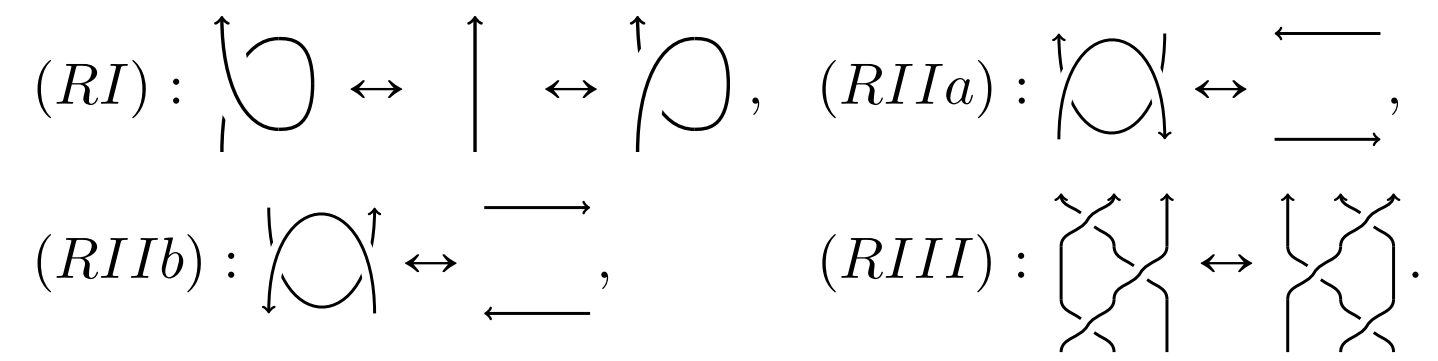

Hence, we show the assignment $[\![-]\!]$ defined in Section 6.2 is invariant under the above list of Reidemeister moves.

## A.1. **Invariance under RIa move.**

**Lemma A.1.** *We have the following isomorphism in the category* $\mathrm{K}^{\mathrm{b}}(\mathrm{Mat}(\mathcal{COB}^{gr}_{G_n/L}))$

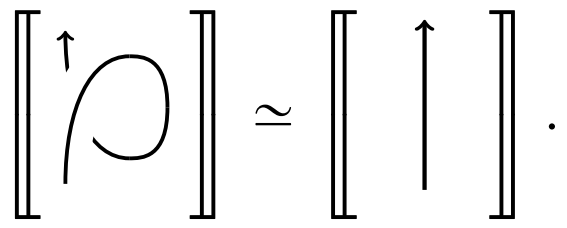

*Proof.* Consider the diagram

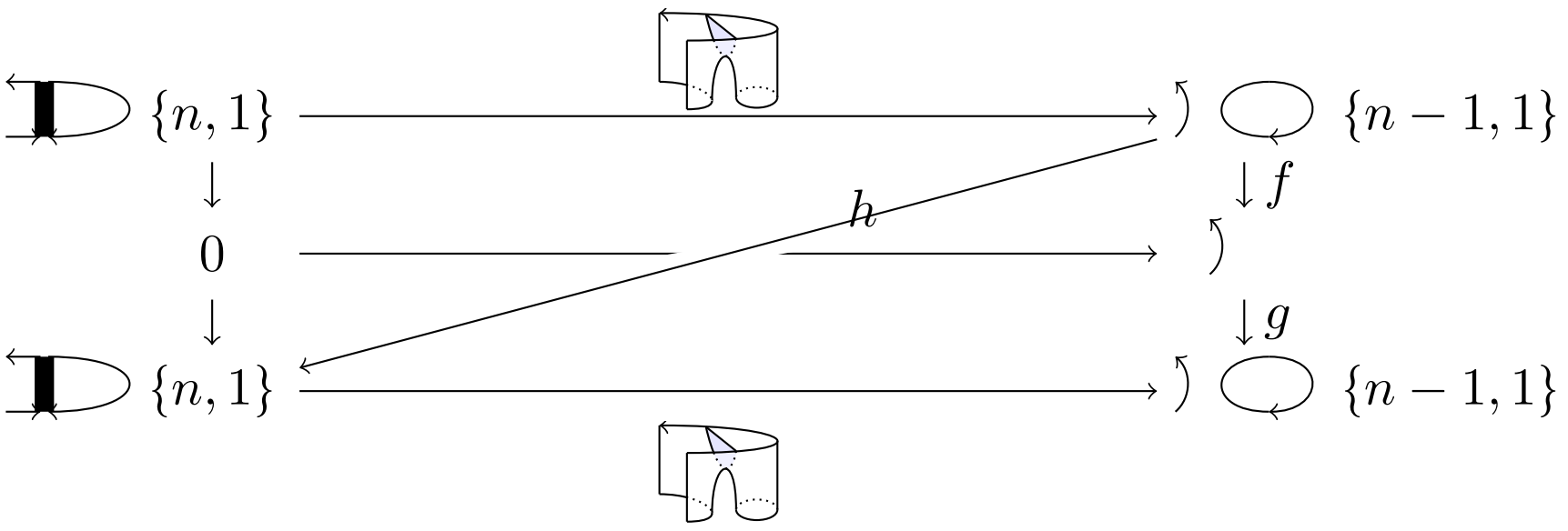

where

$$f = \sum_{\substack{j+k=n-1 \\ j,k\geqslant 0}} (-1)^{nj}$$ 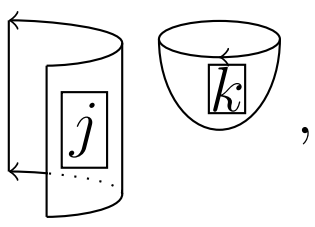 , $$g = (-1)^{n-1} n$$ 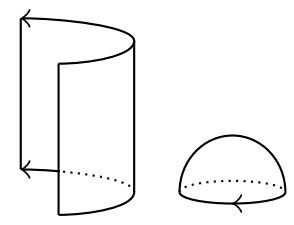 ,

$$h = n \sum_{\substack{j+k+l=n-2 \\ j,k,l\geqslant 0}} (-1)^{n(j+1)-1}$$ 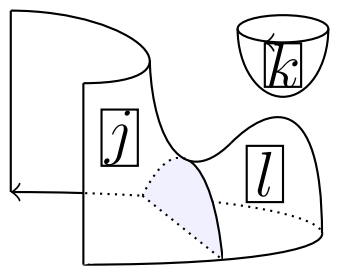 .

- $f$ and $g$ are morphisms of complexes, i.e. $f \circ d = 0$: This is just Lemma 5.17 e).
- $f \circ g \simeq id$: We even have $f \circ g = id$. This follows from (R1) and Lemma 5.17 a).
- $g \circ f \simeq id$:
  - $h \circ d = id$: This is just (R5).
  - $g \circ f + d \circ h = id$: This follows from first applying (R4) and then applying (R3) backwards.

□

A.2. **Invariance under RIb move.**

**Lemma A.2.** *We have the following isomorphism in the category* $\mathrm{K}^{\mathrm{b}}(\mathrm{Mat}(\mathcal{COB}^{gr}_{G_n/L}))$

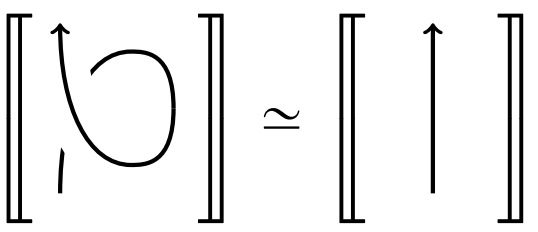

*Proof.* Consider the following diagram:

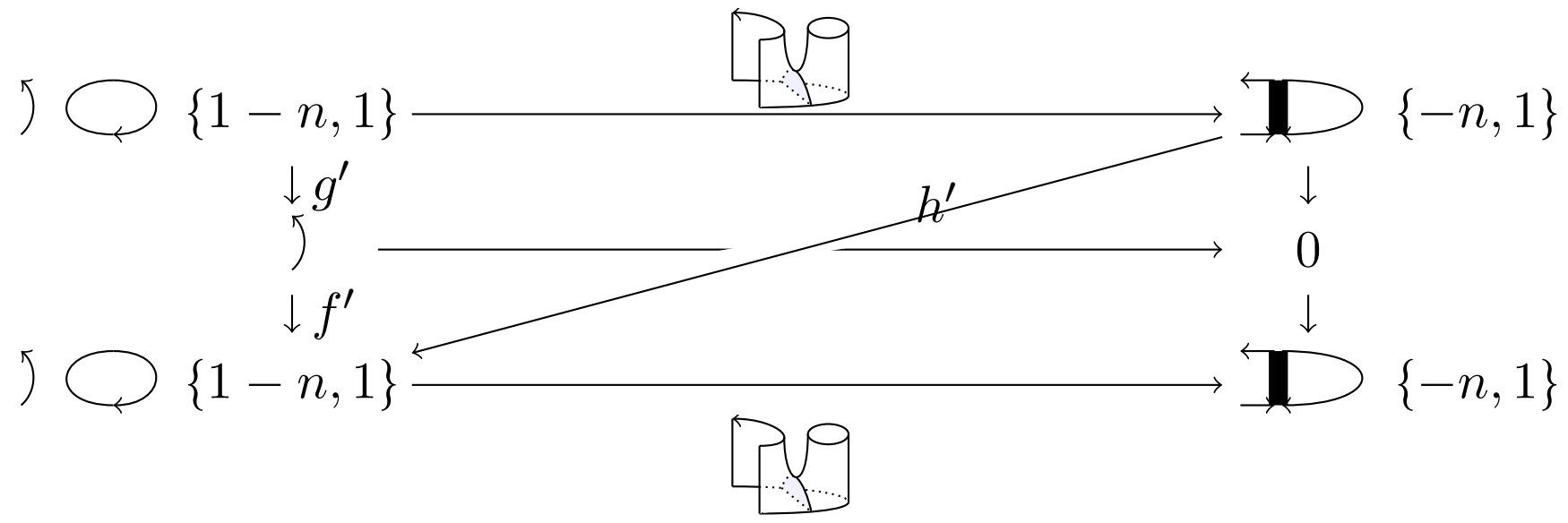

Here $f'$ $g'$ and $h'$ are defined analogously to $f$ and $h$ in the proof of Lemma A.1 and the homotopy equivalence is shown in the same way using Lemma 5.17 f) instead of Lemma 5.17 e). ∎

A.3. **Invariance under RIIa move.**

**Lemma A.3.** *We have the following isomorphism in the category* $\mathrm{K}^{\mathrm{b}}(\mathrm{Mat}(\mathcal{COB}^{gr}_{G_n/L}))$

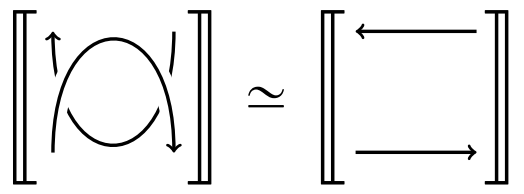

*Proof.* Consider the diagram

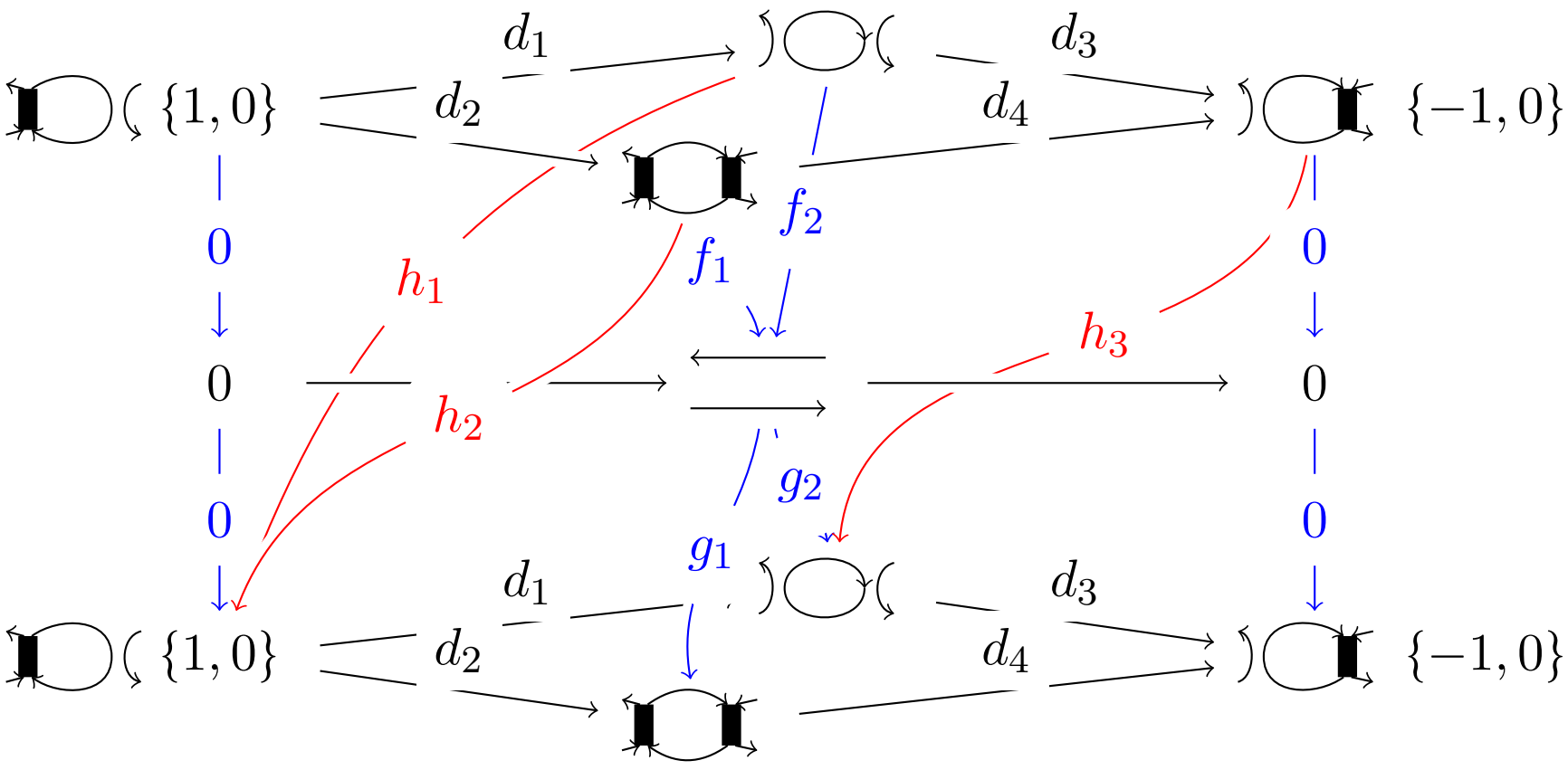

where $d_1 =$ 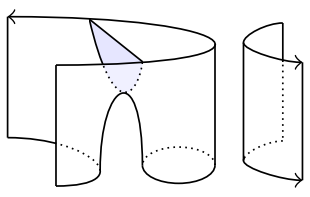 , $d_2 =$ 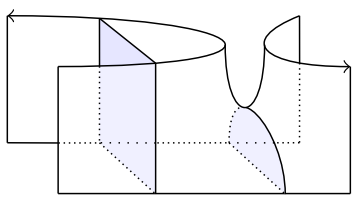 , $d_3 =$ 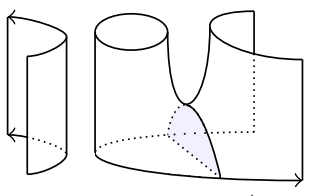 , $d_4 = -$ 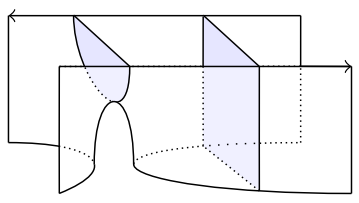 ,

$f_1 =$ 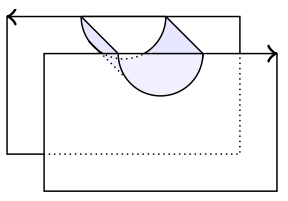 , $f_2 =$ 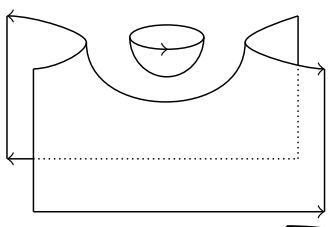 , $g_1 =$ 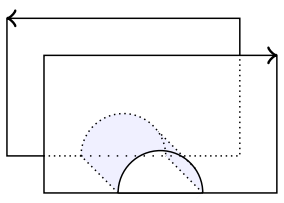 , $g_2 =$ 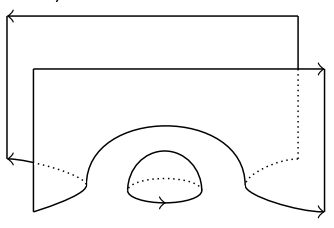 ,

$h_1 = n \sum_{\substack{j+k+l=n-2 \\ j,k,l\geqslant 0}} (-1)^{n(j+l+1)-1}$ 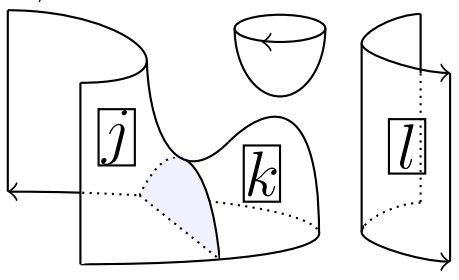 ,

$$h_2 = n\sum_{\substack{i+j+k+l=n-3\\ i,j,k,l\geqslant 0}}(-1)^{n(i+j+1)-1}$$ 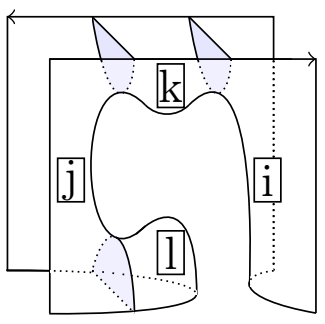 ,

$$h_3 = n\sum_{\substack{j+k+l=n-2\\ j,k,l\geqslant 0}}(-1)^{n(j+1)-1}$$ 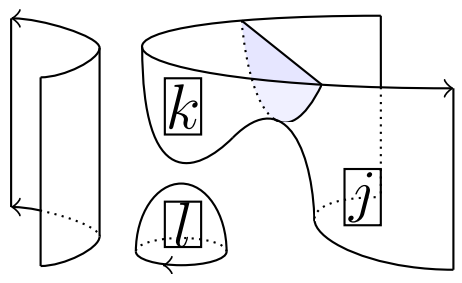 .

- morphism of complexes follows by cobordism isotopy.
- $f_1 \circ g_1 + f_2 \circ g_2 = \mathrm{id}$: follows from (R1) and (R9).
- $g_1 \circ f_1 + d \circ h_1 + h_3 \circ d = \mathrm{id}$: Apply (R3) a) to the first summand and (R4) to the other two summands. After canceling, apply (R3) b) in the vertical direction.
- $g_2 \circ f_1 + d \circ h_1 = 0$: This is exactly Lemma 5.18 a) after using a cobordism isotopy.
- $g_1 \circ f_2 + d \circ h_2 + h_3 \circ d = 0$: Apply (R4) to the second summand. After canceling most of the terms of the second summand, some of the remaining cancel with the third summand and the others cancel with the first summand by Lemma 5.18 b) and a cobordism isotopy.
- $g_2 \circ f_2 + d \circ h_2 = \mathrm{id}$: This follows by a cobordism isotopy and (R8).
- $h_1 \circ d + h_2 \circ d = \mathrm{id}$: Apply (R4) to $h_2 \circ d$, then some summands of $h_2 \circ d$ cancel, also with $h_1 \circ d$, until the right hand side of (R5) with an extra sheet remains.
- $d \circ h_3 = \mathrm{id}$: By putting (R5) a) into 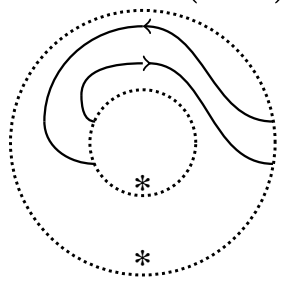 we can rotate it and get exactly the identity we need.

□

## A.4. **Invariance under RIIb move.**

**Lemma A.4.** *We have the following isomorphism in the category* $\mathrm{K}^{\mathrm{b}}(\mathrm{Mat}(\mathcal{COB}^{gr}_{G_n/L}))$

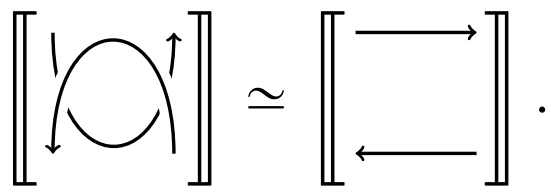

*Proof.* Consider the following diagram:

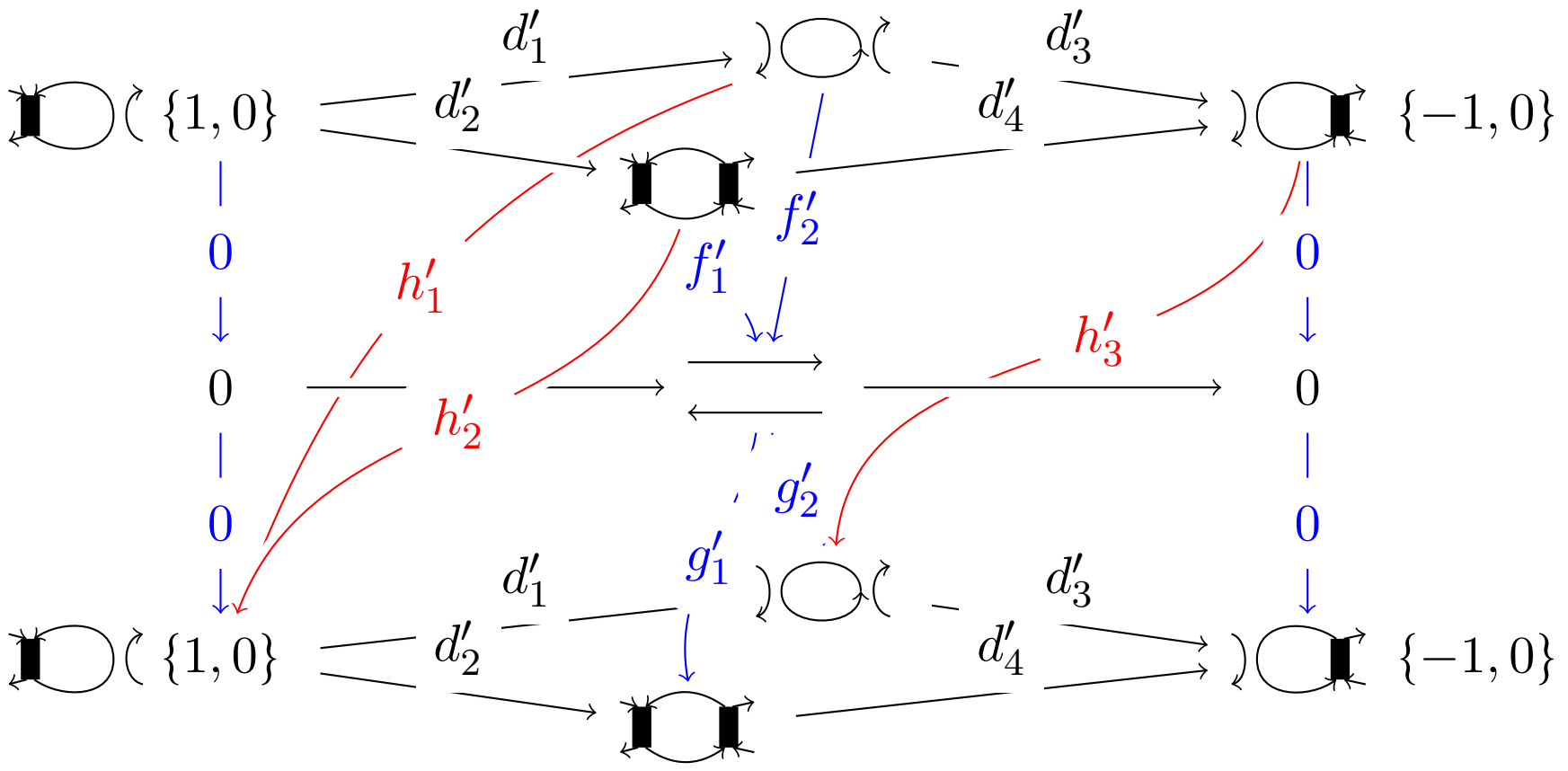

Here, $d'_1, \ldots, d'_4$, $f'_1, f'_2$, $g'_1, g'_2$ and $h'_2$ are exactly defined as $d_1, \ldots, d_4$, $f_1, f_2$, $g_1, g_2$ resp. $h_2$ just with the other orientation while $h'_1$ and $h'_3$ also differ from $h_1$ resp. $h_3$ by a sign:

$$h'_1 = n \sum_{\substack{j+k+l=n-2 \\ j,k,l \geqslant 0}} (-1)^{n(i+j)}$$ 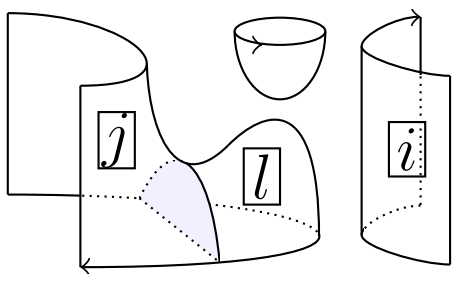 ,

$$h'_2 = n \sum_{\substack{i+j+k+l=n-3 \\ i,j,k,l \geqslant 0}} (-1)^{n(i+j+1)-1}$$ 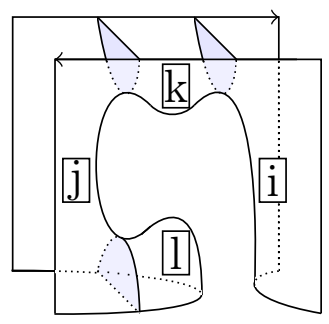 ,

$$h'_3 = n \sum_{\substack{j+k+l=n-2 \\ j,k,l \geqslant 0}} (-1)^{ni}$$ 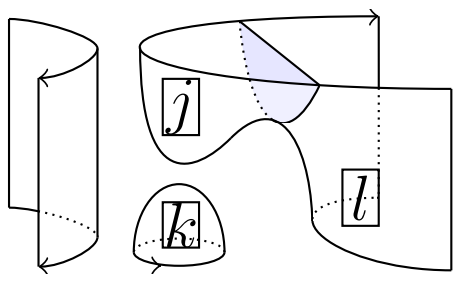 .

The proof is totally analogous to the one of Lemma A.5 using Lemma 5.18. We use Lemma 5.18 c), d) instead of a), b). □

A.5. **Invariance under RIII move.**

**Lemma A.5.** *We have the following isomorphism in the category* $\mathrm{K}^{\mathrm{b}}(\mathrm{Mat}(\mathcal{COB}^{gr}_{G_n/L}))$

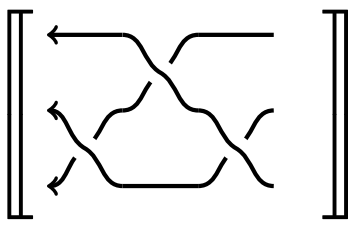 $\simeq$ 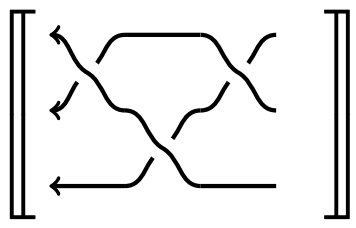 .

*Proof.* Consider the diagram

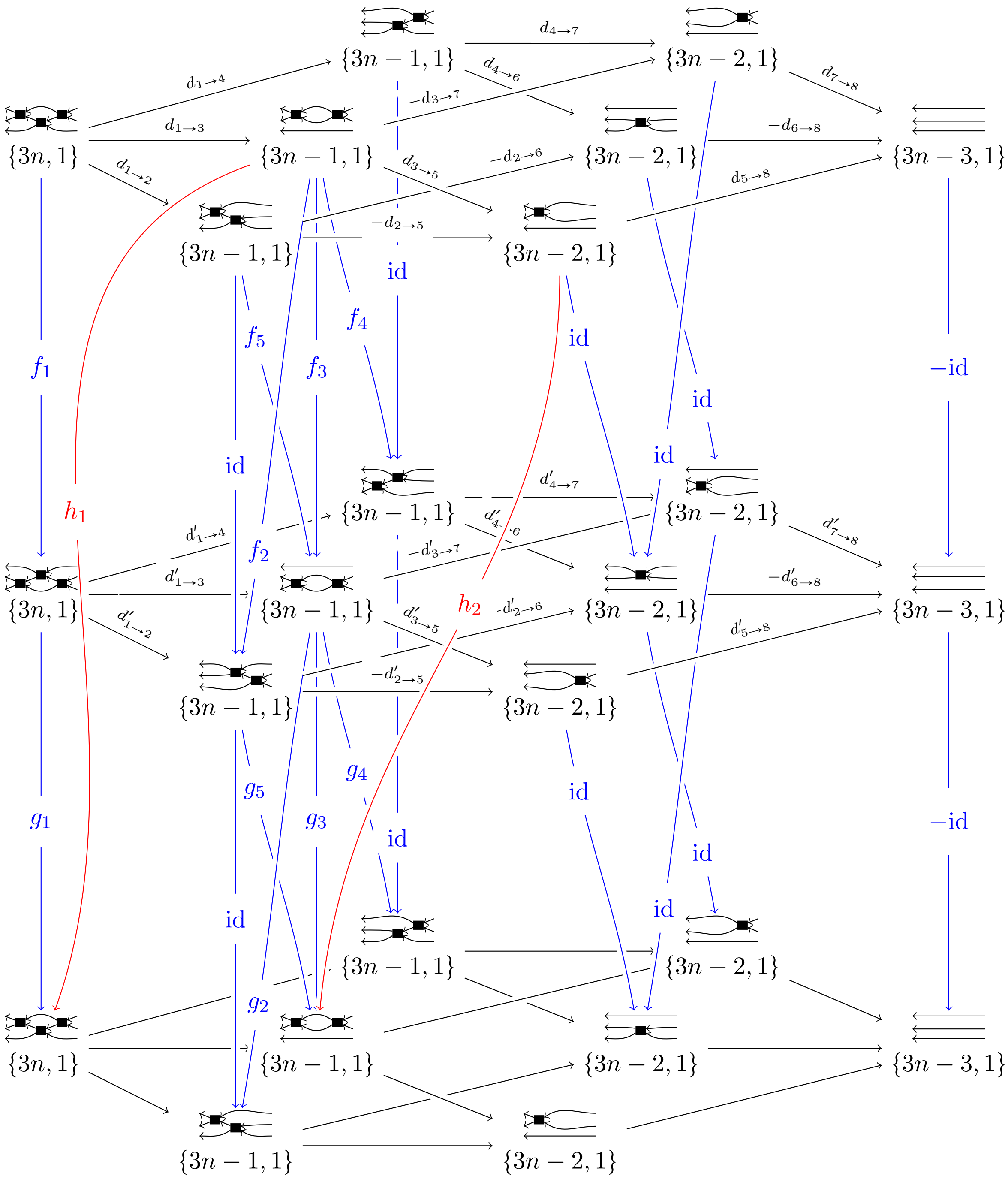

where

$f_1 = \Bigg\{$ 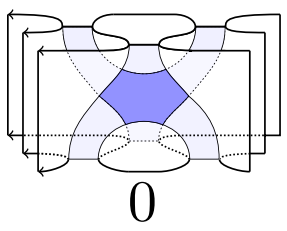 if $n \geqslant 3$; 0 if $n = 2$, $f_2 =$ 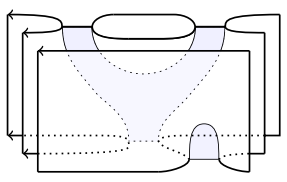, $f_3 =$ 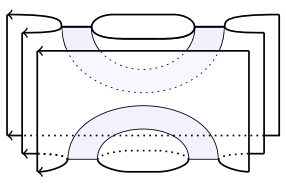, $f_4 =$ 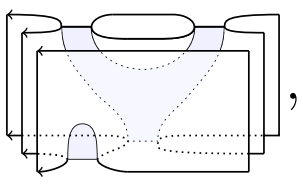,

$f_5 =$ 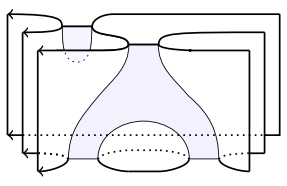, $g_1 = \Bigg\{$ 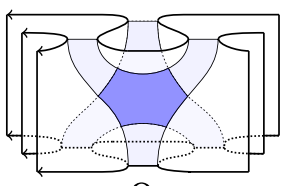 if $n \geqslant 3$; 0 if $n = 2$, $g_2 =$ 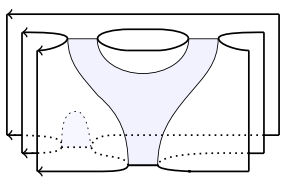, $g_3 =$ 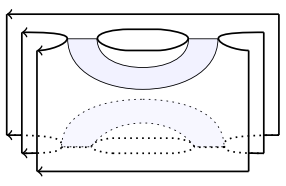,

$g_4 =$ 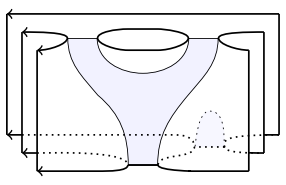, $g_5 =$ 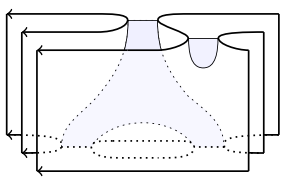, $h_1 = -$ 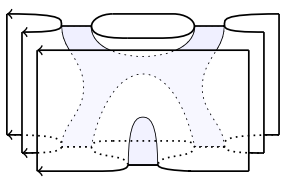, $h_2 =$ 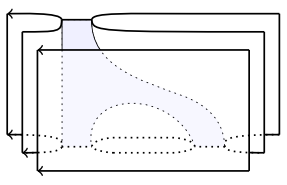,

and the differentials $d_{i\to j}$ and $d'_{i\to j}$ are given by applying 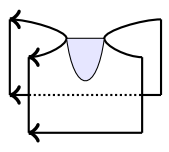 at the appropriate places (with the signs as indicated).

- morphism of complexes: Consider the upper morphisms first. Here the notation $(a \to b)$ means that we consider (commutative) squares starting at $a$ and ending at $b$, where the labeling is from left to right.
  - $(1 \to 2)$: $d'_{1\to 2}f_1 = \mathrm{id}d_{1\to 4} + f_2 d_{1\to 3}$ follows from Lemma 5.17 g) and a cobordism isotopy.
  - $(1 \to 3)$: $d'_{1\to 3}f_1 = f_5 d_{1\to 4} + f_3 d_{1\to 3}$ follows from (R13)b) and a cobordism isotopy
  - $(1 \to 4)$: $d'_{1\to 4}f_1 = f_4 d_{1\to 3} + \mathrm{id}d_{1\to 2}$ follows from Lemma 5.17 h) and a cobordism isotopy
  - $(3 \to 6)$: $-d'_{2\to 6}f_2 + d'_{4\to 6}f_4 = \mathrm{id}d_{3\to 5} + \mathrm{id}(-d_{3\to 7})$ follows from Lemma 5.20 a).
  - The rest of commutative diagrams composed of morphisms can be obtained by cobordism isotopies.

  The lower morphisms are a morphism of complexes analogously, using (R13)a) instead of (R13)b), and a reflected version of Lemma 5.20 a), which is proved exactly the same way.
- homotopy $g_\bullet f_\bullet \overset{h}{\sim} \mathrm{id}_\bullet$:
  - $(1 \to 1)$: $\mathrm{id}_{1\to 1} - g_1 f_1 = h_1 d_{1\to 3}$ follows from (R12)b).
  - $(2 \to 2)$: $\mathrm{id}_{2\to 2} - \mathrm{idid} - g_2 f_5 = 0$ follows from (R7)
  - $(3 \to 2)$: $0 - (-\mathrm{id})f_2 - g_2 f_3 = d'_{1\to 2}h_1$ follows from (R7) and a cobordism isotopy
  - $(3 \to 3)$: $id_{3\to 3} - g_5 f_2 - g_3 f_3 = d'_{1\to 3}h_1 + h_2 d_{3\to 5}$ follows from Lemma 5.20 b).
  - $(3 \to 4)$: $0 - \mathrm{id}f_4 - g_4 f_3 = d_{1\to 4}h_1$ follows from (R7) and a cobordism isotopy
  - The rest of homotopic relations can be obtained by cobordism isotopies.

  The homotopy for the other composition is given by the maps

  $h'_1 = -$ 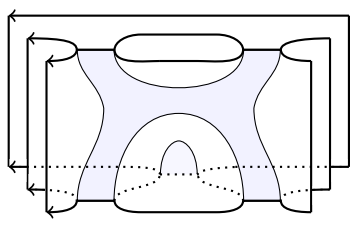, $h'_2 =$ 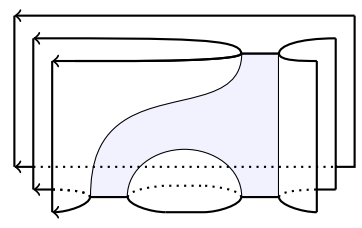

  These cobordisms are the analogs to $h_1$ and $h_2$ and the homotopy is shown analogously by using the same relations and calculations as above but (R12)a) instead of (R12)b) and a reflected version of Lemma 5.20 b), which is also proved exactly the same way.

□

*E-mail address*: `gischaef@math.uni-bonn.de`

Department of Mathematics, University of Bonn, Bonn, Germany

*E-mail address*: yasuyoshi.yonezawa@math.nagoya-u.ac.jp

Graduate school of Mathematics, Nagoya University, Nagoya, Japan